\documentclass[twoside,11pt]{article}
\usepackage{JMLR}
\usepackage{titletoc}

\renewenvironment{proof}[1][\relax]{
\par
\ifx
#1\relax \noindent{\bf Proof\ }
\else
\noindent{\bf #1\ }
\fi
}
{
\hfill\BlackBox\\[2mm]
}

\usepackage{times}
\usepackage{amsfonts}
\usepackage{amsmath, nccmath}
\usepackage{geometry}

\usepackage{comment}
\usepackage{xcolor}
\usepackage{enumitem}
\usepackage{color}

\usepackage{microtype}
\usepackage{graphicx}

\usepackage{bbold}
\usepackage{xfrac}

\usepackage{thmtools}
\usepackage{thm-restate}
\usepackage{cleveref}

\crefname{appsec}{appendix}{appendices}

\newtheorem{assumption}{Assumption}

\usepackage[ruled]{algorithm2e}
\usepackage{algorithmic}

\usepackage{booktabs}

\usepackage{mathtools}
\usepackage{titletoc}
\usepackage{hyperref}
\hypersetup{
    colorlinks = true,
    citecolor = blue,
    linkcolor = blue
}
\newcommand{\floor}[1]{\left\lfloor #1 \right\rfloor}
\newcommand{\bh}{\boldsymbol{h}}
\newcommand{\btau}{\boldsymbol{\tau}}

\newcommand{\be}{\boldsymbol{e}}
\newcommand{\bomega}{\boldsymbol{\omega}}
\newcommand{\bx}{\boldsymbol{x}}
\newcommand{\bw}{\boldsymbol{w}}
\renewcommand{\d}{\,\mathrm{d}}
\newcommand{\bb}{\boldsymbol{b}}
\newcommand{\ba}{\boldsymbol{a}}
\newcommand{\bz}{\boldsymbol{z}}

\newcommand{\bT}{\boldsymbol{T}}

\newcommand{\bm}{\boldsymbol{m}}
\newcommand{\bc}{\boldsymbol{c}}

\newcommand{\bW}{\boldsymbol{W}}
\newcommand{\bv}{\boldsymbol{v}}
\newcommand{\bU}{\boldsymbol{U}}

\newcommand{\bu}{\boldsymbol{u}}
\newcommand{\bg}{\boldsymbol{g}}
\newcommand{\bzeta}{\boldsymbol{\zeta}}

\DeclareMathOperator{\Var}{Var}

\newcommand{\bfI}{\mathbf{I}}
\newcommand{\bfA}{\mathbf{A}}
\newcommand{\tC}{\texttt{C}}

\newcommand{\bfQ}{\mathbf{Q}}
\newcommand{\bq}{\boldsymbol{q}}

\newcommand{\ball}{B}
\newcommand{\sphere}{\partial B}

\newcommand{\Prob}{\mathbf{P}}
\newcommand{\bL}{\bar{L}}

\newcommand{\Exp}{\mathbf{E}}

\newcommand{\bbR}{\mathbb{R}}
\newcommand{\bbN}{\mathbb{N}}
\newcommand{\com}{\Theta}

\newcommand{\class}[1]{\mathcal{#1}}

\newcommand{\parent}[1]{\left( #1 \right)}

\newcommand{\norm}[1]{\left\lVert#1\right\rVert}
\newcommand{\abs}[1]{\left\lvert #1\right\rvert}

\newcommand{\enscond}[2]{\left\{ #1 \, : \, #2\right\}}

\newcommand{\cst}{\texttt{A}}
\newcommand{\scalar}[2]{\left\langle #1\,, #2 \right\rangle}
\DeclareMathOperator*{\argmin}{arg\,min}
\newcommand{\proj}{\ensuremath{\text{\rm Proj}}}

\newcommand{\eg}{{\em e.g.,}}

\newcommand{\resp}{{\em resp.~}}

\newcommand{\wrt}{{\em w.r.t.~}}
\newcommand{\iid}{{\rm i.i.d.~}}

\renewcommand{\enspace}{\,}
\DeclareMathOperator*{\sign}{sign}

\usepackage[textwidth=2.0cm,textsize=tiny]{todonotes}
\jmlrheading{1}{2000}{1-48}{4/00}{10/00}{paper00}{jmlr related}

\ShortHeadings{Gradient-free
optimization of highly smooth functions}{Akhavan, Chzhen, Pontil, Tsybakov}
\firstpageno{1}

\definecolor{purple}{rgb}{0., 0., 0.}
\definecolor{brown}{rgb}{0., 0., 0.}

\begin{document}

% \title{Zero Order Optimization for Convex and non-Convex Higher Order Smooth Functions}

% Evg: Title suggestions here
\title{%Zero-order
Gradient-free
optimization of highly smooth functions: improved analysis and a new algorithm}

\author{\name Arya Akhavan \email aria.akhavanfoomani@iit.it \\
      \addr CSML, Istituto Italiano di Tecnologia\\
      CMAP, École Polytechnique, IP Paris
      \AND
      \name Evgenii Chzhen \email evgenii.chzhen@cnrs.fr \\
      \addr
      CNRS, LMO, 
      Université Paris-Saclay
      \AND
      \name Massimiliano Pontil \email massimiliano.pontil@iit.it \\
      \addr 
    CSML, Istituto Italiano di Tecnologia\\
      University College London
      \AND
      \name Alexandre B. Tsybakov \email alexandre.tsybakov@ensae.fr \\
      \addr CREST, ENSAE, IP Paris
      }

\editor{
% to be determined
}
\maketitle
% \doparttoc % Tell to minitoc to generate a toc for the parts
% \faketableofcontents % Run a fake tableofcontents command for the partocs

% \part{} % Start the document part
% \parttoc % Insert the document TOC
% \startcontents[mainsections]

\begin{abstract}
%\evg{Arya, do not forget to change affiliation.}
This work studies minimization problems with zero-order noisy oracle information under the assumption that the objective function is highly smooth and possibly satisfies additional properties. We consider two kinds of zero-order projected gradient descent algorithms, which differ in the form of the gradient estimator. The first algorithm uses a gradient estimator based on randomization over the $\ell_2$ sphere due to \cite{BP2016}.  We present an improved analysis of this algorithm on the class of highly smooth and strongly convex functions studied in the prior work, and we derive rates of convergence for two more general classes of non-convex functions. Namely, we consider highly smooth functions satisfying the Polyak-Łojasiewicz condition and the class of highly smooth functions with no additional property. The second algorithm is based on randomization over the $\ell_1$ sphere, and it extends to the highly smooth setting the algorithm that was recently proposed for Lipschitz convex functions in \cite{Akhavan_Chzhen_Pontil_Tsybakov22a}. We show that, in the case of noiseless oracle, this novel algorithm enjoys better bounds on bias and variance than the $\ell_2$ randomization and the commonly used Gaussian randomization algorithms, while in the noisy case both $\ell_1$  and $\ell_2$  algorithms benefit from similar improved theoretical guarantees.
The improvements are achieved thanks to a new proof techniques based on Poincaré type inequalities for uniform distributions on the $\ell_1$ or $\ell_2$ spheres. The results are established under weak (almost adversarial) assumptions on the noise. Moreover, we provide minimax lower bounds proving optimality or near optimality of the obtained upper bounds in several cases. 
\end{abstract}

\section{Introduction}
In this work, we study the problem of gradient-free optimization for certain types of smooth functions. Let $f: \bbR^d \to \bbR$ and $\com \subseteq \bbR^d$. We are interested in solving the following optimization problem
\begin{align*}
    f^\star \triangleq \inf_{\bx \in \com} f(\bx)\enspace,
\end{align*}
and we assume that $f^\star$ is finite. One main theme of this paper is to exploit higher order smoothness properties of the underlying function $f$ in order to improve the performance of the optimization algorithm. We consider that the algorithm has access to a zero-order stochastic oracle, which, given a point $\bx \in \bbR^d$ returns a noisy value of $f(\bx)$, under a general noise model. 

We study two kinds of zero-order projected gradient descent algorithms, which differ in the form of the gradient estimator.
Both algorithms can be written as an iterative update of the form
\begin{align*}
    \bx_1 \in \com\qquad\text{and}\qquad \bx_{t + 1} = \proj_\com(\bx_t - \eta_t \bg_t)\qquad t=1, 2,\dots \enspace,
\end{align*}
where $\bg_t$ is a gradient estimator at the point $\bx_t$, $\eta_t>0$ is a step size, and $\proj_\com(\cdot)$ is the Euclidean projection operator onto the set $\com$.
In either case, the gradient estimator is built from two noisy function values, that are queried at two random perturbations of the current guess for the solution, and it involves an additional randomization step. Also, both algorithms invoke smoothing kernels in order to take advantage of higher order smoothness, following the approach initially suggested in \citep{PT90}. The first algorithm uses a form of $\ell_2$ randomization introduced by \cite{BP2016} and it has been studied in the context of gradient-free optimization of strongly convex functions in \citep{akhavan2020,Gasnikov}. The second algorithm is an extension, by incorporating smoothing kernels, of the approach proposed and analysed in~\cite{Akhavan_Chzhen_Pontil_Tsybakov22a} for online minimization of Lipschitz convex functions. It is based on an alternative randomization scheme, which uses $\ell_1$-geometry in place of the $\ell_2$ one.

A principal goal of this paper is to derive upper bounds on the expected optimization error of both algorithms under different assumptions on the underlying function $f$. These assumptions are used to set the step size in the algorithms and the perturbation parameter used in the gradient estimator. Previous works on gradient-free optimization of highly smooth functions considered mostly the strongly convex case~\citep{PT90,BP2016,akhavan2020,akhavan2021distributed,Gasnikov}.  
In this paper, we provide a refined analysis of the strongly convex case, improving the dependence on the dimension $d$ and the strong convexity parameter $\alpha$ for the algorithm with $\ell_2$ randomization, and showing analogous results for the new method with $\ell_1$ randomization. 
%In the basic particular case of 
For the special case of strongly convex functions with Lipschitz gradient, we find the minimax optimal dependence of the bounds on all the three parameters of the problem 
%-- the horizon $T$, the dimension $d$ and $\alpha$, -- 
(namely, the horizon $T$, the dimension $d$ and $\alpha$) and we show that both algorithms attain the minimax rate, which equals $\alpha^{-1}d/\sqrt{T}$. This finalizes the line of work starting from \citep{PT90}, where it was proved that optimal dependence on $T$ is of the order $1/\sqrt{T}$, and papers proving that optimal dependence on $d$ and $T$ scales as $d/\sqrt{T}$ (\cite{Shamir13} establishing a lower bound and \cite{akhavan2020} giving a complete proof, see the discussion in Section \ref{sec:LB}). Furthermore, we complement these results by considering highly smooth but not necessary convex functions $f$, and highly smooth functions $f$, which additionally satisfy the gradient dominance (Polyak-Łojasiewicz) condition. To this end, we develop unified tools that can cover a variety of gradient estimators, and then apply them to the algorithm with $\ell_2$ randomization and to our new algorithm based on $\ell_1$ randomization. We show that, in the case of noiseless oracle, this novel algorithm enjoys better bounds on bias and variance than its $\ell_2$ counterpart, while in the noisy case, both algorithms benefit from similar theoretical guarantees.
The improvements in the analysis are achieved thanks to a new method of evaluating the bias and variance of both algorithms based on Poincaré type inequalities for uniform distributions on the $\ell_1$ or $\ell_2$ spheres. Moreover, we establish all our upper bounds under very weak (almost adversarial) assumptions on the noise. 
%Finally, we provide minimax lower bounds proving optimality or near optimality of the obtained upper bounds in several cases. \massi{should we remove this last sentence as we already spoke of lower bounds}

\subsection{Summary of the upper bounds}
%\massi[inline]{this subsection is good, what it misses is a little discussion/intuition about how/why the different assumptions lead to different rates. Also briefly mention the lower bound}

%To give a high-level overview of the obtained results, in this section, we state, in an informal way, the main  contributions of this work. 

In this subsection, we give a high-level overview of the main contributions of this work. Apart from the improved guarantees for the previously studied function classes, %of objective function $f$, 
one of the main novelties of our work is the analysis in the case of a non-convex smoothness objective function $f$, for which we provide a convergence rate to a stationary point. %objective function $f$. In particular, we provide a convergence rate to a stationary point of $f$, assuming only its smoothness. 
Furthermore, we study the case of $\alpha$-gradient dominant $f$, a popular relaxation of strong convexity, which includes non-convex functions. To the best of our knowledge, the analysis of noisy zero-order optimization in these two cases is novel.
In Section~\ref{sec:LB} we derive minimax lower bounds and discuss the optimality or near optimality of our convergence rates. 

%In each of the following paragraphs
In the following we highlight the guarantees that we derive for the two analysed algorithms.
Each of the guarantees differs in the dependency on the main parameters of the problem, which is a consequence of the different types of available properties of the objective function. Let us also mention that we mainly deal with the unconstrained optimization case, $\com = \bbR^d$. This is largely due to the fact that the Polyak-Łojasiewicz inequality is mainly used in the unconstrained case and the exertion of this condition to the constrained case is still an active area of research~\cite[see e.g.,][and references therein]{balashov2020gradient}. Meanwhile, for the strongly convex case, as in previous works \citep{BP2016,akhavan2020,Gasnikov}, we additionally treat the constrained optimization. In this section, we only sketch our results in the case $\com = \bbR^d$.

We would like to emphasize that we do not assume the measurement noise to have zero mean. Moreover, the noise can be non-random and no independence between noises on different time steps is required, so that the setting can be considered as \emph{almost} adversarial.

\paragraph{Rate of convergence under only smoothness assumption.}Assume that $f$ is a $\beta$-H\"older function with Lipschitz continuous gradient, where $\beta\ge2$. Then, after $2T$ oracle queries both considered algorithms provide a point $\bx_S$  satisfying
        \begin{align*}
            \Exp\left[\norm{\nabla f(\bx_S)}^2\right] \lesssim \parent{\frac{d^2}{T}}^{\frac{\beta - 1}{2\beta - 1}} \text{ under the assumption that } T \geq d^{\frac{1}{\beta}}\enspace,
        \end{align*}
        where $S$ is a random variable with values in $\{1,\dots,T\}$, $\norm{\cdot}$ denotes the Euclidean norm, and the sign $\lesssim$ conceals a multiplicative constant that does not depend on $T$ and $d$. 
         To the best of our knowledge, this result is the first convergence guarantee for the zero-order stochastic optimization under the considered noise model. In a related development,
\cite{Ghadimi2013,balasubramanian2021zeroth} study zero-order optimization of non-convex objective function with Lipschitz gradient, which corresponds to $\beta=2$. They assume querying two function values with identical noises, in which case the analysis and the convergence rates are essentially analogous to the particular case of our setting with no noise (see the discussion in Section \ref{sec:classical stoch opt} below). The work of \cite{duchi} studies noiseless optimization of highly smooth functions assuming that the derivatives up to higher order are observed, and \cite{azjev} consider stochastic optimization with first order oracle. These papers cannot be directly compared with our work as the settings are different. 

\paragraph{Rate of convergence under smoothness and Polyak-Łojasiewicz assumptions.} {Assume that $f$ is a $\beta$-H\"older function with Lipschitz continuous gradient, with a Lipschitz constant $\bL$. Additionally, let $\beta \geq 2$, and suppose that $f$ satisfies the Polyak-Łojasiewicz inequality with a constant $\alpha$.} Then, after $2T$ oracle queries, both considered algorithms provide a point $\bx_T$, for which the expected optimization error satisfies  
        \begin{align*}
        \Exp[f(\bx_T)-f^\star] \lesssim %\frac{d }{\alpha T}+
        \frac{1}{\alpha }\left(\frac{{\color{purple}\mu}d^2}{T}\right)^{\frac{\beta-1}{\beta}} \text{ under the assumption that } T \gtrsim d^{2 - \frac{\beta}{2}}\enspace,
        \end{align*}
        where {\color{purple} $\mu = \bL/\alpha$}, and the signs $\lesssim$ and $\gtrsim$ conceal multiplicative constants that do not depend on $T$, $d$ and $\alpha$. The Polyak-Łojasiewicz assumption was introduced in the context of first order optimization by~\cite{PB63} who showed that it implies linear convergence of the gradient descent algorithm. Years later, this condition received attention in the machine learning and optimization community following the work of~\cite{10.1007/978-3-319-46128-1_50}. To the best of our knowledge, zero-order optimization under the considered noise model with the Polyak-Łojasiewicz assumption was not previously studied. Very recently~\cite{Rando_Molinari_Villa_Rosasco22} studied a related problem under the Polyak-Łojasiewicz assumption when querying two function values with identical noises, which can be compared with the analysis in the particular case of our setting with no noise (see the discussion in Section \ref{sec:classical stoch opt}). Unlike our work, ~\cite{Rando_Molinari_Villa_Rosasco22} %deploy a per-coordinate gradient estimator similar to \citep{akhavan2021distributed} allowing queries only in orthogonal directions, 
        do not deal with higher order smoothness and do not derive the dependency of the bounds on $d$, $\mu$ and $\alpha$.

\paragraph{Rate of convergence under smoothness and strong convexity.} 
Assume that $f$ is a $\beta$-H\"older function with Lipschitz continuous gradient, where $\beta\ge2$, and satisfies $\alpha$-strong convexity condition. Then, after $2T$ oracle queries, both considered algorithms provide a point $\hat\bx_T$ such that
        \begin{align*}
        \Exp[f(\hat\bx_T)-f^\star] \lesssim %\frac{d }{\alpha T}+
        \frac{1}{\alpha}\left(\frac{d^2}{T}\right)^{\frac{\beta-1}{\beta}} \text{ under the assumption that } T \geq d^{2 - \frac{\beta}{2}}\enspace,
        \end{align*}
        where $\lesssim$ conceals a multiplicative constant that does not depend on $T$, $d$ and $\alpha$.
    % \end{itemize}
% \end{theorem}
The closest result to ours is obtained in~\cite{akhavan2020} and it splits into two cases: $\beta = 2$ (Lipschitz continuous gradient) and $\beta > 2$ (higher order smoothness). For $\beta = 2$, \cite{akhavan2020} deal with a compact $\Theta$ and prove a bound with optimal dependence on the dimension (linear in $d$) but sub-optimal in $\alpha$, while for $\beta > 2$ they derive (for $\Theta=\mathbb{R}^d$ and for compact $\Theta$) the rate with sub-optimal dimension factor $d^2$. Later,~\cite{akhavan2021distributed} and~\cite{Gasnikov} improved the dimension factor to $d^{2 - \sfrac{1}{\beta}}$ for $\beta > 2$, which still does not match the linear dependence as $\beta \to 2$.
% By considering a slightly different definition of smoothness, originally used in this context by~\cite{BP2016},
In contrast, by considering a slightly different definition of smoothness, we provide below a unified analysis leading to the dimension factor $d^{2-\sfrac{2}{\beta}}$ for any $\beta \geq 2$, under constrained and unconstrained $\Theta$; the improvement is both in the rate and in the proof technique.

\subsection{Notation} Throughout the paper, we use the following notation. For any $k \in \bbN$ we denote by $[k]$, the set of first $k$ positive integers. For any $\bx \in \bbR^d$ we denote by $\bx \mapsto \sign(\bx)$ the component-wise sign function (defined at $0$ as $1$). We let $\langle \cdot, \cdot \rangle$ and $\|\cdot\|$ be the standard inner product and Euclidean norm on $\mathbb{R}^d$, respectively. %We denote by $\norm{\cdot}$ the Euclidean norm in $\bbR^d$ and by $\norm{\cdot}_p$ the $\ell_p$-norm in $\bbR^d$.
For every close convex set $\com\subset \mathbb{R}^d$ and $\bx\in \mathbb{R}^d$ we denote by $\proj_\com(\bx) = {\rm argmin} \{ \|\bz-\bx\| \,:\, \bz \in \com\}$ the Euclidean projection of $\bx$ onto $\com$. 
For any $p \in [1, +\infty]$ we let $\norm{\cdot}_p$ be the $\ell_p$-norm in $\bbR^d$ and we introduce the open $\ell_p$-ball 
%$\ball^d_p$ 
and $\ell_p$-sphere, respectively, as
%$\sphere^d_p$ as
\begin{align*}
    \ball^d_p \triangleq \enscond{\bx \in \bbR^d}{\norm{\bx}_p < 1}\qquad\text{and}\qquad\sphere^d_p \triangleq \enscond{\bx \in \bbR^d}{\norm{\bx}_p = 1}\enspace.
\end{align*}
For any $\beta \geq 2$ we let $\lfloor \beta\rfloor$ be the largest integer strictly less than $\beta$. Given a multi-index  $\bm=(m_{1},\ldots,m_{d}) \in \mathbb{N}^d$, we set $\bm!\triangleq m_{1}! \cdots m_{d}!$, $|\bm| \triangleq m_{1}+ \cdots+m_{d}$. 

\subsection{Structure of the paper} The paper is organized as follows. 
%in the following manner. 
In Section~\ref{sec2}, we recall some preliminaries and introduce the classes of functions considered throughout.~In Section~\ref{sec3}, we present the two algorithms that are studied in the paper.
In Section \ref{sec4}, we present the upper bounds for both algorithms, and in each of the considered function classes.
In Section~\ref{sec:LB}, we establish minimax lower bounds for the zero-order optimization problem. The proofs of most of the results are presented in the appendix.

\section{Preliminaries}
\label{sec2}
For any multi-index $\bm \in \bbN^d$, any $|\bm|$ times continuously differentiable function $f : \bbR^d \to \bbR$, and every $\bh=(h_1,\dots, h_d)^\top \in \mathbb{R}^{d}$ we define 
\[
D^{\bm}f(\bx) \triangleq \frac{\partial ^{|\bm|}f(\bx)}{\partial ^{m_{1}}x_{1} \cdots\partial ^{m_{d}}x_{d}}
\,,
\qquad
\bh^{\bm} \triangleq h_{1}^{m_{1}} \cdots h_{d}^{m_{d}}\enspace. 
\]
For any $k$-linear form $A: \parent{\bbR^d}^k \to \bbR$, we define its norm as
\begin{align*}
    \|A\| \triangleq \sup\enscond{\abs{A[\bh_1, \ldots, \bh_k]}}{\|\bh_j\| \leq 1,\,\,j \in [k]}\enspace.
\end{align*}
Whenever $\bh_1 = \ldots = \bh_k = \bh$ we write $A[\bh]^k$ to denote $A[\bh, \ldots, \bh]$.
Given a $k$ times continuously differentiable function $f : 
\bbR^d \to \bbR$ and $\bx \in \bbR^d$ we denote by $f^{(k)}(\bx) : \parent{\bbR^d}^k \to \bbR$ the following $k$-linear form:
% \massi[inline]{the sum qualification looks a bit strange/complicated, maybe it can be simplified}
\begin{align*}
    f^{(k)}(\bx)[\bh_1, \ldots, \bh_k]
    &=
    \sum_{|\bm_1| = \cdots = |\bm_k| = 1} D^{\bm_1 + \cdots + \bm_k}f(\bx) \bh_1^{\bm_1}  \cdots  \bh_k^{\bm_k}\,,\quad\, \forall \bh_1, \ldots, \bh_k \in \bbR^d\enspace,
\end{align*}
% {
%\evg{This is equivalent for $k$-times cont. diff. Not sure that simpler though.}
%\begin{align*}    f^{(k)}(\bx)[\bh_1, \ldots, \bh_k] = \frac{\partial^k }{\partial t_1\ldots\partial t_k}\bigg|_{t_1 = \ldots = t_k = 0}f(\bx + t_1 \bh_1 + \ldots + t_k \bh_k)\,,\quad\, \forall \bh_1, \ldots, \bh_k \in \bbR^d\enspace.\end{align*}
% }
% \begin{align*}
%     f^{(k)}(\bx)[\bh_1, \ldots, \bh_k]
%     &=
%     \sum_{|\bm| = k} D^{\bm}f(\bx) \bh_1^{\bm_1}  \cdots  \bh_k^{\bm_k}\,,\quad\, \forall \bh_1, \ldots, \bh_k \in \bbR^d\enspace,
% \end{align*}
where $\bm_1, \ldots, \bm_k \in \bbN^d$. We note that since $f$ is $k$ times continuously differentiable in $\bbR^d$, then $f^{(k)}(\bx)$ is symmetric for all $\bx \in \bbR^d$. 
% \evg[inline]{Changed here the definition of the ball. It was a closed ball, now it is an open ball. The change is made to fit the Stokes' theorem. It does not modify any proof.}

\subsection{Classes of functions}
We start this section by stating all the relevant definitions and assumptions related to the target function $f$.
Following~\cite[Section 1.3]{nemirovski2000topics}, we recall the definition of higher order H\"older smoothness.

%\massi[inline]{related to the next remark for integer beta, we should recall here or in the notation that $\lfloor{\beta}\rfloor$  is the largest integer strictly smaller than $\beta$}
\begin{definition}[Higher order smoothness]
    \label{hh}
    Fix some $\beta \geq 2$ and $L > 0$.
    We denote by ${\mathcal F_\beta}(L)$ the set of all functions $f:\mathbb{R}^d\to \mathbb{R}$ that are $\ell=\lfloor \beta\rfloor$ times continuously differentiable and satisfy, for all $\bx,\bz \in \mathbb{R}^{d}$, the H\"older-type condition
        \begin{align*}
            \norm{f^{(\ell)}(\bx) - f^{(\ell)}(\bz)} \leq L \norm{\bx - \bz}^{\beta - \ell}\enspace.
        \end{align*}
\end{definition}
\begin{remark}
\label{rem:smoothness} Definition \ref{hh} of higher order smoothness was used by~\cite{BP2016} who considered only integer $\beta$, while 
\cite{PT90,akhavan2020,akhavan2021distributed} use a slightly different definition. Namely, they consider a class  $\class{F}_{\beta}'(L)$ defined as the set of all $\ell$ times continuously differentiable functions $f$ satisfying for all $\bx, \bz \in \bbR^d$, the condition
\begin{align*}
    |f(\bx) - T^{\ell}_{\bz}(\bx)| \leq L \norm{\bx - \bz}^{\beta}\enspace,
\end{align*}
where $T^{\ell}_{\bz}(\cdot)$ is the Taylor polynomial of order $\ell=\lfloor \beta\rfloor$ of $f$ around $\bz$. If $f \in \class{F}_{\beta}(L)$, then $f \in \class{F}_{\beta}'(L/\ell!)$ (cf. Appendix~\ref{app:A}). Thus, the results obtained for classes $\class{F}_{\beta}'$ hold true for $f \in \class{F}_{\beta}$ modulo a change of constant. Moreover, if $f$ is convex and
$\beta = 2$ (Lipschitz continuous gradient) the properties defining the two classes are equivalent to within constants, cf. \cite[Theorem 2.1.5]{nesterov-book}. 
\end{remark}
% \evg[inline]{Example that it is strictly larger?} 
% The above definition has been considered by

Since we study the minimization of highly smooth functions, in what follows, we will always assume that $f$ belongs to $\mathcal{F}_{\beta}(L)$ for some $\beta \geq 2$ and $L > 0$. We additionally require that $f \in \mathcal{F}_2(\bar{L})$ for some $\bL > 0$, that is, the gradient of $f$ is Lipschitz continuous.

\begin{assumption}\label{nat} The function $f \in \mathcal{F}_{\beta}(L) \cap \mathcal{F}_{2}(\bar{L})$ for some $\beta \geq 2$ and  $L, \bar{L} > 0$.
\end{assumption}

We will start our analysis 
%with the assumption that $f \in \mathcal{F}_{\beta}(L) \cap \mathcal{F}_{2}(\bar{L})$, 
by providing rates of convergence to a stationary point of $f$ under Assumption \ref{nat}. The first additional assumption that we consider is the Polyak-Łojasiewicz condition, which is also referred to as $\alpha$-gradient dominance. This condition became rather popular since it leads to linear convergence of the gradient descent algorithm without convexity as shown by~\cite{PB63} and further discussed by~\cite{10.1007/978-3-319-46128-1_50}.
\begin{definition}[$\alpha$-gradient dominance]\label{PL}
Let $\alpha>0$. Function $f:\mathbb{R}^{d}\to \mathbb{R}$ is called $\alpha$-gradient dominant on $\mathbb{R}^{d}$, if $f$ is differentiable on $\mathbb{R}^{d}$ and satisfies Polyak-Łojasiewicz inequality, 
\begin{equation}\label{eq:def-PL}
 2\alpha(f(\bx)-f^\star)\leq \norm{\nabla f(\bx)}^{2}\,,\qquad   \forall \bx \in \bbR^d \enspace.   
\end{equation}
\end{definition}
Finally, we consider the second additional condition, which is the $\alpha$-strong convexity. 
\begin{definition}[$\alpha$-strong convexity]\label{def:strong_conv}
Let $\alpha>0$. Function $f:\mathbb{R}^{d}\to \mathbb{R}$ is called $\alpha$-strongly convex on $\bbR^d$, if it is differentiable on $\mathbb{R}^{d}$ and satisfies
\[
f(\bx) \geq f(\bx') + \scalar{\nabla f (\bx')}{\bx - \bx'} + \frac{\alpha}{2}\norm{\bx - \bx'}^2\,,\qquad   \forall \bx, \bx' \in \bbR^d \enspace.
\]
\end{definition}
We recall that $\alpha$-strong convexity implies \eqref{eq:def-PL} ~\cite[see, \eg][Theorem 2.1.10]{nesterov-book}, and thus it is a more restrictive property than $\alpha$-gradient dominance.

An important example of family of functions satisfying the $\alpha$-dominance condition is given by composing strongly convex functions with a linear transformation. Let $n \in \bbN$, $\bfA \in \mathbb{R}^{n \times d}$ and define
\[
\mathcal{F}(\bfA)= \big\{f \,:\, f(\bx)=g(\bfA\bx) \text{,\,\, g} \text{ is }\alpha\text{-strongly convex}\big\}\enspace.
\]
Note that if $\bfA^\top \bfA$ is not invertible then the functions in $\mathcal{F}(\bfA)$ are not necessarily strongly convex. %For example, the case of $g(\bx) = \norm{\bx - \bb}^2$, with $\bb \in \mathbb{R}^{m}$ ($g$ is $2$-strongly convex). 
However, it can be shown that 
%\massi{You mean to say that any $f\in \mathcal{F}(\bfA)$ is $\alpha \gamma$-gradient dominant?} 
any $f \in \mathcal{F}(\bfA)$ is an $\alpha \gamma$-gradient dominant function, where $\gamma$ is the smallest non-zero singular value of $A$ \cite[see, \eg][]{10.1007/978-3-319-46128-1_50}. Alternatively, we can consider the following family of functions
\[
\mathcal{F}'(\bfA) = \big\{f ~:~ f(\bx)=g(\bfA\bx),\quad g \in C^2(\mathbb{R}^d),\quad g ~\text{ is strictly convex}\big\}\enspace,
\]
which is a set of $\alpha$-gradient dominant functions on any compact subset of $\mathbb{R}^{d}$, for some $\alpha>0$. A popular example of such a function appearing in machine learning applications is the logistic loss defined as
\(
g(\bfA\bx)=\sum_{i=1}^{n}\log(1+\exp(\ba_{i}^\top \bx)),
%\enspace,
\)
where for $1\leq i \leq n$, $\ba_{i}$ is $i$-th row of $\bfA$, and $\bx\in \mathbb{R}^{d}$. For this and more examples, see e.g.~\citep{lr} and references therein.

%\massi{Maybe this should be in the text, not a remark}
%\begin{remark}[On risk functions]
In what follows, we consider three different scenarios: \emph{(i)} the case of only smoothness assumption on $f$, \emph{(ii)} smoothness and $\alpha$-gradient dominance assumptions, \emph{(iii)} smoothness and $\alpha$-strong convexity assumptions.
Let $\tilde{\bx}$ be an output of any algorithm. For the first scenario, we obtain stationary point guarantee, that is, a bound on $\Exp[\norm{\nabla f (\tilde{\bx})}^2]$. For the second and the third scenarios, we provide bounds for the optimization error $\Exp[f(\tilde{\bx}) - f^{\star}]$. 

\begin{remark}\label{remark:estimation-risk}
  Note that under $\alpha$-strong convexity and the fact that $\nabla f(\bx^\star) = 0$, as well as under $\alpha$-gradient dominance~\cite[see, \eg][Appendix A]{10.1007/978-3-319-46128-1_50}, for any $\bx \in \bbR^d$ we have
\begin{align}\label{eq:opt-err-to-est-err}
    %\frac1{2\alpha}\norm{\nabla f (\bx)}^2 \ge 
    f(\bx) - f^{\star} \geq \frac{\alpha}{2}\norm{\bx - \bx^{*}}^2\enspace,
\end{align}
where $\bx^{*}$ is the Euclidean projection of $\bx$ onto the solution set $\argmin_{\bx \in \bbR^d} f(\bx)$ of the considered optimization problem, which is a singleton in case of strong convexity.
Thus, our upper bounds on $\Exp[f(\tilde{\bx}) - f^{\star}]$ obtained under $\alpha$-strong convexity or $\alpha$-gradient dominance imply immediately upper bounds for $\Exp[\norm{\tilde\bx - \bx^{*}}^2]$ with an extra factor $2/\alpha$.  
\end{remark}

%\end{remark}

\subsection{Classical stochastic optimization versus our setting}\label{sec:classical stoch opt}

%It is useful to highlight the differences between the usual zero-order stochastic optimization framework and our setting.

The classical zero-order stochastic optimization (CZSO) setting considered by \cite{NY1983,Nesterov2011,Ghadimi2013,duchi2015,
NS17,balasubramanian2021zeroth,Rando_Molinari_Villa_Rosasco22} among others assumes that there is a function $F:\bbR^d \times \bbR \to \bbR $ such that our target function is its expectation over the second argument,
$$
f(\bx)= \Exp [F(\bx,\xi)],
$$
where $\xi$ is a random variable. In order to find a minimizer of $f$, at each step $t$ of the algorithm one makes two queries and gets outputs of the form $(F(\bx_t,\xi_t), F(\bx_t+ h_t\bzeta_t,\xi_t))$, $t=1,2,\dots,$ where $\xi_t$'s are independent identically distributed (iid) realizations of random variable $\xi$, and $\bx_t$, $\bx_t+ h_t\bzeta_t$ are query points at step $t$. Here, $h_t>0$ is a perturbation parameter and $\bzeta_t$ is a random or deterministic perturbation. 
We emphasize two features of this setting:
\begin{itemize}
    \item[(a)] the two queries are obtained with the same random variable $\xi_t$;
    \item[(b)] the random variables $\xi_t$ are iid over $t$\footnote{Some papers, for example, \cite{Ghadimi2013,gasnikov_gradient-free_2016} relax this assumption.}.
\end{itemize}
Both (a) and (b) are not assumed in our setting. On the other hand, we assume additive noise structure. That is, at step $t$, we can only observe the values $F(\bz_t,\xi_t) = f(\bz_t) + \xi_t$ for any choice of query points $\bz_t$ depending only on the observations at the previous steps, but we do not assume that two queries are obtained with the same noise $\xi_t$ or are iid.  We deal with almost adversarial noise, see Assumption \ref{ass1} below. In particular, we do not assume that the noise is zero-mean. Thus, in general, we can have $\Exp [F(\bx,\xi_t)]\ne f(\bx)$. 

In the CZSO setting, the values $F(\bx_t,\xi_t)$ and $ F(\bx_t+ h_t\bzeta_t,\xi_t)$ are used to obtain gradient approximations involving the divided differences $(F(\bx_t+ h_t\bzeta_t,\xi_t) - F(\bx_t,\xi_t))/h_t$. A popular choice is the gradient estimator with Gaussian randomization suggested by \cite{Nesterov2011}:
\begin{align}
    \label{eq:grad_gauss}
    \bg_t^{G} &\triangleq \frac{1}{h_t}(F(\bx_t+ h_t\bzeta_t,\xi_t) - F(\bx_t,\xi_t))\bzeta_t\,, \quad \text{with}\quad \bzeta_t\sim \mathcal{N}(0,I_d)\enspace,
\end{align}
where $\mathcal{N}(0,I_d)$ denotes the standard Gaussian distribution in $\bbR^d$. 
In the case of additive noise, the divided differences are equal to $(f(\bx_t+ h_t\bzeta_t) - f(\bx_t))/h_t$, that is, the analysis of these algorithms reduces to that of noiseless (deterministic) optimization setting. When $\xi_t$'s are not additive, the assumptions that are often made in the literature on CZSO are such that the rates of convergence are the same as in the additive noise case, which is equivalent to noiseless case due to the above remark. 

We can summarize this discussion as follows:
\begin{itemize}
    \item the results obtained in the literature on CZSO, as well as some tools~\cite[e.g., averaging in Algorithm 1 of][]{balasubramanian2021zeroth}, do not apply in our framework;
    \item our upper bounds on the expected optimization error in the case of no noise ($\sigma=0$) imply identical bounds in the CZSO setting when the noise is additive. Under some assumptions made in the CZSO literature~\cite[e.g., under Assumption B of][]{duchi2015}, these bounds for $\sigma=0$ also extend to the general CZSO setting, with possible changes only in constants and not in the rates. In other cases the rates in CZSO setting can only be slower than ours~\cite[e.g., under Assumptions A1b, A2 of][]{Ghadimi2013}.
    % \evg{To discuss, might be not true.} \evg{In general it is probably only true when we can do online to batch conversion. For stationary point it is not obviously true.}
\end{itemize}

\section{Algorithms}\label{sec3}
Given a closed convex set $\com \subseteq \bbR^d$, we consider the following optimization scheme
\begin{align}
    \label{eq:algo_general}
    \bx_1 \in \Theta\qquad\text{and}\qquad \bx_{t+1}= \proj_{\com}(\bx_{t}- \eta_t \bg_t)\qquad t \geq 1\enspace,
\end{align}
where $\bg_t$ is an update direction, 
%which, ideally, mimic 
approximating the gradient direction $\nabla f(\bx_t)$ and $\eta_t > 0$ is a step-size.
Allowing one to perform two function evaluations per step, we consider two 
%distinct 
gradient estimators $\bg_t$ which are based on different randomization schemes. They both 
%estimators 
employ a smoothing kernel $K: [-1,1] \to \mathbb{R}$ which we assume to satisfy, for $\beta \geq 2$ and $\ell = \floor{\beta}$, the conditions
\begin{equation}\label{ker1}
\int K(r) \d r {=}0,  \int r K(r) \d r {=}1,  \int r^j K(r) \d r {=}0, \ j{=}2,\dots, \ell,~
\kappa_\beta {\triangleq} \int |r|^{\beta} |K(r)| \d r<\infty\enspace.
\end{equation}
In \citep{PT90} it was suggested to construct such kernels employing Legendre polynomials, in which case $\kappa_\beta \leq 2\sqrt{2}\beta$, cf. \citealp[Appendix A.3]{BP2016}.

We are now in a position to introduce the two estimators. 
%two procedures announced above.
Similarly to earlier works dealing with $\ell_2$-randomized methods~\citep[see e.g.,][]{NY1983,flaxman2004,BP2016,akhavan2020} we use gradient estimators based on a result, which is sometimes referred to as Stokes' theorem. A general form of this result, not restricted to the $\ell_2$ geometry, can be found in~\citealp[Appendix A]{Akhavan_Chzhen_Pontil_Tsybakov22a}.

\paragraph{\texorpdfstring{Gradient estimator based on $\ell_2$}{} randomization.}
At time $t \geq 1$, let $\bzeta_t^{\circ}$ be distributed uniformly on the $\ell_2$-sphere $\sphere^d_2$, let $r_t$ be uniformly distributed on $[-1, 1]$, and $h_t > 0$. Query two points:
\begin{align*}
    y_t = f(\bx_t+h_tr_t\bzeta_t^{\circ})+ \xi_t\qquad\text{and}\qquad y_t' = f(\bx_t-h_tr_t\bzeta_t^{\circ})+ \xi_t' \enspace,
\end{align*}
where $\xi_t,\xi'_t$ are noises.
Using the above feedback, define the gradient estimator as
\begin{align}
    \label{eq:grad_l2}
    \texttt{($\ell_2$ randomization)}\qquad\qquad&\bg_t^{\circ} \triangleq \frac{d}{2h_t}(y_t - y'_t)\bzeta_t^{\circ} K(r_t)\enspace.
\end{align}
We use the superscript $\circ$ to emphasize the fact that $\bg_t^{\circ}$ is based on the $\ell_2$ randomization. 
% Algorithm~\ref{algose} summarizes the iterative procedure based on $\bg_t^{\circ}$.

\paragraph{Gradient estimator based on \texorpdfstring{$\ell_1$}{} randomization.}

At time $t \geq 1$, let $\bzeta_{t}^{\diamond}$ be distributed uniformly on the $\ell_1$-sphere $\sphere^d_1$, let $r_t$ be uniformly distributed on $[-1, 1]$, and $h_t > 0$. Query two points:
\begin{align*}
    y_t = f(\bx_t+h_tr_t\bzeta_{t}^{\diamond})+ \xi_t\qquad\text{and}\qquad y_t' = f(\bx_t-h_tr_t\bzeta_{t}^{\diamond})+ \xi_t' \enspace.
\end{align*}
Using the above feedback, define the gradient estimator as
\begin{align}
    \label{eq:grad_l1}
    \texttt{($\ell_1$ randomization)}\qquad\qquad&\bg_t^{\diamond} \triangleq \frac{d}{2h_t}(y_t - y'_t)\sign(\bzeta_{t}^{\diamond}) K(r_t)\enspace.
\end{align}
%\evg[inline]{As a motivation for our new method we might want to look at: \url{https://openreview.net/pdf?id=BJe-DsC5Fm}}
We use the superscript $\diamond$ reminiscent of the form of the $\ell_1$-sphere in order to emphasize the fact that $\bg_t^{\diamond}$ is based on the $\ell_1$ randomization. The idea of using an $\ell_1$ randomization (different from \eqref{eq:grad_l1}) was probably first invoked  by~\cite{gasnikov_gradient-free_2016}. We refer to~\cite{Akhavan_Chzhen_Pontil_Tsybakov22a} who highlighted the potential computational and memory gains of another $\ell_1$ randomization gradient estimator compared to its $\ell_2$ counterpart, as well as its advantages in theoretical guarantees. The estimator of~\cite{Akhavan_Chzhen_Pontil_Tsybakov22a} differs from \eqref{eq:grad_l1} as it does not involve the kernel $K$ but the same computational and memory advantages remain true for the estimator \eqref{eq:grad_l1}. %

Throughout the paper, we impose the following assumption on the noises $\xi_t,\xi_t'$ and on the random variables that we generate in the estimators~\eqref{eq:grad_l2} and~\eqref{eq:grad_l1}. %which intuitively forces the oracle to select noise variables before observing the current query points.
\begin{assumption}
\label{ass1}
For all $t \in \{1,\dots,T\}$, it holds that:
\begin{enumerate}
    \item[(i)] the random variables $\xi_t$ and $\xi_t'$ are independent from $\bzeta_{t}^{\circ}$ (\resp $\bzeta_{t}^{\diamond}$) and from $r_t$ conditionally on $\bx_t$, and the random variables $\bzeta_t^{\circ}$ (\resp $\bzeta_{t}^{\diamond}$) and $r_t$ are independent;
    % \evg[inline]{Shouldn't we introduce the filtration and write this independence conditionally on the filtration??}
    \item[(ii)] $ \Exp [\xi_t^2]\le \sigma^2$ and  $\Exp [(\xi_t')^2]\le \sigma^2$, where $\sigma\ge 0$.
    % \evg[inline]{This one we also need conditionally on the filtration. No?}
\end{enumerate}
%\begin{itemize}   
%\item[(i)] For all $t=1,\dots,T$, the random variables $\xi_t$ and $\xi_t'$ are independent from $\zeta_t$ and from $r_t$, and the random variables $\zeta_t$ and $r_t$ are independent.
%\item[(ii)] It holds that $ \displaystyle{\max_{t=1,\dots,T}} \Exp [\xi_t^2]\le \sigma^2,$ and  $\displaystyle{\max_{t=1,\dots,T}} \Exp [(\xi_t')^2]\le \sigma^2$, where $\sigma\ge 0$. %\end{itemize}
\end{assumption}
Let us emphasize that we do not assume $\xi_t$ and $\xi_t'$ to have zero mean. Moreover, they can be non-random and no independence between noises on different time steps is required, so that the setting can be considered as \emph{almost} adversarial.
% \evg[inline]{I added this remark, which could be interesting for people in the individual sequence community.}
%\massi[inline]{"stochastic nature" and "full knowledge of the algorithm" are not very clear to me}
%Interestingly, since we do not necessarily need to assume the stochastic nature of the noise, 
Nevertheless, the first part of Assumption~\ref{ass1} does not permit a completely adversarial setup. Indeed, the oracle is not allowed to choose the noise variable depending on the current query points (i.e., the two perturbations of $\bx_t$). However, Assumption~\ref{ass1} encompasses the following protocol: before running the algorithm,
%any query is performed, 
the oracle fixes an arbitrary bounded (by $\sigma$) sequence $(\xi_t, \xi_t')_{t = 1}^T$ of noise pairs, possibly with full knowledge of the algorithm employed by the learner, and reveals this sequence query by query. 
%Note that the first part of the assumption does not permit a \emph{completely} adversarial setup (the Oracle is not allowed to choose the noise variable depending on the current query points). 
%As it was noted by~\cite{akhavan2020}, such a relaxed set of assumptions is possible because of randomization.

In the next two subsections, we study the bias and variance of the two estimators. 
As we shall see, the $\ell_1$ randomization can be more advantageous in the noiseless 
case than its $\ell_2$ counterpart (cf. Remark~\ref{rem:advantage_l1}).
%In what follows we will show that the algorithm based on the $\ell_1$-randomization can be more advantageous in the noiseless case than its $\ell_2$-counterpart. Namely, in Remark~\ref{rem:advantage_l1} we will show that if $\sigma = 0$, then both bias and variance of $\ell_1$-randomization are smaller than that of the $\ell_2$-randomization.

\subsection{Bias and variance of \texorpdfstring{$\ell_2$}{} randomization}
% \evg[inline]{The whole section is rewritten}
The next results allow us to control the bias and the second moment of gradient estimators $\bg_1^{\circ},\dots,\bg_T^{\circ}$, and play a crucial role in our analysis. 
\begin{lemma}[Bias of $\ell_2$ randomization]
\label{lem:bias_sphere}
Let Assumption~\ref{ass1} be fulfilled.
Suppose that $f \in \mathcal{F}_{\beta}(L)$ for some $\beta \geq 2$ and $L>0$. Let  $\bg_{t}^{\circ}$ be defined in~\eqref{eq:grad_l2} at time $t \geq 1$. Let $\ell = \floor{\beta}$. Then,
\begin{align}
    \norm{\Exp[\bg_{t}^{\circ}\mid \bx_t]-\nabla f(\bx_{t})} \leq \kappa_{\beta}\frac{L}{(\ell - 1)!}\cdot\frac{d}{d+\beta-1}h_{t}^{\beta-1}\enspace.
\end{align}
%where we recall that $\ell = \floor{\beta}$.
\end{lemma}
Intuitively, the smaller $h_t$ is, the more accurately $\bg_t$ estimates the gradient. A result analogous to Lemma~\ref{lem:bias_sphere} with a bigger constant was claimed in~\cite[Lemma 2]{BP2016} but the proof was not provided. The proof of Lemma~\ref{lem:bias_sphere} is presented in the appendix. It relies on the fact that $\bg_{t}^{\circ}$ is an unbiased estimator of some surrogate function, which is strongly related to the original $f$. The factor in front of $h_{t}^{\beta-1}$ in Lemma~\ref{lem:bias_sphere} is $O(1)$ as function of $d$. It should be noted that for $\beta > 2$ the bounds on the bias obtained in~\cite{akhavan2020} and~\cite{Gasnikov}, where the factors scale as $O(d)$ and $O(\sqrt{d})$, respectively, cannot be directly compared to Lemma~\ref{lem:bias_sphere}. This is due to the fact that those bounds are proved under a different notion of smoothness, cf. 
Remark~\ref{rem:smoothness}. Nevertheless, if $f$ is convex and $\beta = 2$ both notions of smoothness coincide, and 
Lemma~\ref{lem:bias_sphere} improves upon the bounds in \citep{akhavan2020} and \citep{Gasnikov} by factors of order $d$ and $\sqrt{d}$, respectively.
% It is implicitly present in~\cite{akhavan2020}.

% \subsubsection{New bound on the variance for $\ell_2$-randomizaiton}
%The next lemma emphasizes the trade-off between the bias and the variance term which does not permit taking $h_t$ arbitrary small.
% \evg[inline]{Inaccuracy: under assumption B, we need to replace $2\sigma^2$ by $\Exp[(\xi_t - \xi_t')^2 \mid \bx_t]$! Otherwise, need to assume that the second moment is bounded conditionally. Another way to alleviate thus issue is to take total expectation on both sides, then the statement becomes correct. Preferences?}
\begin{lemma}[Variance of $\ell_2$ randomization]
\label{lem:var_sphere_v2}
Let Assumption~\ref{ass1} hold and $f \,{\in}\,\mathcal{F}_{2}(\bar{L})$ for some $\bar{L}>0$. Then, for any $d \geq 2$,
%Let Assumption~\ref{ass1} be fulfilled. Assume that $f \in \mathcal{F}_{2}(\bar{L})$, then if $d \geq 2$
\begin{align*}
    \Exp[\|\bg_{t}^{\circ}\|^2]
    &\leq
   \frac{d^2\kappa}{d-1}\Exp\left[\parent{\norm{\nabla f (\bx_t)} + \bL h_t}^2\right] + \frac{d^2 \sigma^2\kappa}{h^2_t}\enspace,
\end{align*}
where $\kappa = \int_{-1}^1 K^2(r) \d r$.
\end{lemma}
The result of Lemma~\ref{lem:var_sphere_v2} can be simplified as
\begin{align}
    \label{eq:var_l2_for_proofs}
    \Exp[\|\bg_{t}^{\circ}\|^2]
    &\leq
   4d\kappa\Exp[\norm{\nabla f (\bx_t)}^2] + 4d\kappa\bL^2 h_t^2 + \frac{d^2 \sigma^2\kappa}{h^2_t}\,,\qquad d \geq 2\enspace.
\end{align}
Let us provide some remarks about this inequality. First, the leading term of order $d^2h_t^{-2}$ in \eqref{eq:var_l2_for_proofs} is the same as in~\cite[Lemma 2.4]{akhavan2020} and in~\cite[Appendix C1, beginning of the proof of Proposition 3]{BP2016},  but we obtain a better constant.  
%that is $d^2h_t^{-2}$ (with slightly better constant, which is almost for free).
The main improvement \wrt to both works lies in the lower order term. Indeed, unlike in those papers, the term $h_t^2$ is multiplied by $d$ instead of $d^2$.
This improvement is crucial for the condition $T \geq d^{2-\beta/2}$, under which we obtain our main results below. In particular, there is no condition on the horizon $T$ whenever $\beta\ge4$. On the contrary, we would need $T \geq d^3$ if we would used the previously known versions of the variance bounds \citep{BP2016,akhavan2020}. The proof of Lemma~\ref{lem:var_sphere_v2} presented below relies on the Poincaré inequality for the uniform distribution on $\sphere^d_2$. %which was exploited by~\cite{Akhavan_Chzhen_Pontil_Tsybakov22a} for the $\ell_1$-randomization.

\vspace{.2truecm}

\begin{proof}[Proof of Lemma~\ref{lem:var_sphere_v2}]
    For simplicity we drop the subscript $t$ from all the quantities.
    Using Assumption~\ref{ass1} and the fact that
\begin{equation}\label{eq:newl2_1_0}
  \Exp[f(\bx + hr \bzeta^{\circ}){-}f(\bx - hr \bzeta^{\circ}) \mid \bx, r] = 0  
\end{equation}
we obtain
       \begin{equation}
        \label{eq:newl2_1}
    \begin{aligned}
        \Exp[\|\bg^{\circ}\|^2]
        &=
        \frac{d^2}{4h^2}\Exp\left[\parent{f(\bx + hr \bzeta^{\circ})-f(\bx - hr \bzeta^{\circ}) + (\xi - \xi')}^2K^2(r)\right]\\
        &\leq
        \frac{d^2}{4h^2}\parent{\Exp\left[(f(\bx + hr \bzeta^{\circ})-f(\bx - hr \bzeta^{\circ})^2 K^2(r)\right] + 4\kappa\sigma^2}\enspace.
    \end{aligned}
    \end{equation}
    Since $f \in \mathcal{F}_2(\bL)$ and \eqref{eq:newl2_1_0} holds, then
    using Wirtinger-Poincaré inequality~\cite[see, \eg][(3.1)] %or Theorem 2, respectively]
    {Osserman78} %, Beckner89}
    we get
    \begin{align}
        \label{eq:newl2_2}
        \Exp\left[(f(\bx {+} hr \bzeta^{\circ}){-}f(\bx {-} hr \bzeta^{\circ}))^2 \big|~\bx, r\right]
        \leq \frac{h^2}{d{-}1}\Exp\big[\norm{\nabla f(\bx {+}hr \bzeta^{\circ}) {+} \nabla f(\bx {-} hr \bzeta^{\circ})}^2 \big|~\bx, r\big].
    \end{align}
    The fact that $f \in \class{F}_2(\bL)$ and the triangle inequality imply that
    \begin{align}
        \label{eq:newl2_3}
        \Exp\left[\norm{\nabla f(\bx +hr \bzeta^{\circ}) + \nabla f(\bx - hr \bzeta^{\circ})}^2 ~\big|~\bx, r\right]
        &\leq
        4\parent{\norm{\nabla f(\bx)} + \bL h}^2\enspace.
        % &=
        % 4\norm{\nabla f(\bx)}^2 + 8h\bL\norm{\nabla f(\bx)} + 4\bL h^2\\
        % &\leq
        % 8\norm{\nabla f(\bx)}^2 + 8\bL h^2\enspace.
    \end{align}
    Combining~\eqref{eq:newl2_1} -- \eqref{eq:newl2_3} proves the lemma.
\end{proof}

Note that 
%which was, for instance, used by~\cite{Gasnikov} to analyse the same algorithm, 
Eqs.~\eqref{eq:newl2_2}--\eqref{eq:newl2_3} imply, in particular, Lemma~9 of~\cite{Shamir17}. Our method based on Poincaré's inequality yields explicitly the constants in the bound. In this aspect, we improve upon~\cite{Shamir17}, where a concentration argument leads only to non-specified constants.

\subsection{Bias and variance of \texorpdfstring{$\ell_1$}{} randomization}
This section analyses the gradient estimate based on the $\ell_1$ randomization. The results below display a very different bias and variance behavior compared to the $\ell_2$ randomization.

\begin{lemma}[Bias of $\ell_1$ randomization]
\label{lem:bias_l1}\label{lem:bias_simplex}
Let Assumption~\ref{ass1} be fulfilled and
$f \in \mathcal{F}_{\beta}(L)$ for some $\beta \geq 2$ and $L>0$. Let  $\bg_{t}^{\diamond}$ be defined in~\eqref{eq:grad_l1}  at time $t \geq 1$. Let $\ell = \floor{\beta}$. Then,
\begin{align}
    \norm{\Exp[\bg_{t}^{\diamond}\mid \bx_t]-\nabla f(\bx_{t})} \leq
    % \kappa_{\beta}h_{t}^{\beta-1} \frac{L}{(\ell - 1)!} \frac{d^{\frac{\beta - 1}{2}}\Gamma(\beta)\Gamma(d)}{\Gamma(d + \beta - 1)}{
    % \leq
    L{\color{purple}c_{\beta}}\kappa_{\beta} \ell^{\beta-\ell} d^{\frac{1-\beta }{2}}h_{t}^{\beta-1}\enspace.
\end{align}
%where we recall that $\ell = \floor{\beta}$.
{When $2\leq \beta <3$, then $c_{\beta}=2^{\frac{\beta-1}{2}}$, and if $\beta \geq 3$ we have $c_{\beta}=1$.}
\end{lemma}
Notice that Lemma~\ref{lem:bias_l1} gives the same dependence on the discretization parameter $h_t$ as Lemma~\ref{lem:bias_sphere}.
However, unlike the bias bound in Lemma~\ref{lem:bias_sphere}, which is {dimension} independent, the result of Lemma~\ref{lem:bias_l1} depends on the dimension in a {favorable} way. In particular, the bias is controlled by a decreasing function of the dimension and this dependence becomes more and more favorable for smoother functions.
Yet, the price for such a favorable control of the bias is an inflated bound on the variance, which is established below.

    {
\begin{lemma}[Variance of $\ell_1$ randomization]
\label{lem:var_l1_v2}
Let Assumption~\ref{ass1} be fulfilled and 
$f \in \mathcal{F}_{2}(\bar{L})$ for some $\bar{L}>0$. Then, for any $d \geq 3$,
\begin{align*}
    \Exp[\|\bg_{t}^{\diamond}\|^2]
    &\leq
   % {\frac{{\bar{\tC}}_{1}{\color{purple}d^2}\kappa }{d-2}\Exp\norm{\nabla f(\bx_t)}^2 +
        % \frac{\bar{\tC}_{2}{\color{purple}d^2}\kappa\bL^2 h^2_t}{(d-2){\color{purple}(d + 1)}} + \frac{d^3\sigma^2\kappa}{h^2_t}}\enspace,
         \frac{8d^3\kappa}{(d+1)(d-2)}\Exp\left[\parent{\norm{\nabla f(\bx_t)} + \bL h_t\sqrt{\frac{{2}}{{d}}}}^2\right] + \frac{d^3\sigma^2\kappa}{h^2_t}\enspace,
\end{align*}
where $\kappa = \int_{-1}^1 K^2(r) \d r$.
% where $\bar{\tC}_{1} \leq 8\parent{1 + \sqrt{2}}^2$ and $\bar{\tC}_{2} = 8(2+\sqrt{2})^2$.
% Lemma~\ref{lem:shamir_likel1}.
\end{lemma}
% Let us first discuss the constant $\bar{\tC}_{1}$. First of all, it is important to observe that both $\bar{\tC}_{1}$ can be computed in practice as it only depends on the dimension. This fact is important for the eventual choice of $\eta_t$ and $h_t$ for the estimator~\eqref{eq:grad_l1} . 
% % Secondly, while $\bar{\tC}_{1}$ will be explicitly taken into account by the Algorithm~\ref{algo:simplex}, $\bar{\tC}_{2}$ will only have an implicit effect only appearing in the eventual upper bounds that we derive. Due to this fact it is important to highlight their behaviour.
% Note that, asymptotically we have
% \begin{align*}
%     &\lim_{d \rightarrow \infty} \bar{\tC}_{1} = 8(1 + \sqrt{2})^2 \leq 46.63\enspace.
%     % \qquad\text{and}\qquad
%     % \lim_{d \rightarrow \infty} \bar{\tC}_{2} = 16(3 + \sqrt{8}) \leq 93.26\enspace.
% \end{align*}
% Furthermore, $\bar{\tC}_{1}$ is increasing (hence, remains upper-bounded by $8(1 + \sqrt{2})^2$) with the growth of $d$.
% Note that $\bar{\tC}_{2} \leq 244.5$ for all $d \geq 3$. The subscript $d$ is chosen to highlight the asymptotic behaviour of this constant, which, for $d \rightarrow \infty$ behave as $\bar{\tC}_{2} \rightarrow 16(3 + \sqrt{8}) \leq 94$.
% Furthermore, note that for all $d \geq 3$ we have $\tfrac{d^2}{(d-2)(d - 1)} \leq 4.5$ and $\tfrac{d^3}{(d-2)(d^2-1)} \leq 3.5$. Hence,
Combined with the facts that $ a^2/((a-2)(a+1)) \leq 2.25$ for all $a \geq 3$ and $(a+b)^2 \leq 2a^2 + 2b^2$, the inequality of Lemma~\ref{lem:var_l1_v2} can be further simplified as
\begin{equation}
\label{eq:var_l1_for_proofs}
\begin{aligned}
    \Exp[\|\bg_{t}^{\diamond}\|^2]
    &\leq
    \frac{16d^3\kappa}{(d+1)(d-1)}\Exp\norm{\nabla f(\bx_t)}^2 +  \frac{32d^2\kappa\bL^2 h^2_t}{(d+1)(d-1)} + \frac{d^3\sigma^2\kappa}{h^2_t} \\
    &\leq
  36d\kappa\Exp\norm{\nabla f(\bx_t)}^2 +
        72{\kappa\bL^2 h^2_t} + \frac{d^3\sigma^2\kappa}{h^2_t}\,,\qquad d\geq 3\enspace.
\end{aligned}
\end{equation}
% Yet, since $\bar{\tC}_{1}$ is known explicitly, it can be implemented in practice and used directly for the choice of the discretization $h_t$ and the step-size $\eta_t$. For this reason, we keep the derived bound on the variance as is and do not rely on its simplified version from Eq.~\eqref{eq:var_l1_for_proofs}.

% \begin{remark}
% After seeing the proof of Lemma~\ref{lem:var_sphere_v2}, a very natural attempt to prove Lemma~\ref{lem:var_l1_v2} is to use a Poincaré inequality for the uniform measure on $\sphere^d_1$: for some $C_{d}$ that depends on $d$
% \begin{align*}
%     \Exp[g(\bzeta^{\diamond})^2] \leq C_{d} \Exp[ \norm{\nabla g (\bzeta^{\diamond})}^2 ]\qquad\text{for all}\qquad \Exp[g(\bzeta^{\diamond})] = 0\enspace
% \end{align*}
% A weighted version of the above inequality was recently derived in~\cite[Lemma 3]{Akhavan_Chzhen_Pontil_Tsybakov22a}, which we use here.
% \end{remark}
The proof of Lemma~\ref{lem:var_l1_v2}  is based on the following lemma,  which gives a version of Poincaré's inequality improving upon the previously derived in~\cite[Lemma 3]{Akhavan_Chzhen_Pontil_Tsybakov22a}. We provide sharper constants and an easier to use expression.
% \evg{I computed an improved version of Poincaré inequality. Main improvement here is that the eventual constants will be sharper}
\begin{lemma}
    \label{lem:true_poincare}
    Let $d \geq 3$. Assume that $G : \bbR^d \to \bbR$  is a continuously differentiable function, and $\bzeta^{\diamond}$ is distributed uniformly on $\sphere^d_1$. Then
    \begin{align*}
        \Var(G(\bzeta^{\diamond})) \leq \frac{4}{d-2}\Exp\left[\|\nabla G(\bzeta^{\diamond})\|^2\norm{\bzeta^{\diamond}}^2\right]\enspace.
    \end{align*}
     Furthermore, if $G : \bbR^d \to \bbR$ is an $L$-Lipschitz function w.r.t. the $\ell_2$-norm then 
    \begin{align*}
        \Var(G(\bzeta^{\diamond})) \leq \frac{8L^2}{(d+1)(d-2)} \leq \frac{18L^2}{d^2}\enspace.
    \end{align*}
\end{lemma}
The proof of Lemma~\ref{lem:true_poincare} is given in the Appendix. %While the  argument is similar to~\cite{Akhavan_Chzhen_Pontil_Tsybakov22a} the claimed improvement is achieved thanks to Lemma~\ref{lem:perturbation} that may be of independent interest. 
%Equipped with the above improved version of the Poincaré inequality, we are in position to prove Lemma~\ref{lem:var_l1_v2}.

\begin{proof}[Proof of Lemma~\ref{lem:var_l1_v2}]
For simplicity we drop the subscript $t$ from all the quantities. 
Similarly to the proof of Lemma~\ref{lem:var_sphere_v2}, using Assumption~\ref{ass1} we deduce that
    \begin{align}
        \label{lem:var_l1_v2_1}
        \Exp[\|\bg^{\diamond}\|^2]
        &\leq
        % \frac{d^3}{4h^2}\Exp\left[\parent{f(\bx + hr \bzeta^{\diamond})-f(\bx - hr \bzeta^{\diamond}) + (\xi - \xi')}^2K^2(r)\right]\\
        % &=
        \frac{d^3}{4h^2}\parent{\Exp[(f(\bx + hr \bzeta^{\diamond}) - f(\bx - hr \bzeta^{\diamond}))^2K^2(r)] + 4\sigma^2\kappa}\enspace.
        % &\leq
        % \frac{d^3}{4h^2}\parent{4\Exp[(f(\bx + hr \bzeta^{\diamond}) - \bar{f})^2K^2(r)] + \Exp[(\xi - \xi')^2K^2(r)]}\\
        % &\leq
        % \frac{d^3\kappa}{4h^2}\parent{\frac{\texttt{C}_{d,1}h^2}{(d-2)(d-1)}\norm{\nabla f(\bx)}^2 +
        % \frac{\texttt{C}_{d, 2}\bL^2 h^4}{(d-2)(d^2 - 1)} + 2\sigma^2}\enspace,
    \end{align}
    Consider $G : \bbR^d \to \bbR$ defined for all $\bu \in \bbR^d$ as $G(\bu) = f(\bx + hr\bu) - f(\bx - hr\bu)$. Using the fact that $f \in \mathcal{F}_2(\bL)$ we obtain for all $\bu \in \bbR^d$
    \begin{align*}
        \norm{\nabla G(\bu)}^2 \leq 4h^2\parent{\norm{\nabla f(\bx)} + \bL h\norm{\bu}}^2\enspace.
    \end{align*}
    % \begin{align*}
    %     \norm{\nabla G(\bu)}^2 \leq 8h^4\bL^2\norm{\bu}^2 + 8h^2\norm{\nabla f(\bx)}^2\enspace.
    % \end{align*}
     Applying Lemma~\ref{lem:true_poincare} to the function $G$ defined above we deduce that
    \begin{align*}
        \Exp\left[\left(G(\bzeta^{\diamond})\right)^2 \mid \bx, r\right]
        &\leq
        \frac{16h^2}{d-2}\Exp\left[\parent{\norm{\nabla f(\bx)} + \bL h\norm{\bzeta^{\diamond}}}^2{\norm{\bzeta^{\diamond}}^2}\right]\enspace.
        \end{align*}
    Lemma~\ref{lem:moments_of_strange_things} provided in the Appendix gives upper bounds on the expectations appearing in the above inequality for all $d \geq 3$. Its application yields:
    \begin{align*}
        \Exp \left[\big(f(\bx + hr \bzeta^{\diamond}) - f(\bx - hr \bzeta^{\diamond}) \big)^2 \mid \bx, r\right]
        \leq
        \frac{32h^2}{(d+1)(d-2)}\parent{\norm{\nabla f(\bx)} + \bL h\sqrt{\frac{2}{{d}}}}^2\enspace.
        % \frac{32h^2}{{\color{purple}d}(d-2)}&\parent{1 + \sqrt{\tfrac{2d}{d+1}}}^2\norm{\nabla f(\bx)}^2\\
        % &+
        % \frac{32(2+\sqrt{2})^2\bL^2 h^4}{{\color{purple}d}(d-2){\color{purple}(d + 1)}}
        % \enspace.
    \end{align*}
 We conclude the proof by combining this bound with~\eqref{lem:var_l1_v2_1}.
    % where the last inequality is due to Lemma~\ref{lem:var_l1_v3_intermediate} applied conditionally on $r$ with $\tC_{d, 1}$ and $\tC_{d, 2}$ defined in the statement of Lemma~\ref{lem:var_l1_v3_intermediate}.
\end{proof}
Modulo absolute constants, the leading term \wrt $h_t$ in Lemma~\ref{lem:var_l1_v2} is the same as for the $\ell_2$ randomization in Lemma~\ref{lem:var_sphere_v2}. However, for the $\ell_2$ randomization this term involves only a quadratic dependence on the dimension $d$, while in the case of $\ell_1$ randomization the dependence is cubic. Looking at the $h_t^2$ term that essentially matters only when $\sigma=0$, we note that 
the factor in front of this term in Lemma~\ref{lem:var_l1_v2} is bounded as the dimension grows. In contrast, the corresponding term in Lemma~\ref{lem:var_sphere_v2} depends linearly on the dimension. We summarize these observations in the following remark focusing on the noiseless case.

\begin{remark}[On the advantage of $\ell_1$ randomization]
\label{rem:advantage_l1}
In the noiseless case ($\sigma = 0$) both bias and variance under the $\ell_1$ randomization are strictly smaller than under the $\ell_2$ randomization. Indeed, if $\sigma = 0$
\begin{align*}
    \begin{cases}
         \norm{\Exp[\bg_{t}^{\circ}\mid \bx_t]-\nabla f(\bx_{t})} \lesssim h_{t}^{\beta-1}\\
         \Exp[\|\bg_{t}^{\circ}\|^2]
    \lesssim
   d\Exp[\norm{\nabla f (\bx_t)}^2] + (\sqrt{d}h_t)^2
    \end{cases}\quad\text{and}\quad
    \begin{cases}
         \norm{\Exp[\bg_{t}^{\diamond}\mid \bx_t]-\nabla f(\bx_{t})} \lesssim  \parent{\frac{h_{t}}{\sqrt{d}}}^{\beta-1}\\
         \Exp[\|\bg_{t}^{\diamond}\|^2]
    \lesssim d\Exp[\norm{\nabla f (\bx_t)}^2] + h_t^2
    \end{cases},
    \end{align*}
where the signs $\lesssim$ hide multiplicative constants that do not depend on $h_t$ and $d$.
% As a thought experiment, substitute $d=10^6$ and $\beta = 2$.
Notice that given an $h_t$ for $\ell_2$ randomization, $\sqrt{d}h_t$ used with $\ell_1$ randomization gives the same order of magnitude for both bias and variance. This is especially useful in the floating-point arithmetic, where very small values of $h_t$ (on the level of machine precision) can result in high rounding errors. Thus, in the absence of noise, $\ell_1$ randomization can be seen as a more numerically stable alternative to the $\ell_2$ randomization.
\end{remark}
For comparison, the corresponding bounds for gradient estimators with Gaussian randomization defined in \eqref{eq:grad_gauss}  are proved for $\beta=2$ and have the form (cf. \cite{Nesterov2011} or \cite[Theorem 3.1]{Ghadimi2013}):
\begin{align}\label{eq:bias-var-gaussian}
    \begin{cases}
         \norm{\Exp[\bg_{t}^{G}\mid \bx_t]-\nabla f(\bx_{t})} \lesssim d^{3/2} h_{t}\\
         \Exp[\|\bg_{t}^{G}\|^2]
    \lesssim
   d\Exp[\norm{\nabla f (\bx_t)}^2] + d^3h_t^2,
    \end{cases}
    \end{align}
    where the signs $\lesssim$ hide multiplicative constants that do not depend on $h_t$ and $d$.
Setting $\beta=2$ in Remark \ref{rem:advantage_l1} we see that the dependence on $h_t$ in \eqref{eq:bias-var-gaussian} is of the same order as for $\ell_1$ and $\ell_2$ randomizations with $\beta=2$, but the dimension factors are substantially bigger. {Also, the bounds in \eqref{eq:bias-var-gaussian} for the Gaussian randomization are tight. Thus, the Gaussian randomization is less efficient than its $\ell_1$ and $\ell_2$ counterparts in the noiseless setting.}

% \paragraph{A comment of the bias-variance trade-off.} In what follows, the discretization step $h_t$ will be chosen to provide an optimal balance between the variance and the squared bias. It is informative to observe that (keeping only the dominant terms and taking into account only the dependency on the dimension) $h_t$ should be selected to balance
% \begin{align*}
%     &h_t^{2\beta - 2} + d^2 h^{-2}_t &&\text{for $\ell_2$-randomization}\enspace,\\
%     &d^{1 - \beta}h_{t}^{2\beta - 2} + d^3 h^{-2}_t &&\text{for $\ell_1$-randomization}\enspace.
% \end{align*}
% Leading to $h_t \asymp d^{\frac{1}{\beta}}$ in the former case and to $h_t \asymp d^{\frac{1}{\beta} + \frac{1}{2}}$ in the latter case. Substituting these values, we deduce that in both cases the balancing choice of $h_t$ would lead to \emph{exactly} the same dependency on the dimension: $d^{\frac{2\beta - 2}{\beta}}$.

%%%%%%%%%%%%%%%%%%%%%%

% \newpage
%!TEX root = main.tex

\section{Upper bounds}\label{sec4}

In this section, we present convergence guarantees for the two considered gradient estimators and for three classes of objective functions $f$. Each of the following subsections is structured similarly: first, we define the choice of $\eta_t$ and $h_t$ involved in both algorithms and then, for each class of the objective functions, we state the corresponding convergence guarantees.

Throughout this section, we assume that  $f \in \mathcal{F}_2(\bar{L}) \cap \mathcal{F}_\beta(L)$ for some $\beta \geq 2$. Under this assumption, in Section~\ref{sec:non_convex} we establish a guarantee for the stationary point. In Section~\ref{sec:non_convex1} we additionally assume that $f$ is $\alpha$-gradient dominant and provide upper bounds on the optimization error.
In Section~\ref{sec5} we additionally assume that $f$ is $\alpha$-strongly convex and provide upper bounds on the optimization error for both constrained and unconstrained cases. %improving the upper bound derived by~\cite{akhavan2020} and~\cite{Gasnikov}.
Unless stated otherwise, the convergence guarantees presented in this section hold under the assumption that the number of queries $T$ is known before running the algorithms.

\subsection{Only smoothness assumptions}
\label{sec:non_convex}
In this subsection, we only assume that the objective function $f : \bbR^d \to \bbR$ satisfies Assumption~\ref{nat}. In particular, since there is no guarantee of the existence of the minimizer, our goal is modest -- we only want to obtain a nearly stationary point.

The plan of our study is as follows. We first obtain guarantees for algorithm~\eqref{eq:algo_general}
with gradient estimator $\bg_t$ satisfying some general assumption, and then concretize the results for the
gradient estimators~\eqref{eq:grad_l2} and~\eqref{eq:grad_l1}. We use the following assumption.
\begin{assumption}
\label{ass:grad_general_general}
Assume that there exist two positive sequences $b_t, v_t : \bbN \to [0, \infty)$ and ${\sf V}_1 \geq 0$ such that for all $t \geq 1$ it holds almost surely that
    \begin{align*}
        \norm{\Exp[\bg_t \mid \bx_t] - \nabla f(\bx_t)} \leq b_t\qquad\text{and}\qquad \Exp[\norm{\bg_t}^2] \leq  {\sf V}_1 \Exp[\norm{\nabla f(\bx_t)}^2]+v_t\enspace.
    \end{align*}
    \end{assumption}
    Note that Assumption \ref{ass:grad_general_general} holds for the
gradient estimators~\eqref{eq:grad_l2} and~\eqref{eq:grad_l1} with $b_t, v_t$ and ${\sf V}_1$ specified in Lemmas~\ref{lem:bias_sphere}--\ref{lem:var_l1_v2} (see also Assumption \ref{ass:grad_general} and Table \ref{table:par} below).

The results of this subsection will be stated on a randomly sampled point along the trajectory of the algorithm. The distribution over the trajectory is chosen carefully, in order to guarantee the desired convergence.
The distribution that we are going to use is defined in the following lemma.

\begin{lemma}\label{HighSF} Let $f \in \mathcal{F}_2(\bar{L})$ for some $\bar{L}>0$, $\com = \bbR^d$ and $f^{\star} > -\infty$.
    Let $\bx_t$ be defined  by algorithm~\eqref{eq:algo_general} with $\bg_t$ satisfying Assumption \ref{ass:grad_general_general}.
    Assume that $\eta_t$ in~\eqref{eq:algo_general} is chosen such that $\bL \eta_t {\sf V}_1 < 1$.
    Let $S$ be a random variable with values in $[T]$, which is independent from $\bx_1, \ldots, \bx_T, \bg_1, \ldots, \bg_T$ and distributed with the law
    \begin{align*}
        \Prob(S = t) = \frac{\eta_t\parent{1 - \bL\eta_t{\sf V}_1}}{\sum_{t = 1}^T\eta_t\parent{1 - \bL\eta_t{\sf V}_1}}, \quad t\in [T] \enspace.
    \end{align*}
    Then, 
    \begin{align*}
        \Exp[\norm{\nabla f(\bx_S)}^2] \leq \frac{2(\Exp[f(\bx_1)] - f^{\star}) + \sum_{t = 1}^T\eta_t\parent{b_t^2 + \bL\eta_t v_t}}{\sum_{t = 1}^T\eta_t\parent{1 - \bL\eta_t{\sf V}_1}}\enspace.
    \end{align*}
\end{lemma}
Lemma \ref{HighSF} is obtained by techniques similar to~\cite{Ghadimi2013}.  However, the paper \cite{Ghadimi2013} considers only a particular choice of $\bg_t$ defined via a Gaussian randomization, and a different setting (cf. the discussion in Section \ref{sec:classical stoch opt}), under which $v_t$ does not increase as the discretization parameter $h_t$ ($\mu$ in the notation of~\cite{Ghadimi2013}) decreases. In our setting, this situation happens only when there is no noise ($\sigma=0$), while in the noisy case $v_t$ increases as $h_t$ tends to 0.

Note that the distribution of $S$ in Lemma~\ref{HighSF} depends on the choice of $\eta_t$ and ${\sf V}_1$. 
In the following results, we are going to specify the exact values of $\eta_t$. We also provide the values of ${\sf V}_1$ for the gradient estimators~\eqref{eq:grad_l2} and~\eqref{eq:grad_l1}. Regarding these two estimators, it will be convenient to use the following instance of Assumption \ref{ass:grad_general_general}. 
\begin{assumption}
\label{ass:grad_general}
There exist positive numbers $b, {\sf V}_1, {\sf V}_2, {\sf V}_3$ such that for all $t \geq 1$ the gradient estimators $\bg_t$ satisfy almost surely the inequalities
 \begin{align*}
        \norm{\Exp[\bg_t \mid \bx_t] - \nabla f(\bx_t)} \leq bLh_t^{\beta-1}\quad\text{and}\quad \Exp[\norm{\bg_t}^2] \leq {{\sf V}_1} \Exp[\norm{\nabla f(\bx_t)}^2] + {{\sf V}_2}\bar{L}^2h_t^{2} + {{\sf V}_3}\sigma^2h_t^{-2}.
\end{align*}
\end{assumption}
It follows from Lemmas~\ref{lem:bias_sphere}--\ref{lem:var_l1_v2} that Assumption~\ref{ass:grad_general} holds for gradient estimators~\eqref{eq:grad_l2} and~\eqref{eq:grad_l1} with the values that are indicated in Table~\ref{table:par}.
Note that the bounds for the variance in those lemmas do not cover the case $d=1$ for the $\ell_2$ randomization
and $d=1,2$ for the $\ell_1$ randomization. Nevertheless, it is straightforward to check that in these cases Assumption~\ref{ass:grad_general} remains valid with  ${\sf V}_j$'s given in Table 1.

\begin{table}[h]
\centering
\begin{tabular}{@{}lllll@{}}
\toprule
                $\phantom{10}$Estimator       & $b$                                                    & ${\sf V}_1\phantom{1000000000}$  & ${\sf V}_2\phantom{1000000000}$  & ${\sf V}_3\phantom{10}$ \\ \midrule
\texttt{$\phantom{10}\ell_2$ randomization} & $\frac{\kappa_{\beta}}{(\ell - 1)!}\cdot\frac{d}{d+\beta-1}$ & $4d\kappa$   & $4d\kappa$   & $d^2\kappa$ \\[.3cm]
\texttt{$\phantom{10}\ell_1$ randomization} & ${\color{purple}c_{\beta}}\kappa_{\beta}\ell^{\beta-\ell}d^{\frac{1-\beta}{2}}$       & $36d\kappa$ & $72\kappa$ & $d^3\kappa$
\end{tabular}
\caption{Factors in the bounds for bias and variance of both gradient estimators, $\ell=\lfloor \beta \rfloor$, $d\ge1$. }\label{table:par}
\end{table}

The next theorem requires a definition of algorithm-dependent parameters, which are needed as an input to our algorithms. We set
\begin{align}
\label{eq:constants_only_smoothness}
    (\mathfrak{y}, \mathfrak{h})
    % (\mathfrak{h}, \mathfrak{n})
    = 
    \begin{cases}
    \bigg((8\kappa\bL)^{-1},\,d^{\frac{1}{2\beta - 1}}\bigg)
    &\qquad\text{for $\ell_2$ randomization,}\\
    % \\
    \bigg((72\kappa\bar{L})^{-1},\, d^{\frac{2\beta + 1}{4\beta - 2}}\bigg)
    &\qquad\text{for $\ell_1$ randomization.}
    \end{cases}
\end{align}
\begin{theorem}\label{nonc}
Let Assumptions \ref{nat} and \ref{ass1} hold, and $\com = \bbR^d$. Let $\bx_t$ be defined  by algorithm~\eqref{eq:algo_general} with gradient estimator~\eqref{eq:grad_l2} or~\eqref{eq:grad_l1}, where the  parameters $\eta_t$ and $h_{t}$ are set for $t = 1, \ldots, T$, as
\begin{align*}
    \eta_t =\min\left(\frac{\mathfrak{y}}{d},\, d^{-\frac{2(\beta-1)}{2\beta-1}}T^{-\frac{\beta}{2\beta-1}}\right) \qquad\text{and}\qquad 
    h_t = \mathfrak{h}\, T^{-\frac{1}{2(2\beta-1)}}\enspace, 
\end{align*}
and the constants $\mathfrak{y}$ and $\mathfrak{h}$ are given in~\eqref{eq:constants_only_smoothness}. 
Assume that $\bx_{1}$ is deterministic and $T\geq d^{\frac{1}{\beta}}$. Then, for the random variable $S$ defined in Lemma~\ref{HighSF}, we have
\begin{align*}
     \Exp[\norm{\nabla f(\bx_S)}^2] \leq  \Big(\cst_1(f(\bx_1)-f^\star)+\cst_2\Big)\left(\frac{d^2}{T}\right)^{\frac{\beta-1}{2\beta-1}}\enspace,
\end{align*}
where the constants $\cst_1, \cst_2>0$ depend only on $\sigma, L, \bar{L}$, $\beta$, and on the choice of the gradient estimator.
\end{theorem}

{\color{brown}
In the case $\sigma = 0$, the result of this theorem can be improved.  As explained in Section \ref{sec:classical stoch opt}, this case is analogous to the CZSO setting, and it is enough to assume that $\beta=2$ since higher order smoothness does not lead to improvement in the main term of the rates. Due to Remark~\ref{rem:advantage_l1} (or Assumption~\ref{ass:grad_general}) one can set $h_t$ for both methods as small as one 
wishes, and thus sufficiently small to make the sum over $t$ in the numerator of the inequality of Lemma~\ref{HighSF}  
less than an absolute constant. Then, choosing $\eta_t = (2{\sf V}_1\bL)^{-1}$ and recalling that, for both algorithms, ${\sf V}_1$ 
scales as $d$ up to a multiplicative constant (cf. Table \ref{table:par}) we get the following result.

\begin{theorem}\label{thm:OSsigma=0}
Let $f$ be a function belonging to $\mathcal{F}_2(\bar L)$ for some $\bar L>0$, $\com = \bbR^d$, and let Assumptions \ref{nat} and \ref{ass1} hold with $\sigma=0$. 
Let $\bx_t$ be defined by algorithm~\eqref{eq:algo_general} with deterministic $\bx_1$,  and gradient estimators~\eqref{eq:grad_l2} or~\eqref{eq:grad_l1} for $\beta=2$, where $\eta_t = (2{\sf V}_1\bL)^{-1}$ and $h_{t}$ is chosen sufficiently small. 
Then we have 
\begin{align*}
        \Exp[\norm{\nabla f(\bx_S)}^2] \leq \cst(f(\bx_1) - f^{\star}+1)\frac{\bL d}{T}\enspace,
\end{align*}
where $\cst>0$ is an absolute constant depending only the choice of the gradient estimator.
\end{theorem}
}

The rate $O(d/T)$ in Theorem \ref{thm:OSsigma=0} coincides with the rate derived in~\cite[inequality (68)]{NS17} for $\beta=2$ under the classical zero-order stochastic optimization setting, where the authors were using Gaussian rather than $\ell_1$ or $\ell_2$ randomization. 
In a setting with non-additive noise,~\cite{Ghadimi2013} exhibit a slower rate of $O(\sqrt{d/T})$.
% \evg{Doublecheck. I looked briefly---the setting is not the same. Ghadimi and Lan are in stochastic, non-additive setup. Nesterov in noiseless (which is the same as additive + 2 points per noise level). I would say: \cite{Ghadimi2013} exhibit a lower rate of $O(\sqrt{d/T})$ for stochastic case with non-additive noise.}

\subsection{Smoothness and \texorpdfstring{$\alpha$}{}-gradient dominance}
\label{sec:non_convex1}

We now provide the analysis of our algorithms under smoothness and \texorpdfstring{$\alpha$}{}-gradient dominance (Polyak-Łojasiewicz) conditions.

\begin{theorem}\label{thm:master_dominant}
Let $f$ be an $\alpha$-gradient dominant function, $\com = \bbR^d$, and let Assumptions \ref{nat} and \ref{ass1} hold, {\color{purple}with $\sigma >0$}.
 Let $\bx_t$ be defined by algorithm~\eqref{eq:algo_general} with $\bg_t$ satisfying Assumption~\ref{ass:grad_general}, 
 deterministic $\bx_1$ and  
\begin{align*}
    \eta_t = \min\left((2\bar{L}{{\sf V}_1})^{-1},\, \frac{4}{\alpha t}\right),
    \qquad
    h_t = \left( \frac{4\bar{L}{\color{purple}\sigma^2}{\sf V}_3}{b^2{\color{purple}L^2}\alpha}\right)^\frac{1}{2\beta}\cdot
 \begin{cases}
        t^{-\frac{1}{2\beta}} &\text{ if \,$\eta_t = \frac{4}{\alpha t}$}\\
        T^{-\frac{1}{2\beta}} &\text{ if \,$\eta_t = \frac{1}{2\bar{L}{{\sf V}_1}}$}
    \end{cases}\enspace.
\end{align*}
 Then
\begin{align*}
    \Exp[f(\bx_T)-f^\star]
    \leq 
    \cst_1&\,\frac{{\bL{\sf V}_1}}{\alpha T}(f(\bx_1)-f^\star)\\&+\frac{\cst_2}{\alpha}\left({{\sf V}_3}\left(\frac{{{\sf V}_3}}{b^2{\color{purple}L^2}}\right)^{-\frac{1}{\beta}}+{{\sf V}_2}{\color{purple}\bL^2}\left(\frac{{{\sf V}_3}}{b^2{\color{purple}L^2}}\right)^{\frac{1}{\beta}}\left(\frac{\alpha T}{\bL{\color{purple}\sigma^2}}\right)^{-\frac{2}{\beta}}{\color{purple}\sigma^{-2}}\right)\left(\frac{\alpha T}{\bL{\color{purple}\sigma^2}}\right)^{-\frac{\beta-1}{\beta}},
\end{align*}
where $\cst_{1}, \cst_2>0$ depend only on $\beta$.

\end{theorem}
Theorem~\ref{thm:master_dominant} provides a general result for any gradient estimator that satisfies Assumption~\ref{ass:grad_general}. By taking the values ${\sf V}_j$ from Table~\ref{table:par} we immediately obtain the following corollary for our $\ell_1$- and $\ell_2$-randomized gradient estimators.
%. 
\begin{corollary}\label{cor:finalres_PL}
Let $f$ be an $\alpha$-gradient dominant function, $\com = \bbR^d$, and let Assumptions \ref{nat} and \ref{ass1} hold, {\color{purple}with $\sigma >0$}. 
Let $\bx_t$ be defined by algorithm~\eqref{eq:algo_general} with deterministic $\bx_1$ and gradient estimators~\eqref{eq:grad_l2} or~\eqref{eq:grad_l1}. Set the parameters $\eta_t$ and $h_{t}$ as in Theorem \ref{thm:master_dominant}, where $b, {{\sf V}_1},{{\sf V}_2},{{\sf V}_3}$ are given in Table \ref{table:par} for each gradient estimator, respectively. Then for any $T\geq d^{2-\frac{\beta}{2}}\frac{\sigma^2}{\alpha L^2}$ we have 
\begin{align*}
    \Exp[f(\bx_T) - f^{\star}] \leq \cst_1\,\frac{\bL d}{\alpha T}\left(f(\bx_1)-f^{\star}\right) + 
    \left({\color{purple}\cst_2 + \cst_3\frac{\bL^2}{\sigma^2}}\right)
    \left(\frac{\bL {\color{purple}\sigma^2}d^2}{\alpha T}\right)^{\frac{\beta-1}{\beta}}{\color{purple}\frac{L^{\frac{2}{\beta}}}{\alpha}}\enspace,
\end{align*}
where $\cst_1,\cst_2,\cst_3>0$ depend only on $\beta$ and on the choice of the gradient estimator.
\end{corollary}
{\color{purple}Note that here we consider $\sigma$ and $L$ as numerical constants. The condition $T\gtrsim d^{2-\frac{\beta}{2}}/\alpha$ mentioned in Corollary \ref{cor:finalres_PL} is satisfied in all reasonable cases since  it is weaker than the condition $T\gtrsim d^2/\alpha$ guaranteeing non-triviality of the bounds.
}

Recall that, in the context of deterministic optimization with first order oracle, the $\alpha$-gradient dominance allows one to obtain the rates of convergence of gradient descent algorithm, which are similar to the case of strongly convex objective function with Lipschitz gradient (\cite{PB63,10.1007/978-3-319-46128-1_50}). A natural question is whether the same property holds in our setting of stochastic optimization with zero-order oracle and higher order smoothness. {\color{brown} Theorem~\ref{thm:master_dominant} shows the rates are only inflated by a multiplicative factor $\mu^{({\beta-1})/{\beta}}$, where $\mu = {\bar L}/\alpha$, compared to the $\alpha$-strongly convex case that will be considered in Section \ref{sec5}. }

Consider now the case $\sigma=0$, which is analogous to the CZSO setting as explained in Section \ref{sec:classical stoch opt}. In this case, we assume that $\beta=2$ since higher order smoothness does not lead to improvement in the main term of the rates. We set the parameters $\eta_t,h_t$ as follows:
\begin{align}\label{eq:eta-h-PL-sigma0}
    \eta_t = \min\left((2\bar{L}{\sf V}_1)^{-1},\, \frac{4}{\alpha t}\right),
    \qquad
    h_t \le \left(\frac{\bL\vee 1}{\alpha\wedge 1}T\left(2b^2{\color{purple}\bL} + \frac{8 {\color{purple}\bL^2}{\sf V}_2}{\alpha}\right)\right)^{-\frac{1}{2}}\enspace.
\end{align}
\begin{theorem}\label{thm:PL-sigma=0}
Let $f$ be an $\alpha$-gradient dominant function belonging to $\mathcal{F}_2(\bar L)$ for some $\bar L>0$, $\com = \bbR^d$, and let Assumptions \ref{nat} and \ref{ass1} hold with $\sigma=0$. 
Let $\bx_t$ be defined by algorithm~\eqref{eq:algo_general} with deterministic $\bx_1$,  and gradient estimators~\eqref{eq:grad_l2} or~\eqref{eq:grad_l1} for $\beta=2$. Set the parameters $\eta_t$ and $h_{t}$ as in \eqref{eq:eta-h-PL-sigma0}. %, \eqref{eq:gamma0}, where $b, {{\sf V}_1},{{\sf V}_2},{{\sf V}_3}$ are given in Table \ref{table:par} for each gradient estimator, respectively. 
Then we have 
\begin{align*}
    \Exp[f(\bx_T) - f^{\star}] \leq 
    \cst_1\frac{\bL d}{\alpha T}\left(\left(f(\bx_1)-f^{\star}\right) + \cst_2\right)\enspace,
\end{align*}
where $\cst_1,\cst_2>0$ are absolute constants depending only the choice of the gradient estimator.
\end{theorem}
Note that, in the CZSO setting, \cite{Rando_Molinari_Villa_Rosasco22} proved the rate $O(T^{-1})$ for the optimization error under $\alpha$-gradient dominance by using an $\ell_2$ randomization gradient estimator. However, unlike Theorem \ref{thm:PL-sigma=0} the bound obtained in that paper does not provide the dependence on the dimension $d$ and on variables $\bL, \alpha$.

\subsection{Smoothness and strong convexity}\label{sec5}
In this subsection, we additionally assume that $f$ is a strongly convex function and denote by $\bx^\star$ its unique minimizer.  
We provide a guarantee on the weighted average point $\hat{\bx}_T$ along the trajectory of the algorithm defined as
\[
    \hat{\bx}_T= \frac{2}{T(T+1)}\sum_{t =1}^T t \bx_t.
\] 
We consider separately the cases of unconstrained and constrained optimization. 

\subsubsection{Unconstrained optimization}
In this part we assume that $\com = \bbR^d$ and the horizon $T$ is known to the learner. 
Similar to Section \ref{sec:non_convex1}, we first state a general result that can be applied to any gradient estimator satisfying Assumption~\ref{ass:grad_general}. 

\begin{theorem}\label{thm:master_str_convex_uncostrained}
Let $f$ be an $\alpha$-strongly convex function, $\com = \bbR^d$, and let Assumptions \ref{nat} and \ref{ass1} hold. 
Let $\bx_t$ be defined by algorithm~\eqref{eq:algo_general} with $\bg_t$ satisfying Assumption~\ref{ass:grad_general}, 
 deterministic $\bx_1$ and  
\begin{align*}
    \eta_t = \min\left(\frac{\alpha}{8\bar{L}^2{{\sf V}_1}},\, \frac{4}{\alpha (t+1)}\right),
    \qquad
    h_t = \left(\frac{4{\color{purple}\sigma^2}{{\sf V}_3}}{b^2{\color{purple}L^2}}\right)^\frac{1}{2\beta}\cdot
 \begin{cases}
        t^{-\frac{1}{2\beta}} &\text{ if $\eta_t = \frac{8}{\alpha (t+1)}$}\\
        T^{-\frac{1}{2\beta}} &\text{ if $\eta_t = \frac{\alpha}{4\bar{L}^2{{\sf V}_1}}$}
    \end{cases}\enspace.
\end{align*}
Then
\begin{align*}
     \Exp[f(\hat{\bx}_T)-f^\star] \leq \cst_1\frac{{\color{purple}\bL^2}{\sf V}_1}{\alpha T}\norm{\bx_1{-}\bx^\star}^{2}+\Bigg\{\cst_2&(b{\color{purple}L})^{\frac{2}{\beta}}({{\sf V}_3}{\color{purple}\sigma^2})^{\frac{\beta-1}{\beta}}\\
     &+
     \cst_3{{\sf V}_2}{\color{purple}\bL^2}\left(\frac{{{\sf V}_3{\color{purple}\sigma^2}}}{b^2{\color{purple}L^2}}\right)^{\frac{1}{\beta}}T^{-\frac{2}{\beta}}\Bigg\}\frac{T^{-\frac{\beta-1}{\beta}}}{\alpha}\,,
\end{align*}
where the constants $\cst_1, \cst_2, \cst_3>0$ depend only on $\beta$.
\end{theorem}
Subsequently, in Corollary \ref{cor:finalres_str_convex}, we customize the above theorem for gradient estimators~\eqref{eq:grad_l2} and~\eqref{eq:grad_l1}, with assignments of $\eta_t, h_t$ that are again selected based on Table \ref{table:par}. We also include a bound for $\Exp[\norm{\hat{\bx}_T-\bx^\star}^2]$, which comes as an immediate consequence due to \eqref{eq:opt-err-to-est-err}. 
\begin{corollary}\label{cor:finalres_str_convex}
Let $f$ be an $\alpha$-strongly convex function, $\com = \bbR^d$, and let Assumptions \ref{nat} and \ref{ass1} hold. 
Let $\bx_t$ be defined by algorithm~\eqref{eq:algo_general} with  gradient estimator~\eqref{eq:grad_l2} or~\eqref{eq:grad_l1}, and parameters $\eta_t$, $h_{t}$ as in Theorem \ref{thm:master_str_convex_uncostrained}, where $b, {{\sf V}_1},{{\sf V}_2},{{\sf V}_3}$ are given in Table \ref{table:par} for each gradient estimator, respectively. Let $\bx_1$ be deterministic.
Then for any $T \geq d^{2-\frac{\beta}{2}}\frac{\sigma^2}{L^2}$ we have
\begin{align}\label{eq-1:cor:finalres_str_convex}
    &\Exp[f(\hat\bx_T) - f^{\star}] \leq \cst_1\,\frac{{\color{purple}\bL^2}d}{\alpha T}\norm{\bx_1-\bx^{\star}}^{2} + \left(\cst_2+\cst_3{\color{purple}\frac{\bL^2}{\sigma^2}}\right)\left(\frac{d^2{\color{purple}\sigma^2}}{T}\right)^{\frac{\beta-1}{\beta}}{\color{purple}L^{\frac{2}{\beta}}}\alpha^{-1},
    \\
    \label{eq-2:cor:finalres_str_convex}
    &\Exp[\norm{\hat{\bx}_T-\bx^\star}^2]  \le 2\cst_1\,\frac{{\color{purple}\bL^2}d}{\alpha^2 T}\norm{\bx_1-\bx^{\star}}^{2} + 2\left(\cst_2+\cst_3{\color{purple}\frac{\bL^2}{\sigma^2}}\right)\left(\frac{d^2{\color{purple}\sigma^2}}{T}\right)^{\frac{\beta-1}{\beta}}{\color{purple}L^{\frac{2}{\beta}}}\alpha^{-2},
\end{align}
where $\cst_1,\cst_2,\cst_3>0$ depend only on $\beta$ and on the choice of the gradient estimator.
\end{corollary}
With a slightly different definition of smoothness class (which coincides with ours for $\beta=2$, cf. Remark \ref{rem:smoothness}), a result comparable to Corollary~\ref{cor:finalres_str_convex}
is derived in~\cite[Theorem 3.2]{akhavan2020}. However, that result imposes an additional condition on $\alpha$ (i.e., $\alpha \gtrsim \sqrt{{d}/{T}}$) and provides a bound with the dimension factor $d^2$ rather than $d^{2-2/\beta}$ in Corollary~\ref{cor:finalres_str_convex}. We also note that earlier \cite{BP2016} analyzed the case $\ell_2$-randomized gradient estimator with integer $\beta>2$ and proved a bound with a slower (suboptimal) rate $T^{-\frac{\beta-1}{\beta+1}}$.
\begin{comment}
\begin{align}\label{eq:eta-h-PL-sigma0}
    \eta_t = \min\left(\frac{\alpha}{2\bar{L}^2{{\sf V}_1}},\, \frac{4}{\alpha t}\right),
    \qquad
    h_t \le \left(\frac{T}{d}\left(b^2 +4\bL^2{\sf V}_2\right)\right)^{-\frac{1}{2}}\enspace.%left( \frac{4\bar{L}\gamma^2{\sf V}_3}{b^2}\right)^\frac{1}{4}\cdot
 %\begin{cases}
       % t^{-\frac{1}{4}} &\text{ if \,$\eta_t = \frac{4}{\alpha t}$}\\
        %T^{-\frac{1}{4}} &\text{ if \,$\eta_t = \frac{1}{2\bar{L}{{\sf V}_1}}$}
   % \end{cases}\enspace.
\end{align}
\end{comment}
%
%
\subsubsection{Constrained optimization}
We now assume that $\com \subset \bbR^d$ is a compact convex set. 
In the present part, we do not need the knowledge of the horizon $T$ to define the updates $\bx_t$. 
We first state the following general theorem valid when $\bg_t$ is any gradient estimator satisfying Assumption~\ref{ass:grad_general}. 
\begin{theorem}
\label{scunt}
   Let $\com \subset \mathbb{R}^d$ be a compact convex set. Assume that $f$ is an $\alpha$-strongly convex function, Assumptions \ref{nat} and \ref{ass1} hold, and $\max_{\bx\in\com}\norm{\nabla f(\bx)}\le G$. 
Let $\bx_t$ be defined by algorithm~\eqref{eq:algo_general} with gradient estimator $\bg_t$ satisfying Assumption~\ref{ass:grad_general} and  
$
    \eta_t = \frac{4}{\alpha (t+1)},
    h_t = \left(\frac{{\color{purple}\sigma^2}{{\sf V}_3}}{b^2{\color{purple}L^2} t}\right)^\frac{1}{2\beta}.
$
Then
\begin{align*}
   \Exp[f(\hat{\bx}_t)-f^\star]&\leq \frac{4{\color{purple}\bL^2}{\sf V}_1 G^2}{\alpha T} + \frac{\cst_1}{\alpha}\left({\sf V}_3{\color{purple}\sigma^2}\left(\frac{{\sf V}_3{\color{purple}\sigma^2}}{b^2{\color{purple}L^2}}\right)^{-\frac{1}{\beta}} + {\sf V}_2{\color{purple}\bL^2}\left(\frac{{\sf V}_3{\color{purple}\sigma^2}}{b^2{\color{purple}L^2}}\right)^{\frac{1}{\beta}}T^{-\frac{2}{\beta}}\right)T^{-\frac{\beta-1}{\beta}}\enspace,
\end{align*}
where the constant $\cst_1>0$ depends only on $\beta$.

\end{theorem}
Using the bounds on the variance and bias of gradient estimators~\eqref{eq:grad_l2} and~\eqref{eq:grad_l1} from Section~\ref{sec3}, Remark~\ref{remark:estimation-risk} and the trivial bounds $\Exp[f(\hat{\bx}_T)-f^\star]\le GB$, $\Exp[\norm{\hat{\bx}_T-\bx^\star}^2]\le B^2$, where $B$ is the Euclidean diameter of $\com$, we immediately obtain  the following corollary.

\begin{corollary}\label{thm:sconvex_cons}
Let $\com \subset \mathbb{R}^d$ be a compact convex set. Assume that $f$ is an $\alpha$-strongly convex function, Assumptions \ref{nat} and \ref{ass1} hold, and $\max_{\bx\in\com}\norm{\nabla f(\bx)}\le G$. 
Let $\bx_t$ be defined by algorithm~\eqref{eq:algo_general} with  gradient estimator~\eqref{eq:grad_l2} or~\eqref{eq:grad_l1}, and parameters $\eta_t$, $h_{t}$ as in Theorem \ref{scunt}, where $b, {{\sf V}_1},{{\sf V}_2},{{\sf V}_3}$ are given in Table \ref{table:par} for each gradient estimator, respectively. %Let $\bx_1$ be deterministic.
Then for any $T \geq d^{2-\frac{\beta}{2}}\frac{\sigma^2}{L^2}$ we have
\begin{align}\label{eq1:thm:sconvex_cons}
    \Exp[f(\hat{\bx}_T)-f^\star] &\leq \min\bigg(GB, \frac{4{\color{purple}\bL^2}{\sf V}_1G^2}{\alpha T}+\left(\cst_1+\cst_2{\color{purple}\frac{\bL^2}{\sigma^2}}\right)\left(\frac{d^2{\color{purple}\sigma^2}}{T}\right)^{\frac{\beta-1}{\beta}}{\color{purple}L^{\frac{2}{\beta}}}\alpha^{-1}\bigg),
    \\
    \label{eq2:thm:sconvex_cons}
    \Exp[\norm{\hat{\bx}_T-\bx^\star}^2] &\leq \min\bigg(B^2, \frac{2GB}{\alpha}, \frac{8{\color{purple}\bL^2}{\sf V}_1 G^2}{\alpha^2 T}+2\left(\cst_1+\cst_2{\color{purple}\frac{\bL^2}{\sigma^2}}\right)\left(\frac{d^2{\color{purple}\sigma^2}}{T}\right)^{\frac{\beta-1}{\beta}}{\color{purple}L^{\frac{2}{\beta}}}\alpha^{-2}\bigg),
\end{align}
 where $B$ is the Euclidean diameter of $\com$, and $\cst_1,\cst_2> 0$ depends only on $\beta$ and on the choice of the gradient estimator.
\end{corollary}
In a similar setting, but assuming independent zero-mean $\xi_t$'s,  \cite{BP2016} considered the case of $\ell_2$ randomization and proved, for integer $\beta>2$, a bound with suboptimal rate $T^{-\frac{\beta-1}{\beta+1}}$.
Corollary~\ref{thm:sconvex_cons} can be also compared to~\cite{akhavan2020,akhavan2021distributed} ($\ell_2$ randomization and coordinatewise radomization) and, for $\beta>2$, to~\cite{Gasnikov} ($\ell_2$ randomization). However, those papers use a slightly different definition of $\beta$-smoothness class (both definitions coincide if $\beta = 2$, see Remark~\ref{rem:smoothness}). Their bounds guarantee the rate $O\big(\tfrac{d^{2-1/\beta}}{\alpha T}\big)$ for $\beta>2$~\cite[Corollary 6]{akhavan2021distributed}, \cite[Theorem 1]{Gasnikov} and $O\big(\tfrac{d}{\sqrt{\alpha T}}\big)$ for $\beta=2$~\cite[Theorem~D.4]{akhavan2020} by using two different approaches for the two cases. In contrast, Corollary~\ref{thm:sconvex_cons} yields $O\big(\tfrac{d^{2-2/\beta}}{\alpha T}\big)$ and $O\big(\tfrac{d}{\alpha\sqrt{ T}}\big)$, respectively, and obtains these rates by a unified approach for all $\beta\ge2$, and simultaneously under $\ell_1$ and $\ell_2$ randomizations.  Note that, under the condition $T\ge d^2$  guaranteeing non-triviality of the bound, and $\alpha\ge 1$ the rate $O\big(\tfrac{d}{\alpha\sqrt{ T}}\big)$ that we obtain in Corollary~\ref{thm:sconvex_cons} for $\beta=2$  matches the minimax lower bound (cf. Theorem \ref{lb} below)
as a function of all the three parameters $T, d$, and $\alpha$. 

% \newpage

%%%%%%%%%%%%%%%%%%%%%%

\section{Lower bounds}
\label{sec:LB}
{\color{purple}In this section we prove minimax lower bounds on the optimization error over all sequential 
strategies with two-point feedback that allow the query points depend on the past.
For $t=1,\dots, T$, we assume that the values $y_t=f(\bz_t)+\xi_t$ and $y'_t=f(\bz'_t)+\xi'_t$ are observed, where $(\xi,\xi'_t)$ are random noises, and 
$(\bz_t, \bz'_t)$ are query points.
We consider all strategies of choosing the query points as 
$\bz_t= \Phi_t\big(\left(\bz_i,y_i\right)_{i=1}^{t-1},\left(\bz'_i,y'_i\right)_{i=1}^{t-1},\btau_t\big)$ and $\bz'_t= \Phi'_t\big(\left(\bz_i,y_i\right)_{i=1}^{t-1},\left(\bz'_i,y'_i\right)_{i=1}^{t-1},\btau_t\big)$ for $t\ge 2$, where $\Phi_t$'s and $\Phi'_t$'s are measurable functions, $\bz_1, \bz'_1 \in \mathbb{R}^d$
are any random variables, and $\{\btau_t\}$ is a sequence of random variables with values in a measurable space $(\mathcal Z, \mathcal U)$, such that $\btau_t$ is independent of $\big(\left(\bz_i,y_i\right)_{i=1}^{t-1},\left(\bz'_i,y'_i\right)_{i=1}^{t-1}\big)$. We denote by $\Pi_T$ the set of all such strategies of choosing query points up to $t=T$.
The class $\Pi_T$ includes the sequential strategy of Algorithm~\eqref{eq:algo_general} with either of the two 
considered gradient estimators~\eqref{eq:grad_l2} and~\eqref{eq:grad_l1}. In this case, $\btau_t=(\bzeta_t,r_t)$, $\bz_t=\bx_t+h_t\bzeta_tr_t$ and $\bz'_t=\bx_t-h_t\bzeta_tr_t$, where $\bzeta_t=\bzeta_t^{\circ}$ or $\bzeta_t=\bzeta_t^{\diamond}$.

To state our assumption on the noises $(\xi,\xi'_t)$, we introduce the squared Hellinger distance $H^2(\cdot, \cdot)$ defined for two probability measures $\mathbf{P}, \mathbf{P}'$ on a measurable space $(\Omega, \mathcal{A})$ as} 
\[
H^2(\mathbf{P}, \mathbf{P}') \triangleq \int (\sqrt{\d \mathbf{P}} - \sqrt{\d \mathbf{P}'})^2\enspace.
\]

\begin{assumption}
\label{ass:lower-bound}
For every $t\ge 1$, the following holds:
\begin{itemize}
    \item The cumulative distribution function $F_t : \bbR^2 \to \bbR$ of random variable $(\xi_t,\xi'_t)$ 
is such that
\begin{equation}
\label{distribution}
H^2 (P_{F_t(\cdot,\cdot)}, P_{F_t(\cdot+v,\cdot+v)})\leq I_{0}v^2\,, \quad\quad |v|\leq v_{0}\enspace,
\end{equation}
for some $0<I_{0}<\infty$, $0<v_{0}\leq \infty$. Here, $P_{F(\cdot,\cdot)}$ denotes the probability measure corresponding to the c.d.f. $F(\cdot,\cdot)$.
\item The random variable $(\xi_t,\xi'_t)$ is independent of $((\bz_i,y_i)_{i=1}^{t-1},(\bz'_i,y'_i)_{i=1}^{t-1},\btau_t)$. 
\end{itemize}
\end{assumption}

Condition~\eqref{distribution} is not restrictive and encompasses a large family of distributions. It is satisfied with small enough $v_0$  for distributions that correspond to regular statistical experiments, \cite[see \eg][Chapter 1]{Ibra}. If $F_t$ is a Gaussian c.d.f. condition~\eqref{distribution}  
is satisfied with $v_0=\infty$. 

To state the lower bounds, we consider a subset of the classes of functions, for which we obtained the upper bounds in Section~\ref{sec4}. Let $\com=\{\bx\in\mathbb{R}^d\,:\, \norm{\bx}\le 1\}$.  For $\alpha, L, \bar L>0$, $ \beta\ge 2$,  let  $\mathcal{F}_{\alpha,\beta}$ denote the set of all $\alpha$-strongly convex functions $f$ that satisfy Assumption~\ref{nat}, attain their minimum over $\mathbb{R}^d$ in $\com$ and such that $\max_{\bx\in\com}\|\nabla f(\bx)\|\le G$,  and the condition $G> \alpha$ is satisfied
\footnote{The condition $G\ge\alpha$ is necessary for the class $\mathcal{F}_{\alpha,\beta}$ to be non-empty. Indeed, due to \eqref{eq:def-PL} and \eqref{eq:opt-err-to-est-err}, for all $f\in\mathcal{F}_{\alpha,\beta}$ and $x\in \com$ we have $\norm{\bx - \bx^{*}}\le G/\alpha$, and thus $2G/\alpha \ge {\rm diam}(\com)=2$.}. 

\begin{restatable}{theorem}{lowerB}
\label{lb}
Let $\com=\{\bx\in\mathbb{R}^d\,:\, \norm{\bx}\le 1\}$ and let Assumption~\ref{ass:lower-bound} hold.
{\color{brown}Then, for any estimator $\tilde\bx_T$ based on the observations $((\bz_t,y_t), (\bz'_t,y'_t), t=1,\dots, T)$, where $((\bz_t,\bz'_t), t=1,\dots, T)$ are obtained by any strategy in the class $\Pi_T$ we have}
\begin{equation}\label{eq1:lb}
\sup_{f \in \mathcal{F}_{\alpha,\beta}}\Exp\big[f(\tilde\bx_T)-f^{\star}\big]\geq C\min\left(\max\left(\alpha,\, T^{-1/2+1/\beta}\right),\, \frac{d}{\sqrt{T}}, \,\frac{d}{\alpha}T^{-\frac{\beta-1}{\beta}}\right)\enspace,
\end{equation}
and
\begin{equation}\label{eq2:lb}
\sup_{f \in \mathcal{F}_{\alpha,\beta}}\Exp[\norm{\bz_{T}-\bx^{*}(f)}^{2}]\geq C\min\left(1,\, \frac{d}{T^{\frac{1}{\beta}}}, \,\frac{d}{\alpha^{2}}T^{-\frac{\beta-1}{\beta}}\right)\enspace,
\end{equation}
 where  $C>0$ is a constant that does not depend  of $T,d$, and $\alpha$, and $\bx^\star(f)$ is the minimizer of $f$ on~$\com$. 
\end{restatable}

Some remarks are in order here. First, note that the threshold $T^{-1/2+1/\beta}$ on the strong convexity parameter $\alpha$ plays an important role in bounds \eqref{eq1:lb} and \eqref{eq2:lb}. Indeed, for $\alpha$ below this threshold, the bounds start to be independent of $\alpha$.  Intuitively, it seems reasonable that $\alpha$-strong convexity should be of no added value for small $\alpha$. Theorem \ref{lb} allows us to quantify exactly how small such $\alpha$ should be, namely, $\alpha \lesssim T^{-1/2+1/\beta}$. In particular, for $\beta=2$ the threshold occurs at $\alpha\asymp 1$. Also, quite naturally, the threshold becomes smaller when the smoothness $\beta$ increases.

In the regime below the $T^{-1/2+1/\beta}$ threshold, the rate of \eqref{eq1:lb} becomes $\min(T^{1/\beta},d)/\sqrt{T}$, which is asymptotically  $d/\sqrt{T}$ independently of the smoothness index $\beta$ and on $\alpha$. Thus, we obtain that $d/\sqrt{T}$
is a lower bound over the class of simply convex functions. On the other hand, the achievable rate for convex functions  is shown to be $d^{16}/\sqrt{T}$ in \cite{agarwal2011} and improved to $d^{4.5}/\sqrt{T}$ 
in  \cite{lattimore-gyorgy2021} 
(both results are up to poly-logarithmic factors, and under sub-Gaussian noise $\xi_t$). 
The gap between our lower bound $d/\sqrt{T}$ and these upper bounds is only in the dependence on the dimension, but this gap is substantial. 
In the regime where $\alpha$ is above the $T^{-1/2+1/\beta}$ threshold, our results imply that the gap between upper and lower bounds is much smaller. Thus, our upper bounds in this regime scale as $\frac{d^{2-2/\beta}}{\alpha} T^{-\frac{\beta-1}{\beta}}$ while the lower bound of 
Theorem~\ref{lb} is of the order $\frac{d}{\alpha}T^{-\frac{\beta-1}{\beta}}$.

Consider now the case $\beta=2$. Then the lower bounds \eqref{eq1:lb} and \eqref{eq2:lb} are of order $d/(\max(\alpha, 1)\sqrt{T})$ and $d/(\max(\alpha^2, 1)\sqrt{T})$, respectively, under the condition $T\ge d^2$  guaranteeing non-triviality of the rates. If, in addition, $\alpha\gtrsim 1$ (meaning that $\alpha$ is above the threshold $\alpha\asymp 1$) we obtain 
the lower rates $d/(\alpha\sqrt{T})$ and $d/(\alpha^2\sqrt{T})$, respectively. Comparing this remark to Corollary~\ref{thm:sconvex_cons} we obtain the following result. 
\begin{corollary}\label{cor:minimax rate for beta=2}
Let $\beta=2$ and let the assumptions of Theorem~\ref{lb}  and Corollary~\ref{thm:sconvex_cons} hold. If $\alpha\ge 1$ and $T\ge \max({\color{brown}d^2,}{\color{purple}\bL^2}d,{\color{purple}\bL^4}G^4)$ then there exist positive constants $c,C$ that do not depend on $T,d$, and $\alpha$ such that we have the following bounds on the minimax risks:
\begin{equation}\label{eq1:cor:lb}
c\frac{d}{\alpha\sqrt{T}} \le \inf_{\tilde\bx_T} \sup_{f \in \mathcal{F}_{\alpha,\beta}}\Exp\big[f(\tilde\bx_T)-f^{\star}\big] \le C\frac{d}{\alpha\sqrt{T}},
\end{equation}
and
\begin{equation}\label{eq2:cor:lb}
c\frac{d}{\alpha^2\sqrt{T}} \le \inf_{\tilde\bx_T} \sup_{f \in \mathcal{F}_{\alpha,\beta}}\Exp[\norm{\tilde\bx_T-\bx^{*}(f)}^{2}]\le C \frac{d}{\alpha^2\sqrt{T}},
\end{equation}
 where $\bx^\star(f)$ is the minimizer of $f$ on~$\com$, and  
 %$\displaystyle{\inf_{\tilde\bx_T}}$ denotes 
 the infimum {\color{brown} is over all estimators $\tilde\bx_T$ based on query points obtained via strategies in the class~$\Pi_T$. The minimax rates in \eqref{eq1:cor:lb} and \eqref{eq2:cor:lb} are attained by the 
 estimator $\tilde\bx_T=\hat\bx_{T}$ with parameters as in Corollary~\ref{thm:sconvex_cons}.}
\end{corollary}
Thus, the weighted average estimator $\hat\bx_{T}$  as in Corollary~\ref{thm:sconvex_cons} is minimax optimal with respect to all the three parameters $T, d$, and $\alpha$, both in the optimization error and in the estimation risk. Note that we introduced the condition $T\ge \max({\color{purple}\bL^2}d,{\color{purple}\bL^4}G^4)$ in Corollary~\ref{cor:minimax rate for beta=2} to guarantee that the upper bounds are of the required order, cf. Corollary~\ref{thm:sconvex_cons}. Thus, $G$ is allowed to be a function of $T$ that grows not too fast. Since $G>\alpha$ this condition also prevents $\alpha$ from being too large, that is,  Corollary~\ref{cor:minimax rate for beta=2} does not hold if $\alpha\gtrsim T^{1/4}$.
 
%Next, note that since any $\alpha$-strongly convex function is $\alpha$-gradient dominant, inequality~\eqref{eq1:lb} provides a lower bound for the family of $\alpha$-gradient dominant functions satisfying the assumptions of Corollary~\ref{cor:finalres_PL}. This fact implies the minimax optimality, in terms of $T$, of our algorithms (with parameters as in Corollary~\ref{cor:finalres_PL}) on the class of $\alpha$-gradient dominant $\beta$-Hölder functions. Again, for the case $\beta=2$, under the natural assumptions $\alpha\gtrsim 1$, and $T\ge d^2$ we obtain that the rate
%$d/(\alpha\sqrt{T})$ is minimax optimal as a function of the triplet ($T, d$, $\alpha$) on the class
%of $\alpha$-gradient dominant 2-Hölder functions, and it is attained by our algorithms with parameters as in Corollary~\ref{cor:finalres_PL}. {\color{red}(Sasha: this should be stated as a corollary for $\alpha$-gradient dominant, or included in the corollary for $\alpha$-strongly convex)}

The issue of finding the minimax optimal rates in gradient-free stochastic optimization under strong convexity and smoothness assumptions has a long history. It was initiated in \cite{Fabian,PT90} and more recently developed in \cite{dippon,jamieson2012,Shamir13,BP2016,akhavan2020,akhavan2021distributed}. 
%\sasha{some words about importance of strong convexity?} 
It was shown in \cite{PT90} that the minimax optimal rate on the class of $\alpha$-strong convex and $\beta$-Hölder functions scales as $c(\alpha,d)T^{-{(\beta-1)}/{\beta}}$ for $\beta\ge 2$, where $c(\alpha,d)$ is an unspecified function of $\alpha$ and $d$ (for $d=1$ and integer $\beta\ge 2$ an upper bound of the same order was earlier derived in 
\cite{Fabian}). The issue of establishing non-asymptotic fundamental limits as function of the main parameters of the problem ($\alpha$, $d$ and $T$) was first addressed in \cite{jamieson2012} giving a lower bound $\Omega(\sqrt{d/T})$ for $\beta=2$, without specifying the dependency on $\alpha$. This was improved to $\Omega(d/\sqrt{T})$ when $\alpha\asymp 1$ by \cite{Shamir13} who also claimed that the rate $d/\sqrt{T}$ is optimal for $\beta=2$ 
referring to an upper bound in \cite{agarwal2010}. However, invoking \cite{agarwal2010} in the setting with random noise $\xi_t$ (for which the lower bound of \cite{Shamir13} was proved) is not legitimate because in \cite{agarwal2010} the observations are considered as a Lipschitz function of $t$. The complete proof of minimax optimality of the rate $d/\sqrt{T}$ for $\beta=2$ under random noise was later provided in \cite{akhavan2020}. However, the upper and the lower bounds in \cite{akhavan2020} still differ in their dependence on $\alpha$. Corollary~\ref{cor:minimax rate for beta=2} completes this line of work by establishing the minimax optimality as a function of the whole triplet ($T, d$, $\alpha$) for $\beta=2$. %Moreover, we get analogous conclusion not only for $\alpha$-strongly convex functions but also for $\alpha$-gradient dominant functions that were not covered by the previous work.  

The main lines of the proof of Theorem~\ref{lb} follow~\cite{akhavan2020}. However, in~\cite{akhavan2020} the assumptions on the noise are much more restrictive -- the random variables $\xi_t$ are assumed iid and instead of~\eqref{distribution} a much stronger condition is imposed, namely, a bound on the Kullback–Leibler divergence. 
In particular, in order to use the Kullback–Leibler divergence between two distributions we need one of them to be absolutely continuous with respect to the other. Using the Hellinger distance allows us to drop this restriction. For example, if $F_t$ is a distribution with bounded support then the Kullback–Leibler divergence between $F_t(\cdot)$ and $F_t(\cdot+v)$ is $+\infty$ while the Hellinger distance is finite.

\bibliography{biblio}

\begin{thebibliography}{28}
\providecommand{\natexlab}[1]{#1}
\providecommand{\url}[1]{\texttt{#1}}
\expandafter\ifx\csname urlstyle\endcsname\relax
  \providecommand{\doi}[1]{doi: #1}\else
  \providecommand{\doi}{doi: \begingroup \urlstyle{rm}\Url}\fi

\bibitem[Akhavan et~al.(2020)Akhavan, Pontil, and Tsybakov]{akhavan2020}
A.~Akhavan, M.~Pontil, and A.B. Tsybakov.
\newblock Exploiting higher order smoothness in derivative-free optimization
  and continuous bandits.
\newblock In \emph{Advances in Neural Information Processing Systems 33}, 2020.

\bibitem[Akhavan et~al.(2022)Akhavan, Chzhen, Pontil, and
  Tsybakov]{Akhavan_Chzhen_Pontil_Tsybakov22a}
A.~Akhavan, E.~Chzhen, M.~Pontil, and A.B. Tsybakov.
\newblock A gradient estimator via l1-randomization for online zero-order
  optimization with two point feedback, 2022.

\bibitem[Akhavan et~al.(2021)Akhavan, Pontil, and
  Tsybakov]{akhavan2021distributed}
Arya Akhavan, Massimiliano Pontil, and Alexandre~B Tsybakov.
\newblock Distributed zero-order optimization under adversarial noise.
\newblock \emph{arXiv preprint arXiv:2102.01121}, 2021.

\bibitem[Arjevani et~al.(2019)Arjevani, Carmon, Duchi, Foster, Srebro, and
  Woodworth]{azjev}
Yossi Arjevani, Yair Carmon, John Duchi, Dylan Foster, Nathan Srebro, and Blake
  Woodworth.
\newblock Lower bounds for non-convex stochastic optimization.
\newblock 12 2019.

\bibitem[Bach and Perchet(2016)]{BP2016}
F.~Bach and V.~Perchet.
\newblock Highly-smooth zero-th order online optimization.
\newblock In \emph{Proc. 29th Annual Conference on Learning Theory}, pages
  1--27, 2016.

\bibitem[Balashov et~al.(2020)Balashov, Polyak, and
  Tremba]{balashov2020gradient}
MV~Balashov, BT~Polyak, and AA~Tremba.
\newblock Gradient projection and conditional gradient methods for constrained
  nonconvex minimization.
\newblock \emph{Numerical Functional Analysis and Optimization}, 41\penalty0
  (7):\penalty0 822--849, 2020.

\bibitem[Balasubramanian and Ghadimi(2021)]{balasubramanian2021zeroth}
Krishnakumar Balasubramanian and Saeed Ghadimi.
\newblock Zeroth-order nonconvex stochastic optimization: Handling constraints,
  high dimensionality, and saddle points.
\newblock \emph{Foundations of Computational Mathematics}, pages 1--42, 2021.

\bibitem[Barthe et~al.(2005)Barthe, Gu{\'e}don, Mendelson, and
  Naor]{Barthe_Guedon_Mendelson_Naor05}
Franck Barthe, Olivier Gu{\'e}don, Shahar Mendelson, and Assaf Naor.
\newblock {A probabilistic approach to the geometry of the Lpn-ball}.
\newblock \emph{The Annals of Probability}, 33\penalty0 (2):\penalty0 480 --
  513, 2005.

\bibitem[Beckner(1989)]{Beckner89}
William Beckner.
\newblock A generalized poincaré inequality for gaussian measures.
\newblock \emph{Proceedings of the American Mathematical Society}, 105\penalty0
  (2):\penalty0 397--400, 1989.

\bibitem[Bernstein et~al.(2018)Bernstein, Wang, Azizzadenesheli, and
  Anandkumar]{bernstein2018signsgd}
Jeremy Bernstein, Yu-Xiang Wang, Kamyar Azizzadenesheli, and Animashree
  Anandkumar.
\newblock signsgd: Compressed optimisation for non-convex problems.
\newblock In \emph{International Conference on Machine Learning}, pages
  560--569. PMLR, 2018.

\bibitem[Bubeck(2015)]{bubeck2014}
S.~Bubeck.
\newblock Convex optimization: Algorithms and complexity.
\newblock \emph{Foundations and Trends in Machine Learning}, 8(3-4):231--257,
  2015.

\bibitem[Carmon et~al.(2017)Carmon, Duchi, Hinder, and Sidford]{duchi}
Yair Carmon, John Duchi, Oliver Hinder, and Aaron Sidford.
\newblock Lower bounds for finding stationary points i.
\newblock \emph{Mathematical Programming}, 184, 10 2017.

\bibitem[Feng(2010)]{fengqi}
Qi~Feng.
\newblock \emph{Bounds for the ratio of two gamma functions}.
\newblock Journal of Inequalities and Applications, 2010.

\bibitem[Flaxman et~al.(2005)Flaxman, Kalai, and McMahan]{flaxman2004}
A.~D. Flaxman, A.~T. Kalai, and H.~B. McMahan.
\newblock Online convex optimization in the bandit setting: gradient descent
  without a gradient.
\newblock In \emph{Proc. 16th Annual ACM-SIAM Symposium on Discrete algorithms
  (SODA)}, pages 385–--394, 2005.

\bibitem[Garrigos et~al.(2017)Garrigos, Rosasco, and Villa]{lr}
Guillaume Garrigos, Lorenzo Rosasco, and Silvia Villa.
\newblock Convergence of the forward-backward algorithm: Beyond the worst case
  with the help of geometry.
\newblock 03 2017.

\bibitem[Ghadimi and Lan(2013)]{Ghadimi2013}
S.~Ghadimi and G.~Lan.
\newblock Stochastic first- and zeroth-order methods for nonconvex stochastic
  programming.
\newblock \emph{SIAM Journal on Optimization}, 23(4):\penalty0 2341–2368,
  2013.

\bibitem[Karimi et~al.(2016)Karimi, Nutini, and
  Schmidt]{10.1007/978-3-319-46128-1_50}
H.~Karimi, J.~Nutini, and M.~Schmidt.
\newblock Linear convergence of gradient and proximal-gradient methods under
  the polyak-{\l}ojasiewicz condition.
\newblock In \emph{Machine Learning and Knowledge Discovery in Databases},
  pages 795--811, 2016.

\bibitem[Nemirovski(2000)]{nemirovski2000topics}
Arkadi Nemirovski.
\newblock Topics in non-parametric statistics.
\newblock \emph{Ecole d’Et{\'e} de Probabilit{\'e}s de Saint-Flour},
  28:\penalty0 85, 2000.

\bibitem[Nemirovsky and Yudin(1983)]{NY1983}
A.~S. Nemirovsky and D.~B Yudin.
\newblock \emph{Problem Complexity and Method Efficiency in Optimization.}
\newblock Wiley \& Sons, 1983.

\bibitem[Novitskii and Gasnikov(2021)]{Gasnikov}
V.~Novitskii and A.~Gasnikov.
\newblock Improved exploiting higher order smoothness in derivative-free
  optimization and continuous bandit.
\newblock \emph{arXiv preprint arXiv:2101.03821}, 2021.

\bibitem[Osserman(1978)]{Osserman78}
R.~Osserman.
\newblock {The isoperimetric inequality}.
\newblock \emph{Bulletin of the American Mathematical Society}, 84\penalty0
  (6):\penalty0 1182 -- 1238, 1978.

\bibitem[Polyak(1963)]{PB63}
Boris Polyak.
\newblock Gradient methods for the minimisation of functionals.
\newblock \emph{Ussr Computational Mathematics and Mathematical Physics},
  3:\penalty0 864--878, 12 1963.

\bibitem[Polyak and Tsybakov(1990)]{PT90}
T.~B. Polyak and A.~B. Tsybakov.
\newblock Optimal order of accuracy of search algorithms in stochastic
  optimization.
\newblock \emph{Problems of Information Transmission}, 26\penalty0
  (2):\penalty0 45--53, 1990.

\bibitem[Rachev and Ruschendorf(1991)]{Rachev_Ruschendorf91}
S.~T. Rachev and L.~Ruschendorf.
\newblock {Approximate Independence of Distributions on Spheres and Their
  Stability Properties}.
\newblock \emph{The Annals of Probability}, 19\penalty0 (3):\penalty0 1311 --
  1337, 1991.

\bibitem[Schechtman and Zinn(1990)]{Schechtman_Zinn90}
G.~Schechtman and J.~Zinn.
\newblock On the volume of the intersection of two {$L^n_p$} balls.
\newblock \emph{Proc. Amer. Math. Soc.}, 110\penalty0 (1):\penalty0 217--224,
  1990.

\bibitem[Shamir(2017)]{Shamir17}
O.~Shamir.
\newblock An optimal algorithm for bandit and zero-order convex optimization
  with two-point feedback.
\newblock \emph{Journal of Machine Learning Research}, 18\penalty0
  (1):\penalty0 1703--1713, 2017.

\bibitem[Tsybakov(2009)]{Tsybakov09}
A.~Tsybakov.
\newblock \emph{Introduction to Nonparametric Estimation}.
\newblock Springer, New York, 2009.

\bibitem[Zorich(2016)]{zorich2016}
Vladimir~Antonovich Zorich.
\newblock \emph{Mathematical analysis II}.
\newblock Springer, 2016.

\end{thebibliography}


\begin{thebibliography}{40}
\providecommand{\natexlab}[1]{#1}
\providecommand{\url}[1]{\texttt{#1}}
\expandafter\ifx\csname urlstyle\endcsname\relax
  \providecommand{\doi}[1]{doi: #1}\else
  \providecommand{\doi}{doi: \begingroup \urlstyle{rm}\Url}\fi

\bibitem[Agarwal et~al.(2010)Agarwal, Dekel, and Xiao]{agarwal2010}
A.~Agarwal, O.~Dekel, and L.~Xiao.
\newblock Optimal algorithms for online convex optimization with multi-point
  bandit feedback.
\newblock In \emph{Proc. 23rd International Conference on Learning Theory},
  pages 28--40, 2010.

\bibitem[Agarwal et~al.(2011)Agarwal, Foster, Hsu, Kakade, and
  Rakhlin]{agarwal2011}
A.~Agarwal, D.~P. Foster, D.~J. Hsu, S.~M. Kakade, and A.~Rakhlin.
\newblock Stochastic convex optimization with bandit feedback.
\newblock In \emph{Advances in Neural Information Processing Systems},
  volume~25, pages 1035--1043, 2011.

\bibitem[Akhavan et~al.(2020)Akhavan, Pontil, and Tsybakov]{akhavan2020}
A.~Akhavan, M.~Pontil, and A.B. Tsybakov.
\newblock Exploiting higher order smoothness in derivative-free optimization
  and continuous bandits.
\newblock In \emph{Advances in Neural Information Processing Systems 33}, 2020.

\bibitem[Akhavan et~al.(2021)Akhavan, Pontil, and
  Tsybakov]{akhavan2021distributed}
A.~Akhavan, M.~Pontil, and A.~B. Tsybakov.
\newblock Distributed zero-order optimization under adversarial noise.
\newblock In \emph{Advances in Neural Information Processing Systems 34}, 2021.

\bibitem[Akhavan et~al.(2022)Akhavan, Chzhen, Pontil, and
  Tsybakov]{Akhavan_Chzhen_Pontil_Tsybakov22a}
A.~Akhavan, E.~Chzhen, M.~Pontil, and A.B. Tsybakov.
\newblock A gradient estimator via l1-randomization for online zero-order
  optimization with two point feedback.
\newblock In \emph{Advances in Neural Information Processing Systems 35}, 2022.

\bibitem[Arjevani et~al.(2022)Arjevani, Carmon, Duchi, Foster, Srebro, and
  Woodworth]{azjev}
Y.~Arjevani, Y.~Carmon, J.~Duchi, D.~Foster, N.~Srebro, and B.~Woodworth.
\newblock Lower bounds for non-convex stochastic optimization.
\newblock \emph{Mathematical Programming}, 2022.

\bibitem[Bach and Perchet(2016)]{BP2016}
F.~Bach and V.~Perchet.
\newblock Highly-smooth zero-th order online optimization.
\newblock In \emph{Proc. 29th Annual Conference on Learning Theory}, 2016.

\bibitem[Balashov et~al.(2020)Balashov, Polyak, and
  Tremba]{balashov2020gradient}
M.~V. Balashov, B.~T. Polyak, and A.~A. Tremba.
\newblock Gradient projection and conditional gradient methods for constrained
  nonconvex minimization.
\newblock \emph{Numerical Functional Analysis and Optimization}, 41\penalty0
  (7):\penalty0 822--849, 2020.

\bibitem[Balasubramanian and Ghadimi(2021)]{balasubramanian2021zeroth}
K.~Balasubramanian and S.~Ghadimi.
\newblock Zeroth-order nonconvex stochastic optimization: Handling constraints,
  high dimensionality, and saddle points.
\newblock \emph{Foundations of Computational Mathematics}, pages 1--42, 2021.

\bibitem[Barthe et~al.(2005)Barthe, Gu{\'e}don, Mendelson, and
  Naor]{Barthe_Guedon_Mendelson_Naor05}
F.~Barthe, O.~Gu{\'e}don, S.~Mendelson, and A.~Naor.
\newblock {A probabilistic approach to the geometry of the {$L_p^n$} ball}.
\newblock \emph{The Annals of Probability}, 33\penalty0 (2):\penalty0 480--513,
  2005.

\bibitem[Bubeck(2015)]{bubeck2014}
S.~Bubeck.
\newblock Convex optimization: Algorithms and complexity.
\newblock \emph{Foundations and Trends in Machine Learning}, 8(3-4):231--257,
  2015.

\bibitem[Carmon et~al.(2020)Carmon, Duchi, Hinder, and Sidford]{duchi}
Y.~Carmon, J.~Duchi, O~Hinder, and Aaron Sidford.
\newblock Lower bounds for finding stationary points {I}.
\newblock \emph{Mathematical Programming}, 184:\penalty0 71--120, 2020.

\bibitem[Dippon(2003)]{dippon}
J.~Dippon.
\newblock Accelerated randomized stochastic optimization.
\newblock \emph{Ann. Statist.}, 31\penalty0 (4):\penalty0 1260--1281, 2003.

\bibitem[Duchi et~al.(2015)Duchi, Jordan, Wainwright, and Wibisono]{duchi2015}
J.~C. Duchi, M.~I. Jordan, M.~J. Wainwright, and A.~Wibisono.
\newblock Optimal rates for zero-order convex optimization: The power of two
  function evaluations.
\newblock \emph{IEEE Transactions on Information Theory}, 61\penalty0
  (5):\penalty0 2788--2806, 2015.

\bibitem[Fabian(1967)]{Fabian}
V.~Fabian.
\newblock Stochastic approximation of minima with improved asymptotic speed.
\newblock \emph{The Annals of Mathematical Statistics}, 38\penalty0
  (1):\penalty0 191--200, 1967.

\bibitem[Flaxman et~al.(2005)Flaxman, Kalai, and McMahan]{flaxman2004}
A.~D. Flaxman, A.~T. Kalai, and H.~B. McMahan.
\newblock Online convex optimization in the bandit setting: gradient descent
  without a gradient.
\newblock In \emph{Proc. 16th Annual ACM-SIAM Symposium on Discrete algorithms
  (SODA)}, 2005.

\bibitem[Garrigos et~al.(2023)Garrigos, Rosasco, and Villa]{lr}
G.~Garrigos, L.~Rosasco, and S.~Villa.
\newblock Convergence of the forward-backward algorithm: Beyond the worst case
  with the help of geometry.
\newblock \emph{Mathematical Programming}, 198:\penalty0 937--996, 2023.

\bibitem[Gasnikov et~al.(2016)Gasnikov, Lagunovskaya, Usmanova, and
  Fedorenko]{gasnikov_gradient-free_2016}
A.~Gasnikov, A.~Lagunovskaya, I.~Usmanova, and F.~Fedorenko.
\newblock Gradient-free proximal methods with inexact oracle for convex
  stochastic nonsmooth optimization problems on the simplex.
\newblock \emph{Automation and Remote Control}, 77\penalty0 (11):\penalty0
  2018--2034, 2016.

\bibitem[Ghadimi and Lan(2013)]{Ghadimi2013}
S.~Ghadimi and G.~Lan.
\newblock Stochastic first- and zeroth-order methods for nonconvex stochastic
  programming.
\newblock \emph{SIAM Journal on Optimization}, 23(4):\penalty0 2341--2368,
  2013.

\bibitem[Ibragimov and Khas'minskii(1982)]{Ibra}
I.~A. Ibragimov and R.~Z. Khas'minskii.
\newblock Estimation of the maximum value of a signal in gaussian white noise.
\newblock \emph{Mat. Zametki}, 32\penalty0 (4):\penalty0 746--750, 1982.

\bibitem[Jamieson et~al.(2012)Jamieson, Nowak, and Recht]{jamieson2012}
K.~G. Jamieson, R.~Nowak, and B.~Recht.
\newblock Query complexity of derivative-free optimization.
\newblock In \emph{Advances in Neural Information Processing Systems},
  volume~26, pages 2672--2680, 2012.

\bibitem[Karimi et~al.(2016)Karimi, Nutini, and
  Schmidt]{10.1007/978-3-319-46128-1_50}
H.~Karimi, J.~Nutini, and M.~Schmidt.
\newblock Linear convergence of gradient and proximal-gradient methods under
  the {P}olyak-{Ł}ojasiewicz condition.
\newblock In \emph{Machine Learning and Knowledge Discovery in Databases},
  2016.

\bibitem[Lattimore and György(2021)]{lattimore-gyorgy2021}
T.~Lattimore and A.~György.
\newblock Improved regret for zeroth-order stochastic convex bandits.
\newblock In \emph{Advances in Neural Information Processing Systems 34}, 2021.

\bibitem[Nemirovski(2000)]{nemirovski2000topics}
A.~Nemirovski.
\newblock Topics in non-parametric statistics.
\newblock \emph{Ecole d’Et{\'e} de Probabilit{\'e}s de Saint-Flour 28}, 2000.

\bibitem[Nemirovsky and Yudin(1983)]{NY1983}
A.~S. Nemirovsky and D.~B Yudin.
\newblock \emph{Problem Complexity and Method Efficiency in Optimization.}
\newblock Wiley \& Sons, 1983.

\bibitem[Nesterov(2011)]{Nesterov2011}
Y.~Nesterov.
\newblock Random gradient-free minimization of convex functions.
\newblock Technical Report 2011001, Center for Operations Research and
  Econometrics (CORE), Catholic University of Louvain, 2011.

\bibitem[Nesterov(2018)]{nesterov-book}
Y.~Nesterov.
\newblock \emph{Lectures on Convex Optimization}.
\newblock Springer, 2018.

\bibitem[Nesterov and Spokoiny(2017)]{NS17}
Y.~Nesterov and V.~Spokoiny.
\newblock Random gradient-free minimization of convex functions.
\newblock \emph{Found. Comput. Math.}, 17:\penalty0 527--–566, 2017.

\bibitem[Novitskii and Gasnikov(2022)]{Gasnikov}
V.~Novitskii and A.~Gasnikov.
\newblock Improved exploitation of higher order smoothness in derivative-free
  optimization.
\newblock \emph{Optimization Letters}, 16:\penalty0 2059--2071, 2022.

\bibitem[Osserman(1978)]{Osserman78}
R.~Osserman.
\newblock {The isoperimetric inequality}.
\newblock \emph{Bulletin of the American Mathematical Society}, 84\penalty0
  (6):\penalty0 1182--1238, 1978.

\bibitem[Polyak and Tsybakov(1990)]{PT90}
B.~T. Polyak and A.~B. Tsybakov.
\newblock Optimal order of accuracy of search algorithms in stochastic
  optimization.
\newblock \emph{Problems of Information Transmission}, 26\penalty0
  (2):\penalty0 45--53, 1990.

\bibitem[Polyak(1963)]{PB63}
B.T. Polyak.
\newblock Gradient methods for the minimisation of functionals.
\newblock \emph{USSR Computational Mathematics and Mathematical Physics},
  3:\penalty0 864--878, 1963.

\bibitem[Qi and Luo(2013)]{fengqi}
F.~Qi and Q.-M. Luo.
\newblock {Bounds for the ratio of two gamma functions: from Wendel’s
  asymptotic relation to Elezovi{ć}-Giordano-Pe{č}ari{ć}’s theorem}.
\newblock \emph{Journal of Inequalities and Applications}, 2013.

\bibitem[Rachev and Ruschendorf(1991)]{Rachev_Ruschendorf91}
S.~T. Rachev and L.~Ruschendorf.
\newblock Approximate independence of distributions on spheres and their
  stability properties.
\newblock \emph{The Annals of Probability}, 19\penalty0 (3):\penalty0 1311 --
  1337, 1991.

\bibitem[Rando et~al.(2022)Rando, Molinari, Villa, and
  Rosasco]{Rando_Molinari_Villa_Rosasco22}
M.~Rando, C.~Molinari, S.~Villa, and L.~Rosasco.
\newblock Stochastic zeroth order descent with structured directions.
\newblock \emph{arXiv:2206.05124}, 2022.

\bibitem[Schechtman and Zinn(1990)]{Schechtman_Zinn90}
G.~Schechtman and J.~Zinn.
\newblock On the volume of the intersection of two {$L^n_p$} balls.
\newblock \emph{Proceedings of the American Mathematical Society}, 110\penalty0
  (1):\penalty0 217--224, 1990.

\bibitem[Shamir(2013)]{Shamir13}
O.~Shamir.
\newblock On the complexity of bandit and derivative-free stochastic convex
  optimization.
\newblock In \emph{Proc. 30th Annual Conference on Learning Theory}, pages
  1--22, 2013.

\bibitem[Shamir(2017)]{Shamir17}
O.~Shamir.
\newblock An optimal algorithm for bandit and zero-order convex optimization
  with two-point feedback.
\newblock \emph{Journal of Machine Learning Research}, 18\penalty0
  (1):\penalty0 1703--1713, 2017.

\bibitem[Tsybakov(2009)]{Tsybakov09}
A.~B. Tsybakov.
\newblock \emph{Introduction to Nonparametric Estimation}.
\newblock Springer, New York, 2009.

\bibitem[Zorich(2016)]{zorich2016}
V.~A. Zorich.
\newblock \emph{Mathematical analysis II}.
\newblock Springer, 2016.

\end{thebibliography}
%\newpage

\appendix
% \setcounter{page}{1}
% \section*{Appendix}

% \setcounter{page}{1}
\begin{center}
    {\Large\bf Appendix}
    % \vspace{1cm}
    % {\large CONTENTS}
\end{center}

% \startcontents[appendices]
% \printcontents[appendices]{l}{1}{\setcounter{tocdepth}{1}}

% The appendix is organized as follows. In Appendix ...

In this appendix we first provide some auxiliary results and then prove the results stated in the main body of the paper.  

\paragraph{Additional notation} 
%We also require the following additional notation. 
Let $\bW_1, \bW_2$ be two random variables, we write $\bW_1 \stackrel{d}{=} \bW_2$ to denote their equality in distribution. We also denote by $\Gamma : \bbR_+ \to \bbR_+$ the gamma function defined, for every $z > 0$, as
\(\Gamma(z) = \int_{0}^\infty x^{z - 1} \exp(-x) \d x\).

\section{Consequences of the smoothness assumption}\label{app:A}
Let us first provide some immediate consequences of the smoothness assumption that we consider. 

\begin{remark}
\label{rem:multi-index_juggle}
%One can verify that 
For all $k \in \bbN \setminus \{0\}$ and all $\bh \in \bbR^d$ it holds that
\[
%\begin{align*}
    f^{(k)}(\bx)[\bh]^k
%&
= 
    \sum_{|\bm_1| = \cdots = |\bm_k| = 1} D^{\bm_1 + \cdots + \bm_k}f(\bx) \bh^{\bm_1 + \cdots + \bm_k}
    =
    \sum_{|\bm| = k} \frac{k!}{\bm!}D^{\bm}f(\bx) \bh^{\bm}\enspace.
%\end{align*}
\]
\end{remark}
\begin{proof}~
The first equality of the remark follows from the definition. For the second one it is sufficient to show that for each $\bm = (m_1, \ldots, m_d)^\top \in \bbN^d$ with $|\bm| = k$ there exist exactly $k! / \bm!$ distinct choices of $(\bm_1, \ldots, \bm_k) \in (\bbN^d)^k$ with $|\bm_1| = \ldots = |\bm_k| = 1$ and $\bm_1 + \ldots + \bm_k = \bm$.
% \evg{It should be a classical fact, but I do not have a reference}
To see this, we map $\bm \in \bbN^d$ into a \emph{word} containing letters from $\{a_1, a_2, \ldots, a_d\}$ as
\begin{align*}
    \bm \mapsto W(\bm) \triangleq \underbrace{a_1\ldots a_1}_{m_1-\text{times}}\underbrace{a_2\ldots a_2}_{m_2-\text{times}}\ldots \underbrace{a_d\ldots a_d}_{m_d-\text{times}}\enspace.
\end{align*}
By construction, each letter $a_j$ is repeated exactly $m_j$-times in $W(\bm)$.
Furthermore, if $|\bm| = k$, then $W(\bm)$ contains exactly $k$ letters.
From now on, fix an arbitrary $\bm \in \bbN^d$ with $|\bm| = k$.
Given $(\bm_1, \ldots, \bm_k) \in (\bbN^d)^k$ such that $|\bm_1| = \ldots = |\bm_k| = 1$ and $\bm_1 + \ldots + \bm_k = \bm$, define\footnote{The summation of words is defined as concatenation.}
\begin{align*}
    (\bm_1, \ldots, \bm_k) \mapsto W(\bm_1) + W(\bm_2) + \ldots + W(\bm_k)\enspace.
\end{align*}
We observe that the condition $\bm_1 + \ldots + \bm_k = \bm$, implies that the word  $W(\bm_1) + W(\bm_2) + \ldots + W(\bm_k)$ is a permutation of $W(\bm)$. A standard combinatorial fact states that the number of distinct permutations of $W(\bm)$ is given by the multinomial coefficient, i.e., by $k! / \bm!$. Since the mapping $(\bm_1, \ldots, \bm_k) \mapsto W(\bm_1) + W(\bm_2) + \ldots + W(\bm_k)$ is invertible, we conclude.
\end{proof}

\begin{lemma}
\label{lem:f_to_grad}
Assume that $f \in \mathcal{F}_{\beta}(L)$ for some $\beta \geq 2$ and $L > 0$. Let $\bv \in \mathbb{R}^d$ with $\|\bv\| = 1$ and defined the function $g_{\bv} : \mathbb{R}^d \rightarrow \mathbb{R}$ as $g_{\bv}(x) \equiv \scalar{\bv}{\nabla f(x)}$, $x\in \mathbb{R}^d$. Then $g_{\bv} \in \mathcal{F}_{\beta - 1}(L)$.
% continuously differentiable and that there exist $\alpha \in [0, 1]$ and $L > 0$ such that for all $x, y \in \bbR^d$
% \begin{align}
%     \norm{f^{(\ell)}(x) - f^{(\ell)}(y)} \leq L \|x - y\|^\alpha \enspace.\label{eq:holder_alternative}
% \end{align}
%Then, for any $\bv \in \bbR^d$ with $\|\bv\| = 1$ the function $g_{\bv}(\cdot) \equiv \scalar{\bv}{\nabla f(\cdot)} \in \mathcal{F}_{\beta - 1}(L)$.
% % Then,
% \begin{align*}
%     \norm{g^{(\ell - 1)}_v(x) - g^{(\ell - 1)}_v(y)} \leq L \|x - y\|^\alpha\enspace.
% \end{align*}
\end{lemma}
\begin{proof}
Set $\ell \triangleq \floor{\beta}$.
Note that since $f$ is $\ell$ times continuously differentiable, then $g_{\bv}$ is $\ell - 1$ times continuously differentiable.
Furthermore, for any $\bh^1, \ldots, \bh^{\ell - 1} \in \bbR^d$
\begin{align*}
    g^{(\ell - 1)}_{\bv}(\bx)[\bh^1, \ldots, &\bh^{\ell - 1}]
    =
    \sum_{|\bm_1| = \ldots = |\bm_{\ell - 1}| = 1} D^{\bm_1 + \ldots + \bm_{\ell-1}} g_{\bv}(\bx) h_1^{\bm_1} \cdot \ldots \cdot \bh_{\ell - 1}^{m_{\ell - 1}}\\
    &=
    % \sum_{|\bm_1| = \ldots = |\bm_{\ell - 1}| = 1} D^{\bm_1 + \ldots + \bm_{\ell-1}} \parent{\sum_{|\bm_\ell| = 1}D^{\bm_\ell}f(\bx)\bv^{\bm_\ell}} \bh_1^{\bm_1} \cdot \ldots \cdot \bh_{\ell - 1}^{\bm_{\ell - 1}}\\
    % &=
     \sum_{|\bm_1| = \ldots = |\bm_{\ell}| = 1} D^{\bm_1 + \ldots + \bm_\ell} f(\bx) h_1^{\bm_1} \cdot \ldots \cdot \bh_{\ell - 1}^{\bm_{\ell - 1}}\bv^{\bm_\ell}\\
    &=
    f^{(\ell)}(\bx)[\bh^1, \ldots, \bh^{\ell - 1}, \bv]\enspace.
\end{align*}
Hence, for any $\bx, \bz \in \bbR^d$ we can write by definition of the norm of a $\ell{-}1$-linear form
\begin{align*}
    &\norm{g^{(\ell - 1)}_{\bv}(\bx) - g^{(\ell - 1)}_{\bv}(\bz)}\\
        &\phantom{\leq}= \sup\enscond{\abs{g^{(\ell - 1)}_{\bv}(\bx)[\bh^1, \ldots, \bh^{\ell-1}] - g^{(\ell - 1)}_{\bv}(\bz)[\bh^1, \ldots, \bh^{\ell-1}]}}{\|\bh^j\| = 1\,\,j \in [\ell-1]}\\
        &\phantom{\leq}=
        \sup\enscond{\abs{f^{(\ell)}(\bx)[\bh^1, \ldots, \bh^{\ell - 1}, \bv] - f^{(\ell)}(\bz)[\bh^1, \ldots, \bh^{\ell - 1}, \bv]}}{\|\bh^j\| = 1\,\,j \in [\ell-1]}\\
        % &\leq
        % \sup\enscond{\abs{f^{(\ell)}(\bx)[\bh^1, \ldots, \bh^{\ell - 1}, \bh^\ell] - f^{(\ell)}(\bz)[\bh^1, \ldots, \bh^{\ell - 1}, \bh^\ell]}}{\|\bh^j\| = 1\,\,j \in [\ell]}\\
        &\phantom{\leq}\leq
        \norm{f^{(\ell)}(\bx) - f^{(\ell)}(\bz)} \leq L \|\bx - \bz\|^{\beta - \ell}\enspace.
\end{align*}
% The proof is concluded.
\end{proof}

\begin{lemma}
\label{lem:holder_to_taylor}
Fix some real $\beta \geq 2$ and assume that $f \in \class{F}_{\beta}(L)$.
% {$f : \bbR^d \to \bbR$ is $(\beta, L)$-(modified) H\"older.}
Then, for all $\bx, \bz \in \bbR^d$
\begin{equation*}
\bigg|f(\bx)-\sum_{0\le |\bm|\leq \ell} \frac{1}{\bm!}D^{\bm}f(\bz)(\bx-\bz)^{\bm} \bigg|\leq \frac{L}{\ell!} \|\bx-\bz\|^{\beta}\enspace.
\end{equation*}
\end{lemma}
\begin{proof}
Fix some $\bx, \bz \in \bbR^d$. Taylor expansion yields that, for some $c \in (0, 1)$, 
\begin{align*}
    f(\bx) = \sum_{0\le |\bm|\leq \ell -1} \frac{1}{\bm!}D^{\bm}f(\bz)(\bx-\bz)^{\bm} + \sum_{|\bm| = \ell} \frac{1}{\bm!}D^{\bm}f(\bz + c(\bx - \bz))(\bx-\bz)^{\bm}\enspace.
\end{align*}
Thus, using Remark~\ref{rem:multi-index_juggle} and the fact that $f \in \class{F}_{\beta}(L)$, we get
\begin{align*}
    \bigg|f(\bx)-\sum_{|\bm|\leq \ell} \frac{1}{\bm!}D^{\bm}f(\bz)(\bx-\bz)^{m} \bigg|
    &=
    \abs{\sum_{ |\bm| = \ell} \frac{1}{\bm!}\left(D^{\bm}f(\bz + c(\bx - \bz)) - D^{\bm}f(\bz)\right)(\bx-\bz)^{\bm}}\\
    &=
    \frac{1}{\ell!}\abs{f^{(\ell)}(\bz + c(\bx - \bz))[\bx - \bz]^{\ell} - f^{(\ell)}(\bz)[\bx - \bz]^\ell}\\
    &\leq
    \frac{L}{\ell!} \norm{\bx - \bz}^{\ell} \norm{c(\bx - 
    \bz)}^{\beta - \ell} \leq \frac{L}{\ell!}\norm{\bx - \bz}^\beta\enspace.
\end{align*}
% which concludes the proof.
% \begin{align*}
%     f(x) = \sum_{k = 0}^\ell \frac{1}{k!}f^{(k)}(y)[x - y]^k
% \end{align*}
\end{proof}

\section{Bias and variance of the gradient estimators}

% Let us first recall the following result, which can be found in~\cite[Section 13.3.5, Exercise 14a]{zorich2016}. It was often used without a reference in the context of zero-order optimization with noisy Oracle~\citep{flaxman2004,BP2016,NY1983}.

% \begin{theorem}[Integration by parts in a multiple integral]
% \label{thm:ipp}
%     Let $D$ be an open connected subset of $\bbR^d$ with a piecewise smooth boundary $\partial D$ oriented by the outward unit normal $\bn = (n_1, \ldots, n_m)^\top$. Let $f$ be a smooth function in $D \cup \partial D$, then
%     \begin{align*}
%         \int_{D} \nabla f(\bx) \d \bx = \int_{\partial D} f(\bx) \bn(\bx) \d S(\bx)\enspace.
%     \end{align*}
% \end{theorem}
% % \evg[inline]{Added the above result. In Zorich $D$ should be an open set, it seems to me that it is not an issue to replace $D$ by its closure (since $\d\bx$ is the Lebesgue measure). Zorich calls the above formula as: Integration by parts in a multiple integral and not Stokes' theorem.}
% \begin{remark}
% We refer to~\cite[Section 12.3.2, Definitions 4 and 5]{zorich2016} for the definition of piecewise smooth surfaces and their orientations respectively.
% \end{remark}

\subsection{Gradient estimator with \texorpdfstring{$\ell_2$}{} randomization}

In this section we prove Lemma~\ref{lem:bias_sphere} that provides a bound on the 
bias of the gradient estimator with $\ell_2$ randomization. The variance of this estimator is evaluated in Lemma~\ref{lem:var_sphere_v2} in the main body of the paper.   

We will need the following auxiliary lemma.

\begin{lemma}
\label{tru}
Let $f : \mathbb{R}^d \to \mathbb{R}$ be a continuously differentiable function.
% has a directional derivative which belongs to $\mathcal{F}_{\beta-1}$(L), for $\beta \geq 2$.
Let $r, \bU^{\circ}, \bzeta^{\circ}$ be uniformly distributed on $[-1,1], \ball^d_2$, and $\sphere^d_2$, respectively. Then, for any $h >0$, we have 
\[
\Exp[\nabla f(\bx+hr\bU^{\circ})rK(r)]=\frac{d}{h}\Exp[f(\bx + hr\bzeta^{\circ})\bzeta^{\circ} K(r)]\enspace.
\]
\end{lemma}
\begin{proof}
Fix $r \in [-1,1] \setminus \{0\}$. Define $\phi:\mathbb{R}^d \to \mathbb{R}$ as $\phi (\bu) = f(\bx + hr\bu) K(r)$ and note that $\nabla \phi(\bu) = hr\nabla f(\bx + hr\bu)K(r)$. Hence, we have
\begin{align*}
    \Exp[\nabla f(\bx+hr\bU^{\circ})K(r) \mid r] = \frac{1}{hr}\Exp[\nabla \phi(\bU^{\circ}) \mid r]
    &=
    \frac{d}{hr}\Exp[\phi(\bzeta^{\circ})\bzeta^{\circ} \mid r]\\
    &=
    \frac{d}{hr}K(r)\Exp[f(\bx+hr\bzeta^{\circ})\bzeta^{\circ} \mid r]\enspace,
\end{align*}
where the second equality is obtained from a version of Stokes' theorem~\cite[see e.g.,][Section 13.3.5, Exercise 14a]{zorich2016}. Multiplying by $r$ from both sides, using the fact that $r$ follows a continuous distribution, and taking the total expectation concludes the proof.
\end{proof}

\begin{proof}[Proof of Lemma~\ref{lem:bias_sphere}]
Using Lemma~\ref{tru}, the fact that $\int_{-1}^{1} rK(r)\d r = 1$, and the variational representation of the Euclidiean norm, we can write 
\begin{align}
    \label{eq:bias_new1}
    \norm{\Exp[\bg_t^{\circ} \mid \bx_t]-\nabla f(\bx_{t})} = \sup_{\bv \in \sphere^d_2}\Exp[\big(\nabla_{\bv} f(\bx+h_tr_t\bU^{\circ}) - \nabla_{\bv}f(\bx)\big)r_tK(r_t)]\enspace,
\end{align}
where we recall that $\bU^{\circ}$ is uniformly distributed on $\ball^d_2$. Lemma~\ref{lem:f_to_grad} asserts that for any $\bv \in \sphere^d_2$ the directional gradient $\nabla_{\bv}f(\cdot)$ is $(\beta-1, L)$-H\"older.
Thus, due to Lemma~\ref{lem:holder_to_taylor} we have the following Taylor expansion: 
\begin{small}
\begin{align}
    \label{eq:taylor_grad}
    \nabla_{\bv}f(\bx_t+h_tr_t\bU^{\circ}) = \nabla_{\bv} f(\bx_t) +\sum_{1\leq|\bm|\leq \ell -1}\frac{(r_th_t)^{|\bm|}}{\bm!}D^{\bm}\nabla_{\bv}f(\bx_t)({\bU^{\circ}})^{\bm}+R(h_tr_t\bU^{\circ})\enspace,
\end{align}
\end{small}
where the residual term $R(\cdot)$ satisfies $|R(\bx)|\leq \tfrac{L}{(\ell - 1)!} \norm{\bx}^{\beta - 1}$.

Substituting~\eqref{eq:taylor_grad} in~\eqref{eq:bias_new1} and using the ``zeroing-out'' properties of the kernel $K$, we deduce that
\begin{align*}
    \norm{\Exp[\bg_t^{\circ} \mid \bx_t]- \nabla f(\bx_{t})}\leq \kappa_{\beta}h_{t}^{\beta-1}\frac{L}{(\ell - 1)!}\Exp\norm{\bU^{\circ}}^{\beta-1} = \kappa_{\beta}h_{t}^{\beta-1} \frac{L}{(\ell - 1)!} \frac{d}{d+\beta-1}\enspace,
\end{align*}
where the last equality is obtained from the fact that $\Exp\norm{\bU^{\circ}}^q = \frac{d}{d+q}$, for any $q \geq0$.
\end{proof}

% \newpage
\subsection{Gradient estimator with \texorpdfstring{$\ell_1$}{} randomization}

In this section, we prove Lemma~\ref{lem:bias_l1} that gives a bound 
on the bias of our gradient estimator with $\ell_1$ randomization, and Lemma~\ref{lem:true_poincare},
that provides a Poincaré type inequality crucial for the control of its variance.

Let $\zeta$ be a real valued random variable with $\Exp[\zeta^2] \leq 4\sigma^2$ and let $\bzeta^{\diamond}$ be distributed uniformly on $\sphere^d_1$. Assume that $\zeta$ and $\bzeta^{\diamond}$ are independent from each other and from the random variable $r$, which is uniformly distributed on $[-1,1]$. In order to control the bias and variance of the gradient estimator \eqref{eq:grad_l1} for any fixed $t$, it is sufficient to do it for the random vector
\begin{align}
    \label{eq:grad_simplex}
    \bg_{\bx, h}^{\diamond} = \frac{d}{2h}\parent{f(\bx + hr\bzeta^{\diamond}) - f(\bx - hr\bzeta^{\diamond}) + \zeta}\sign(\bzeta^{\diamond}) K(r),
\end{align}
where $\zeta = \xi_t - \xi'_t$.

\subsubsection{Control of the bias}

\begin{lemma}
    \label{cor:simplex_unbiased}
    Let $\bU^{\diamond}$ be uniformly distributed on $\ball^d_1$ and $\bzeta^{\diamond}$ be uniformly distributed on $\sphere^d_1$.
    Fix some $\bx \in \bbR^d$ and $h > 0$, and let Assumption~\ref{ass1} be fulfilled, then the estimator in~\eqref{eq:grad_simplex} satisfies
    \begin{align*}
    \Exp[\bg_{\bx, h}^{\diamond}] = \Exp[\nabla f(\bx + hr\bU^{\diamond})rK(r)]\enspace.
    \end{align*}
\end{lemma}
\begin{proof}
The proof is analogous to that of Lemma~\ref{tru} using~\cite[Theorem 6]{Akhavan_Chzhen_Pontil_Tsybakov22a}.
\end{proof}

In order to obtain a bound on the bias of the estimator in~\eqref{eq:grad_simplex} we need the following result, which controls the moments of the Euclidean norm of $\bU^{\diamond}$.

% \evg[inline]{Here there is a mild inaccuracy. Need to consider the case of $\beta \in [2, 1)$ separately from $\beta > 2$. For $\beta \in [2, 1)$ Jensen's inequality will give the result. To come back!}
%
\begin{lemma}
\label{lem:moments_l1}
Let $\bU^{\diamond} \in \bbR^d$ be distributed uniformly on $\ball^d_1$. Then for any {\color{purple}$\beta \geq 1$ it holds that
\begin{align*}
    \Exp [\norm{\bU^{\diamond}}^\beta ]\leq \frac{c_{\beta+1}d^{\frac{\beta}{2}}\Gamma(\beta + 1)\Gamma(d + 1)}{\Gamma(d + \beta + 1)}\enspace,
\end{align*}
where $c_{\beta+1}=2^{\beta/2}$ for $1\leq \beta <2$ and $c_{\beta+1}=1$ for $\beta \geq 2$.}
\end{lemma}
\begin{proof}
% \evg{I modified here a bit to make it cleaner}
Let $\bW = (W_1, \ldots, W_d), W_{d+1}$ be \iid random variables following Laplace distribution with mean $0$ and scale parameter $1$.
% and $W \sim \exp(1)$ being independent from $U_1, \ldots, U_d$.
Then, following \cite[Theorem 1]{Barthe_Guedon_Mendelson_Naor05} we have
\begin{align*}
    \bU^{\diamond}\stackrel{d}{=} \frac{\bW}{\norm{\bW}_1 + |W_{d+1}|}\enspace,
\end{align*}
where the sign $\stackrel{d}{=} $ stands for equality in distribution. Furthermore, it follows from \cite[Theorem 2]{Barthe_Guedon_Mendelson_Naor05} (see also~\cite{Rachev_Ruschendorf91,Schechtman_Zinn90}) that
\begin{align*}
    \frac{(\bW, |W_{d+1}|)}{\norm{\bW}_1 + |W_{d+1}|}\qquad\text{and}\qquad \norm{\bW}_1 + |W_{d+1}|\enspace,
\end{align*}
are independent.
Hence, 
\begin{align}
    \label{eq:moment_l1_1}
   \Exp [\norm{\bU^{\diamond}}^\beta]= \Exp \left[\parent{\frac{\sum_{j = 1}^d W_j^2}{\parent{\norm{\bW}_1 + |W_{d+1}|}^2}}^{\beta / 2}\right]
    % \leq
    % \Exp \left[\frac{\norm{\bW}^{\beta}}{\norm{(\bW, W_{d+1})}_1^{\beta}}\right]
    = \frac{\Exp[\norm{\bW}^\beta]}{\Exp[\norm{(\bW, W_{d+1})}_1^\beta]}\enspace,
\end{align}
where the equality follows from the independence recalled above. Note that $|W_j|$ is $\exp(1)$ random variable for any $j = 1, \ldots, d$.
{\color{purple}Thus, if $1 \leq \beta < 2$ by Jensen's inequality we can write
\begin{align}\label{eq:moment_l1_2}
    \Exp[\norm{\bW}^\beta] = \Exp\parent{\sum_{j  = 1}^d W_j^2}^{\frac{\beta}{2}} \leq \parent{\sum_{j = 1}^d \Exp[W_j^2]}^{\frac{\beta}{2}} = d^{\frac{\beta}{2}}\parent{\Exp[W_1^2]}^{\frac{\beta}{2}}=d^{\frac{\beta}{2}}\Gamma(3)^{\frac{\beta}{2}}\enspace.
\end{align}
If $\beta \geq 2$, again by Jensen's inequality we have
\begin{align}
\label{eq:moment_l1_2_1}
    \Exp[\norm{\bW}^\beta] = d^{\frac{\beta}{2}}\Exp\parent{\frac{1}{d}\sum_{j  = 1}^d W_j^2}^{\frac{\beta}{2}} \leq d^{\frac{\beta}{2} - 1} \sum_{j = 1}^d \Exp[W_j^\beta] = d^{\frac{\beta}{2}} \Exp[W_1^\beta] = d^{\frac{\beta}{2}} \Gamma(\beta + 1)\enspace.
\end{align}}
It remains to provide a suitable expression for $\Exp[\norm{(\bW, W_{d+1})}_1^\beta]$. We observe that $\norm{(\bW, W_{d+1})}_1$ follows the Erlang distribution with parameters $(d+1, 1)$ (as a sum of $d+1$ \iid $\exp(1)$ random variables).
Hence, using the expression for the density of the Erlang distribution we get
\begin{align}
\label{eq:moment_l1_3}
    \Exp[\norm{(\bW, W_{d+1})}_1^\beta] = \frac{1}{\Gamma(d+1)}\int_{0}^{\infty}x^{d + \beta}\exp(-x) \d x =  \frac{\Gamma(d + \beta + 1)}{\Gamma(d + 1)}\enspace.
\end{align}
Combining~\eqref{eq:moment_l1_1}--\eqref{eq:moment_l1_3}  proves the lemma.
\end{proof}

%We are now in a position to derive an upper bound on the bias of the gradient estimator~\eqref{eq:grad_simplex}.

\begin{proof}[Proof of Lemma~\ref{lem:bias_l1}]
Using Lemma~\ref{cor:simplex_unbiased} and following the same lines as in the proof of Lemma~\ref{lem:bias_sphere} we deduce that
% \evg{here also modified accordingly.}
\begin{align*}
    \norm{\Exp[\bg_t^{\diamond} \mid \bx_t]- \nabla f(\bx_{t})}
    \leq
    \kappa_{\beta}h_{t}^{\beta-1} \frac{L}{(\ell - 1)!}\Exp \norm{\bU^{\diamond}}^{\beta-1}
    \leq
    \kappa_{\beta}h_{t}^{\beta-1} \frac{L}{(\ell - 1)!}\frac{{\color{purple}c_{\beta}}d^{\frac{\beta - 1}{2}}\Gamma(\beta)\Gamma(d + 1)}{\Gamma(d + \beta)}\enspace,
\end{align*}
where the last inequality is due to Lemma~\ref{lem:moments_l1}. {
% \color{black}
Next, recall that the Gamma function satisfies $\Gamma(z+1) = z\Gamma(z)$ for any $z>0$. Applying this relation iteratively and using the fact that $\ell=\lfloor\beta\rfloor$ we get: 
\begin{align*}
    \frac{\Gamma(d + 1)}{\Gamma(d+\beta)}
    &= \frac{\Gamma(d+1)}{\Gamma\big(d+\underbrace{(\beta-\ell)}_{\in (0, 1]}\big)\prod_{i = 1}^{\ell}\big(d+\beta-i\big) }
    \leq
    \frac{(d + \beta - \ell)^{1 - (\beta -\ell)}}{\prod_{i = 1}^{\ell}\big(d+\beta-i\big) }
    % \leq
    % \frac{(d + \beta - \ell)^{\beta - \ell}}{\prod_{i = 1}^{\ell - 1}\big(d+\beta-i\big) }
    % \frac{\Gamma(d+1)}{d^{\ell}\Gamma\big(d+\underbrace{\beta-\ell}_{\in (0, 1]}\big)}
    \leq \frac{1}{d^{\beta-1}}\enspace,
\end{align*}
where the first inequality is obtained from \cite[Remark 1]{fengqi}. Proceeding analogously we obtain that  $\frac{\Gamma(\beta)}{(\ell-1)!}\le \ell^{\beta-\ell}$. Combining this bound with the two preceding displays yields the lemma.
% \begin{align*}
%     \frac{\Gamma(\beta)}{(\ell-1)!} = \frac{\Gamma(\beta)}{\Gamma(\ell)}= \frac{\ell \Gamma(\ell + (\beta-\ell))}{\Gamma(\ell+1)}\leq \ell^{\beta-\ell}\enspace.\qedhere
% \end{align*}
% where we again used \cite[Remark 2.1.1]{fengqi} for the first inequality. Combining the above derived bounds we conclude.
% \begin{align*}
%     \norm{\Exp[\bg_t \mid \bx_t]- \nabla f(\bx_{t})}\leq \kappa_{\beta}h_{t}^{\beta-1} \ell{L} d^{\frac{1-\beta }{2}}\enspace.
% \end{align*}
}
\end{proof}

{
\subsubsection{Poincaré inequality for the control of the variance}
% \color{black}
% \evg[inline]{Better variance for $\ell_1$}

We now prove the Poincaré inequality of Lemma~\ref{lem:true_poincare} used to control the variance of the $\ell_1$-randomized estimator.

\begin{proof}[Proof of Lemma~\ref{lem:true_poincare}]
The beginning of the proof is the same as in~\cite[Lemma~3]{Akhavan_Chzhen_Pontil_Tsybakov22a}.
In particular,  without loss of generality we assume that $\Exp[G(\bzeta)] = 0$, and
consider first the case of continuously differentiable $G$. 
Let $\bW = \left(W_1, \dots, W_d\right)$ be a vector such that the components $W_j$ are i.i.d. Laplace random variables with mean $0$ and scale parameter $1$.  Set $\bT(\bw) = \bw / \norm{\bw}_1$.
% Introduce the notation
% \begin{align*}
%     F(\bw) \eqdef \|\bw\|_1^{\sfrac 1 2}G(\bT(\bw))\enspace.
% \end{align*}
Lemma~1 in~\cite{Schechtman_Zinn90} asserts that, for $\bzeta$ uniformly distributed on $\partial B_1^d$,
\begin{align}
    \label{eq:SZ90}
    \bT(\bW) \stackrel{d}{=} \bzeta\quad\text{and}\quad \bT(\bW)  
    \text{ is independent of } \|\bW\|_1 \enspace.
\end{align}
Furthermore, in the proof of Lemma~3 in~\cite{Akhavan_Chzhen_Pontil_Tsybakov22a}, it is shown that
\begin{align}
\label{eq:poincare_new_10}
   \Var(G(\bzeta)) \leq \frac{4}{{d(d-2)}} \Exp\left[\norm{\bfI - \bT(\bW)\big(\sign(\bW)\big)^{\top}}^2\norm{\nabla G(\bT(\bW))}^2\right],
\end{align}
% and
% \begin{align*}
%     4\|\nabla F(\bW)\|^2
%     &=
%     \frac{d}{\|\bW\|_1}G^2(\bT(\bW)) + \frac{4}{\|\bW\|_1}\norm{\bigg(\bfI - \bT(\bW)\big(\sign(\bW)\big)^{\top}\bigg)\nabla G(\bT(\bW))}^2
% \end{align*}
% From now on, the proof slightly diverges from that of~\cite{Akhavan_Chzhen_Pontil_Tsybakov22a}. First, we also use the fact that
% \begin{align*}
%     \norm{\bigg(\bfI - \bT(\bW)\big(\sign(\bW)\big)^{\top}\bigg)\nabla G(\bT(\bW))}^2 \leq \norm{\bfI - \bT(\bW)\big(\sign(\bW)\big)^{\top}}^2\norm{\nabla G(\bT(\bW))}^2\enspace.
% \end{align*}
where $\bfI$ is the identity matrix, $\|\cdot\|$ applied to matrices denotes the spectral norm, and $\sign(\cdot)$ applied to vectors denotes the vector of signs of the coordinates. 

From this point, the proof diverges from that of Lemma~3 in~of~\cite{Akhavan_Chzhen_Pontil_Tsybakov22a}.
 Instead of bounding the spectral norm of $\bfI - \ba\bb^\top$ by $1 + \norm{\ba}\norm{\bb}$ as it was done in that paper, we compute it exactly, which leads to the main improvement. Namely, Lemma~\ref{lem:perturbation} proved below gives
\begin{align*}
    \norm{\bfI - \bT(\bW)\big(\sign(\bW)\big)^{\top}}^2 = d\norm{\bT(\bW)}^2\enspace.
\end{align*}
Combining this equality with~\eqref{eq:poincare_new_10} we obtain the first bound of the lemma.
% and the independence of $\|\bW\|_1$ and $\bT(\bW)$ (cf. \eqref{eq:SZ90}) yields
% \begin{align*}
%     d\parent{1 - \frac{1}{d-1}}\Var(G(\bzeta)) \leq \frac{4d}{d-1}\Exp\left[\|\nabla G(\bT(\bW))\|^2\|\bT(\bW)\|^2\right]\enspace.
% \end{align*}
% Rearranging, we deduce the first claim of the lemma since $\bT(\bW) \stackrel{d}{=} \bzeta$.
The second bound of the lemma (regarding Lipschitz functions $G$) is deduced from the first one by the same argument as in~\cite{Akhavan_Chzhen_Pontil_Tsybakov22a}.
\end{proof}

\begin{lemma}
    \label{lem:perturbation}
    Let $\ba \in \bbR^d$ be such that $\|\ba\|_1 = 1$.
    Then,
    \begin{align*}
        \norm{\bfI - \ba \big(\sign(\ba)\big)^\top} = \sqrt{d}\norm{\ba}\enspace.
    \end{align*}
\end{lemma}
%Note that the first proof of Poincaré inequality for the uniform measure on the $\ell_1$-sphere, see \cite{Akhavan_Chzhen_Pontil_Tsybakov22a}, used the fact that $\|{\bfI - \ba \big(\sign(\ba)\big)^\top}\| \leq 1 + \sqrt{d}\norm{\ba}$, which resulted in inflated constants. Lemma~\ref{lem:perturbation} shows that the norm $\|{\bfI - \ba \big(\sign(\ba)\big)^\top}\|$ can be computed explicitly. 

\begin{proof}[Proof of Lemma~\ref{lem:perturbation}]
    Let $\bu = \ba / \norm{\ba}$, $\bv = \sign(\ba) / \sqrt{d}$, and $\gamma = \sqrt{d}\norm{\ba}$. Then, since $1 = \norm{\ba}_1 = \scalar{\ba}{\sign(\ba)}$ we have $\scalar{\bu}{\bv} = 1/\gamma$.
    Consider the matrix $\bfQ = [\bu, \bq_2, \ldots, \bq_d]$, such that $\bfQ^\top\bfQ = \bfQ\bfQ^\top = \bfI$. Let $\be_1 = (1, 0, 
    \ldots, 0)^\top$.
    For any matrix ${\bf B}\in \bbR^{d\times d}$ we have $\norm{\bf B} = \norm{\bfQ^\top \bf B \bfQ}$ and $\norm{\bf B}^2 = \norm{\bf B\bf B^\top}$.
    Using these remarks and the fact that $\norm{\bfQ \bv}^2 = 1$, $\bfQ^\top \bu = \be_1$, we deduce that
    \begin{align*}
        \norm{\bfI - \ba \big(\sign(\ba)\big)^\top}^2
        &=
       \norm{{(\bfI - \gamma\be_1 (\bfQ^\top\bv)^\top)(\bfI - \gamma\be_1 (\bfQ^\top\bv)^\top)^\top}}\\
       &=
       \norm{{\bfI - \gamma\be_1 (\bfQ^\top\bv)^\top - \gamma(\bfQ^\top\bv)\be_1 + \gamma^2\be_1\be_1^\top}} = \norm{\bfA}\enspace,
    \end{align*}
    where 
    \begin{align*}
        \bfA = \left[\begin{array}{ c | c }
    \gamma^2 - 1 & -\gamma\bar{\bv}^\top \\
    \hline
    -\gamma\bar{\bv} & \bfI
  \end{array}\right]\enspace,
    \end{align*}
    with $(\bar{\bv})_j = \scalar{\bq_{j+1}}{\bv}$ for $j = 1, \ldots, d-1$.
    % We have sown that
    % \begin{align*}
    %     \|\bfI - \ba \bc^\top\|^2 = \norm{\bfA}\enspace.
    % \end{align*}
    Let us find the eigenvalues of $\bfA$. For any $\lambda \in \bbR$, using the expression for the determinant of a block matrix we get
    \begin{align*}
        \det(\bfA - \lambda \bfI) = (1 - \lambda)^{d - 1}\parent{\gamma^2 - 1 - \lambda -\frac{\gamma^2\norm{\bar{\bv}}^2}{1-\lambda}}\enspace.
    \end{align*}
    Note that $1 = \norm{\bfQ \bv}^2 = \frac{1}{\gamma^2} + \norm{\bar{\bv}}^2$.
    Hence, 
    \begin{align*}
        \det(\bfA - \lambda \bfI)
        &=
        (1 - \lambda)^{d - 2}\parent{(1 - \lambda)(\gamma^2 - 1 - \lambda) - (\gamma^2 - 1)}
        = (1 - \lambda)^{d - 2}(\lambda - \gamma^2)\lambda\enspace.
    \end{align*}
    Thus, $\|\bfI - \ba \big(\sign(\ba)\big)^\top\| = \max\{\gamma, 1\} = \max\{\sqrt{d}\norm{\ba}, 1\}$. We conclude the proof by observing that $\sqrt{d}\norm{\ba} \geq \norm{\ba}_1 = 1 $.
\end{proof}

    Finally, we provide the following auxiliary lemma used in the proof of Lemma~\ref{lem:var_l1_v2}.

\begin{lemma}
\label{lem:moments_of_strange_things}
For all $d \geq 3$ and all $\bc \in \bbR^d, h>0$ it holds that
% \begin{small}
\begin{align*}
    \Exp(\norm{\bc} + h\norm{\bzeta^{\diamond}})^2\norm{\bzeta^{\diamond}}^2 \leq \frac{2}{d+1}\parent{\norm{\bc} + h\sqrt{\frac{2}{{d}}}}^2\enspace.
\end{align*}
% \begin{align*}
%     &\Exp\left[{(1 {+} \sqrt{d}\norm{\bzeta^{\diamond}})^2}\right] \leq {\parent{1 {+} \sqrt{\frac{2d}{d{+}1}}}^2},\,\,
%     &\Exp\left[{(1 {+} \sqrt{d}\norm{\bzeta^{\diamond}})^2\norm{\bzeta^{\diamond}}^2}\right] \leq
%     \frac{(2+\sqrt{2})^2}{d+1}\,.
%     % \frac{12+\frac{88}{d}}{(d-1)(d-2)(d+1)}\enspace.
% \end{align*}
% \end{small}
\end{lemma}
\begin{proof}
% We start the proof with some preparations.
Observe that the vector $|\bzeta^{\diamond}| \triangleq (|\zeta^{\diamond}_1|, \ldots, |\zeta^{\diamond}_d|)^\top$ follows the Dirichlet distribution (i.e., the uniform distribution of the probability simplex on $d$ atoms).
In what follows we will make use of the following expression for the moments of the Dirichlet distribution:
\begin{align}
    \label{eq:moments_dirichlet}
    \Exp[(\bzeta^{\diamond})^{\bm}] = \frac{\Gamma(d)}{\Gamma(d + |\bm|)}\prod_{i = 1}^d\Gamma(m_i + 1) = \frac{(d-1)!\bm!}{(d - 1 + |\bm|)!}\enspace,
\end{align}
for any multi-index $\bm = (m_1, \ldots, m_d) \in \bbN^d$ with even coordinates.

Using~\eqref{eq:moments_dirichlet} we get %first derive upper bounds on $\Exp\norm{\bzeta^{\diamond}}^2$ and $\Exp\norm{\bzeta^{\diamond}}^4$.
%Eq.~\eqref{eq:moments_dirichlet} yields
\begin{align}
    \label{eq:strange_moments1}
    \Exp\norm{\bzeta^{\diamond}}^2 = \frac{2}{d+1}\enspace.
\end{align}
Furthermore, using the multinomial identity and the expression for the moments in~\eqref{eq:moments_dirichlet} we find
    % \begin{align*}
    %     \Exp\norm{\bzeta^{\diamond}}^6
    %     =
    %     \sum_{|\bm|=3}\frac{6}{\bm!}\Exp[(\bzeta^{\diamond})^{2\bm}]
    %     =
    %     \sum_{|\bm|=3}\frac{6}{\bm!}\cdot \frac{(d-1)!(2\bm)!}{(d+5)!}
    %     =
    %     6\frac{(d-1)!}{(d+5)!} \sum_{|\bm| = 3} \frac{(2\bm)!}{\bm!}\enspace,
    % \end{align*}
    % and
    \begin{align*}
        \Exp\norm{\bzeta^{\diamond}}^4
        =
        \sum_{|\bm|=2}\frac{2}{\bm!}\Exp[(\bzeta^{\diamond})^{2\bm}]
        =
        \sum_{|\bm|=2}\frac{2}{\bm!}\cdot \frac{(d-1)!(2\bm)!}{(d+3)!}
        =
        2\frac{(d-1)!}{(d+3)!} \sum_{|\bm| = 2} \frac{(2\bm)!}{\bm!}\enspace.
    \end{align*}
  Direct calculations show that
  % $\sum_{|\bm| = 3} \frac{(2\bm)!}{\bm!} = \frac{4}{3}d((d + 3)^2 + 74)$ and
  $\sum_{|\bm| = 2} \frac{(2\bm)!}{\bm!} = 2d(d + 5)$.
% $\sum_{|\bm| = 2} \frac{(2\bm)!}{\bm!} = 2d(d + 11)$.
    % \begin{align*}
    %     &\sum_{|\bm| = 3} \frac{(2\bm)!}{\bm!} = \frac{4}{3}d((d + 3)^2 + 74)\qquad\text{and}\qquad
    %     \sum_{|\bm| = 2} \frac{(2\bm)!}{\bm!} = 2d(d + 11)\enspace.
    % \end{align*}
    Hence, we deduce that, for all $d \geq 1$,
    \begin{align*}
        % &\Exp\norm{\bzeta^{\diamond}}^6 = \frac{8d!(d+3)^2}{(d+5)!}\parent{1 + \frac{74}{(d+3)^2}}\qquad\text{and}\qquad
        % \leq \frac{8d!(d+3)^2}{(d+5)!}\parent{1 + \frac{74}{36}}\enspace,\\
        \Exp\norm{\bzeta^{\diamond}}^4 = \frac{4d!(d+5)}{(d+3)!}\enspace.
    \end{align*}
    Note that $\frac{d(d+5)}{(d+2)(d+3)} \leq 1$ for all $d \geq 1$. Thus,
    \begin{align}
        \label{eq:strange_moments2}
        \Exp\norm{\bzeta^{\diamond}}^4 = \frac{4d!(d+5)}{(d+3)!} = \frac{4(d+5)}{(d+1)(d+2)(d+3)} \leq \frac{4}{d(d+1)}\enspace.
    \end{align}
    Finally, observe that by the Cauchy-Schwarz inequality,
\begin{align*}
    \Exp(\norm{\bc} + h\norm{\bzeta^{\diamond}})^2\norm{\bzeta^{\diamond}}^2
    &\leq
   h^2\Exp\norm{\bzeta^{\diamond}}^4 + 2h\norm{\bc}\sqrt{\Exp\norm{\bzeta^{\diamond}}^2\Exp\norm{\bzeta^{\diamond}}^4} + \norm{\bc}^2\Exp\norm{\bzeta^{\diamond}}^2\\
   &= \parent{\norm{\bc}\sqrt{\Exp\norm{\bzeta^{\diamond}}^2} + h\sqrt{\Exp\norm{\bzeta^{\diamond}}^4}}^2\enspace.
   % &\leq
   % \frac{1}{d+1}\parent{\sqrt{2}\norm{\bc} + \frac{2h}{\sqrt{d}}}^2
   % \enspace.
\end{align*}
Combining this bound with~\eqref{eq:strange_moments1} and~\eqref{eq:strange_moments2} concludes the proof.
% where
% \begin{align*}
%     \texttt{T}_1 &= h^2\Exp\norm{\bzeta^{\diamond}}^2(1 + \sqrt{d}\norm{\bzeta^{\diamond}})^2\enspace,\\
%     \texttt{T}_2 &= 2h\norm{\bc}\Exp\norm{\bzeta^{\diamond}}(1 + \sqrt{d}\norm{\bzeta^{\diamond}})^2\enspace,\\
%     \texttt{T}_3 &= \norm{\bc}^2\Exp(1 + \sqrt{d}\norm{\bzeta^{\diamond}})^2\enspace.
% \end{align*}
% \paragraph{Handling $\texttt{T}_1$.} Using Cauchy-Schwartz inequality, we deduce that
% \begin{align*}
%     \texttt{T}_1 \leq h^2\parent{\Exp\norm{\bzeta^{\diamond}}^2 + 2\sqrt{d\Exp\norm{\bzeta^{\diamond}}^2\Exp\norm{\bzeta^{\diamond}}^4} + d\Exp\norm{\bzeta^{\diamond}}^4} \leq h^2(1+\sqrt{2})^2\frac{2}{d+1}\enspace.
% \end{align*}
% \paragraph{Handling $\texttt{T}_2$.} Similarly, we obtain
% \begin{align*}
%     \texttt{T}_2 \leq 2h\norm{\bc}\sqrt{\Exp\norm{\bzeta^{\diamond}}^2}\parent{1 + 2\sqrt{d\Exp\norm{\bzeta^{\diamond}}^2} + d \sqrt{\Exp\norm{\bzeta^{\diamond}}^4}} \leq 2h(1+\sqrt{2})^2\sqrt{\frac{2}{d+1}}\norm{\bc}\enspace.
% \end{align*}
% \paragraph{Handling $\texttt{T}_3$.} Again, using similar arguments, we obtain
% \begin{align*}
%     \texttt{T}_3 \leq \norm{\bc}^2 \parent{1 + \sqrt{2}}^2\enspace.
% \end{align*}
% We conclude by combining the bounds on all three terms.
\end{proof}

}

\section{A technical lemma}
In this section, we provide a lemma, which will be useful to handle recursive relations in the main proofs. It is a direct extension of~\cite[Lemma D.1]{akhavan2020}.
\begin{lemma}
\label{lemma:bisev}
    Let $\{\delta_{t}\}_{t \geq 1}$ be a sequence of real numbers such that for all integers $t > t_{0} \geq 1$, 
    \begin{equation}
    \label{ch:0}
        \delta_{t+1} \leq \left(1 - \frac{c}{t} \right)\delta_{t} + \sum_{i=1}^{N}\frac{a_{i}}{t^{p_{i}+1}}\enspace,
    \end{equation}
where $c\geq 1$, $p_{i} \in (0,c)$ and $a_{i}\geq0$ for $i \in [N] $. Then for $t \geq t_{0} \geq c+1$, we have
\begin{equation}
\label{ch:res}
        \delta_{t} \leq \frac{2(t_{0}-1)\delta_{t_{0}}}{t}+\sum_{i=1}^{N}\frac{a_{i}}{(c-p_{i})t^{p_{i}}}\enspace.
    \end{equation}
\end{lemma}
\begin{proof}
For any fixed $t > 0$ the convexity of the mapping $u\mapsto g(u)=(t+u)^{-p}$ implies that  $g(1)-g(0)\ge g'(0)$, i.e., $\tfrac{1}{t^{p}}-\tfrac{1}{(t+1)^{p}}
\leq \tfrac{p}{t^{p+1}}$.
Thus, using the fact that $\tfrac{1}{t^p} - \tfrac{p}{t^{p+1}} = \tfrac{(c - p) + (t - c)}{t^{p+1}} \leq \tfrac{1}{(t+1)^p}$,
\begin{eqnarray}
\label{ch:2}
\frac{a_i}{t^{p+1}}\leq \frac{a_i}{c-p}\left\{\frac{1}{(t+1)^{p}} - \Big(1-\frac{c}{t}\Big)\frac{1}{t^{p}} \right\}\enspace.
\end{eqnarray}
Using~\eqref{ch:0}, \eqref{ch:2} and rearranging terms, for any $t \geq t_{0}$ we get 
\begin{eqnarray}
\nonumber
\delta_{t+1}-\sum_{i=1}^{N}\frac{a_{i}}{(c-p_{i})(t+1)^{p_{i}}}\le\Big(1-\frac{c}{t}\Big)\left\{\delta_{t} - \sum_{i=1}^{N}\frac{a_{i}}{(c-p_{i})t^{p_{i}}}\right\}\enspace.
\end{eqnarray}
Letting $\tau_{t} = \delta_{t} - \sum_{i=1}^{N}\tfrac{a_{i}}{(c-p_{i})t^{p_{i}}}$ we have $\tau_{t+1}\leq (1-\tfrac{c}{t}) \tau_{t}$. Now, if $\tau_{t_{0}} \leq 0$ then $\tau_{t} \leq 0$ for any $t \geq t_{0}$ and thus
(\ref{ch:res}) holds. Otherwise, if $\tau_{t_0}>0$ then for $t\ge t_{0}+1$ we have
\[
\tau_{t}\leq \tau_{t_{0}}\prod_{i={t_{0}}}^{t-1}\Big(1-\frac{c}{i}\Big)\leq \tau_{t_{0}}\prod_{i={t_{0}}}^{t-1}\Big(1-\frac{1}{i}\Big)\le \frac{(t_{0}-1)\tau_{t_{0}}}{t}\le \frac{2(t_{0}-1)\delta_{t_{0}}}{t}\enspace.
\]
Thus, \eqref{ch:res} holds in this case as well.
\end{proof}
\section{Upper bounds} 

\subsection{Upper bounds: Only smoothness assumption}

\begin{proof}[Proof of Lemma~\ref{HighSF}]
    For brevity we write $\Exp_t[\cdot]$ in place of $\Exp[\cdot \mid \bx_t]$.
    Using Lipschitz continuity of $\nabla f$ and the definition of the algorithm in~\eqref{eq:algo_general} we can write
    \begin{align*}
        \Exp_t[f(\bx_{t+1})]
        &\leq
        f(\bx_t) - \eta_t \scalar{\nabla f(\bx_t)}{\Exp_t[\bg_t]} + \frac{\bL\eta_t^2}{2}\Exp_t\left[\norm{\bg_t}^2\right]\\
        &\leq
        f(\bx_t) - \eta_t \norm{\nabla f(\bx_t)}^2 + \eta_t\norm{\nabla f(\bx_t)}\norm{\Exp_t[\bg_t] - \nabla f(\bx_t)}
        +
        \frac{\bL\eta_t^2}{2}\Exp_t\left[\norm{\bg_t}^2 \right]\enspace.
    \end{align*}
    Furthermore, invoking the assumption on the bias and the variance of $\bg_t$ and using the fact that $2ab \leq a^2 + b^2$ we deduce
    \begin{equation}
        \label{eq:lemma10}
    \begin{aligned}
        \Exp_t[f(\bx_{t+1})] - &f(\bx_t)
        \leq
        - \eta_t \norm{\nabla f(\bx_t)}^2 + \eta_tb_t\norm{\nabla f(\bx_t)}
        +
        \frac{\bL\eta_t^2}{2}\parent{v_t + {\sf V}_1 \norm{\nabla f(\bx_t)}^2}\\
        &\leq
         - \eta_t \norm{\nabla f(\bx_t)}^2 + \frac{\eta_t}{2}\parent{b_t^2 + \norm{\nabla f(\bx_t)}^2}
        +
        \frac{\bL\eta_t^2}{2}\parent{v_t + {\sf V}_1 \norm{\nabla f(\bx_t)}^2}\\
        &=
         - \frac{\eta_t}{2}\parent{1 - \bL\eta_t{\sf V}_1}\norm{\nabla f(\bx_t)}^2 + \frac{\eta_t}{2}\parent{b_t^2 + \bL\eta_t v_t}\enspace.
    \end{aligned}
    \end{equation}
     Let $S$ be a random variable with values in $\{1, \ldots, T\}$, which is independent from $\bx_1, \ldots, \bx_T, \bg_1, \ldots, \bg_T$ and such that
    \begin{align*}
        \Prob(S = t) = \frac{\eta_t\parent{1 - \bL\eta_t{\sf V}_1}}{\sum_{t = 1}^T\eta_t\parent{1 - \bL\eta_t{\sf V}_1}}\enspace.
    \end{align*}
    Assume that $\eta_t$ in~\eqref{eq:algo_general} is chosen to satisfy $\bL \eta_t m < 1$ and that $f^{\star} > -\infty$.
    Taking total expectation in \eqref{eq:lemma10} and summing up these inequalities for $t \leq T$, combined with the fact that $f(\bx_{T + 1}) \geq f^{\star}$, we deduce that
    % Then, it holds that
    \begin{align*}
        \Exp\left[\norm{\nabla f(\bx_S)}^2\right] \leq \frac{2(\Exp[f(\bx_1)] - f^{\star}) + \sum_{t = 1}^T\eta_t\parent{b_t^2 + \bL\eta_t v_t}}{\sum_{t = 1}^T\eta_t\parent{1 - \bL\eta_t{\sf V}_1}}\enspace.
    \end{align*}
\end{proof}

\begin{proof}[Proof of Theorem~\ref{nonc}]
The proof will be split into two parts: for gradient estimators \eqref{eq:grad_l2} and \eqref{eq:grad_l1}, respectively. 
Both of these proofs follow from Lemma \ref{HighSF}, which states that
\begin{align}
    \label{eq:proof_only_smooth}
        \Exp\left[\norm{\nabla f(\bx_S)}^2\right] \leq \frac{2\delta_1+ \sum_{t = 1}^T\eta_t\parent{b_t^2 + \bL\eta_t v_t}}{\sum_{t = 1}^T\eta_t\parent{1 - \bL\eta_t{\sf V}_1}}\enspace,
\end{align}
where $\delta_1 = \Exp[f(\bx_1)] - f^{\star}$.
Using the corresponding bounds on the bias $b_t$ and variance $v_t$, we substitute these values in the above inequality with $\eta_t$ and $h_t$ obtained by optimizing the obtained expressions. 

We start with the part of the proof that is common for both gradient estimators.
Introduce the notation
\begin{align*}
    \Xi_T := d^{-\frac{2(\beta-1)}{2\beta-1}}T^{-\frac{\beta}{2\beta-1}}\enspace.
\end{align*}
Using this notation, we consider  algorithm~\eqref{eq:algo_general} with gradient estimators \eqref{eq:grad_l2} or \eqref{eq:grad_l1} such that
\begin{align*}
    \eta_t =\min\left(\frac{\mathfrak{y}}{d},\, \Xi_T\right) \qquad\text{and}\qquad 
    h_t = \mathfrak{h} T^{-\frac{1}{2(2\beta-1)}} \enspace,
\end{align*}
where 
\begin{align*}
    %\label{eq:constants_only_smoothness_proof}
    (\mathfrak{y}, \mathfrak{h})
    % (\mathfrak{h}, \mathfrak{n})
    = 
    \begin{cases}
    \bigg((8\kappa\bL)^{-1},\,d^{\frac{1}{2\beta - 1}}\bigg)
    &\qquad\text{for estimator \eqref{eq:grad_l2}}\\
    % \\
    \bigg((72\kappa \bar{L})^{-1},\, d^{\frac{2\beta + 1}{4\beta - 2}}\bigg)
    &\qquad\text{for estimator \eqref{eq:grad_l1}}
    \end{cases}\enspace.
\end{align*}
Given the values of ${\sf V}_1$ in Table \ref{table:par},
%Furthermore, in the notation of Lemma \ref{HighSF} and, in particular, of Eq.~\eqref{eq:proof_only_smooth} above, the bounds in Lemma~\ref{lem:var_sphere_v2} and Lemma~\ref{lem:var_l1_v2} imply that 
the choice of $\eta_t$ for both algorithms ensures that
\begin{align*}
   \frac{1}{2} \leq 1 - \bL\eta_t {\sf V}_1\enspace.
\end{align*}
Thus we get from~\eqref{eq:proof_only_smooth} that both algorithms satisfy
\begin{align}
    \label{eq:proof_only_smooth1}
        \Exp\left[\norm{\nabla f(\bx_S)}^2\right] \leq
        \parent{\sum_{t = 1}^T\eta_t}^{-1}\left(4\delta_1  + 2\sum_{t = 1}^T\eta_tb_t^2 + 2\bL \sum_{t = 1}^T\eta_t^2v_t\right)\enspace.
\end{align}
Furthermore, since $\eta_t = \min (\mathfrak{y} / d,\, \Xi_T)$, then in both cases we have
\begin{align*}
     \parent{\sum_{t = 1}^T\eta_t}^{-1} = \max\parent{\frac{d}{T\mathfrak{y}},\, \frac{1}{T\Xi_T}} \leq \frac{d}{T\mathfrak{y}} + \frac{1}{T\Xi_T}\enspace.
\end{align*}
Using this bound in~\eqref{eq:proof_only_smooth1} we deduce that
\begin{align*}
        \Exp\left[\norm{\nabla f(\bx_S)}^2\right] \leq
        \parent{\frac{d}{T\mathfrak{y}} + \frac{1}{T\Xi_T}}\left(4\delta_1  + 2\sum_{t = 1}^T\eta_tb_t^2 + 2\bL \sum_{t = 1}^T\eta_t^2v_t\right)\enspace.
\end{align*}
Finally, by the definition of $\eta_t$ we have $\eta_t \leq \Xi_T$ for all $t = 1, \ldots, T$, which yields 
\begin{align}
    \label{eq:proof_only_smooth2}
    \Exp\left[\norm{\nabla f(\bx_S)}^2\right] \leq
        \parent{\frac{d}{\mathfrak{y}} + \frac{1}{\Xi_T}}\frac{4\delta_1}{T} +
        2\parent{\frac{d\Xi_T}{T\mathfrak{y}} + \frac{1}{T}}\sum_{t = 1}^T\left\{b_t^2 + \bL \Xi_Tv_t\right\}\enspace.
\end{align}
In the rest of the proof, we use the algorithm specific bounds on $b_t$ and $v_t$ as well as the particular choice of $\mathfrak{y}$ and $\mathfrak{h}$ in order to get the final results.

% \begin{proof}[Proof of Theorem~\ref{nonc}]

\subsection{Bounds for the gradient estimator (\ref{eq:grad_l2}) - $\ell_2$ randomization}

Lemma \ref{lem:bias_sphere} for the bias and Lemma \ref{lem:var_sphere_v2} for the variance
imply that
% substituted into Eq.~\eqref{eq:proof_only_smooth} yield
% Thanks to Lemma \ref{lem:bias_sphere} and Lemma \ref{lem:var_sphere_v2}, we have
\begin{align*}
    &b_t^2 \leq \left(\frac{\kappa_\beta L}{(\ell-1)!}\right)^2h_{t}^{2(\beta-1)}\quad\text{and}\quad
    v_t = 4d\kappa\bL^2 h_t^{2}+\frac{d^2\sigma^2\kappa}{2h_t^2}\,,
    \quad\text{and}\quad
    {\sf V}_1 = 4d\kappa\enspace.
\end{align*}
Using these bounds in~\eqref{eq:proof_only_smooth2} we get
    \begin{equation}
        \label{noncal1}
   \begin{aligned}
        \Exp\left[\norm{\nabla f(\bx_S)}^2\right] &\leq
        \parent{\frac{d}{\mathfrak{y}} + \Xi_T^{-1}}\frac{4\delta_1}{T}\\
        &\phantom{\leq}+
        \parent{\frac{d\Xi_T}{T\mathfrak{y}} + \frac{1}{T}}\sum_{t = 1}^T\bigg\{
        \cst_3h_{t}^{2(\beta-1)}
        +
        \Xi_Td^2\big(\cst_4d^{-1}h_t^2{+}\cst_5h_t^{-2}\big) \bigg\}\\
        &\leq
        \left(\frac{ d}{\mathfrak{y}}+\Xi_T^{-1}\right)\frac{4\delta_1}{T} + \frac{d\Xi_T{+}1}{T}\sum_{t=1}^{T}\left\{\cst_6h_t^{2(\beta-1)}+\cst_7d^2\Xi_T\big(d^{-1}h_t^2+h_t^{-2}\big)\right\}
    \end{aligned}
    \end{equation}
    where $\cst_3 = \Big(\tfrac{\kappa_\beta L}{(\ell-1)!}\Big)^2$, $\cst_4 = 4\kappa \bL^3$, $\cst_5 = \frac{\kappa \sigma^2\bL}{2}$, and $\cst_6 = 2\cst_3\big({\mathfrak{y}}^{-1}+1\big)$, $\cst_7 = 2\big({\mathfrak{y}}^{-1}+1\big)\big(\cst_4 + \cst_5\big)$.
Since $h_t = h_T$ for $t = 1,\dots,T$, inequality~\eqref{noncal1}  has the form
\begin{align}\label{noncal10}
    \Exp\left[\norm{\nabla f(\bx_S)}^2\right] \leq \left(\frac{d}{\mathfrak{y}}{+}{\Xi_T^{-1} }\right)\frac{4\delta_1}{T} + (d\Xi_T{+}1)\Big(\cst_6h_T^{2(\beta-1)}+\cst_7d^2\Xi_T\big(d^{-1}h_T^2+h_T^{-2}\big)\Big).
\end{align}
After substituting the expressions for $\Xi_T$ and $h_T$ into the above bound, the right hand side of~\eqref{noncal10} reduces to
\begin{align*}
    % \Exp\norm{\nabla f(\bx_S)}^2
    % &\leq \left(\frac{64\kappa\bL d}{T}+\frac{4}{a_T T}\right)\delta_1 + (da_T+1)\left(\cst_6h_T^{2(\beta-1)}+\cst_7d^2a_T\left(d^{-1}h_T^2+h_T^{-2}\right)\right)
    % \\
    % &\leq
    \frac{4d}{T\mathfrak{y}}\delta_1 +\left\{4\delta_1 + \left(\left(\frac{d}{T^{\beta}}\right)^{\frac{1}{2\beta-1}}{+}1\right)\left(\cst_6+\cst_7\left(1+d^{\frac{5-2\beta}{2\beta-1}}T^{-\frac{2}{2\beta-1}}\right)\right)\right\}\left(\frac{d^2}{T}\right)^{\frac{\beta-1}{2\beta-1}}.
\end{align*}
To conclude, we note that the assumption $T \geq d^{\frac{1}{\beta}}$, implies that for all $\beta \geq 2$ we have  $d^{\tfrac{5-2\beta}{2\beta-1}}T^{-\frac{2}{2\beta-1}} \leq 1$ and $\left({d}/{T^{\beta}}\right)^{\frac{1}{2\beta-1}}\leq 1$. Therefore, the final bound takes the form
\begin{align*}
    \Exp\left[\norm{\nabla f(\bx_S)}^2\right]
    &\leq
    \frac{4d}{T\mathfrak{y}}\delta_1 +\Big(4\delta_1 + 2\Big(\cst_6+2\cst_7\Big)\Big)\Big(\frac{d^2}{T}\Big)^{\frac{\beta-1}{2\beta-1}}
    \leq
    \Big(\cst_1\delta_1 + \cst_2\Big)\left(\frac{d^2}{T}\right)^{\frac{\beta-1}{2\beta-1}},
\end{align*}
where $\cst_1 = {4}(\mathfrak{y}^{-1}+1)$ and $\cst_2 = 2\big(\cst_6+2\cst_7\big)$.

\subsection{Bounds for the gradient estimator~(\ref{eq:grad_l1}) - $\ell_1$ randomization}
Lemma \ref{lem:bias_simplex} for the the bias and the bound~\eqref{eq:var_l1_for_proofs} for the variance imply that
\begin{align*}
    b_t^2\leq ({\color{purple}c_{\beta}}\kappa_\beta \ell L)^2h_t^{2(\beta-1)}d^{1-\beta}
    ,\qquad
    v_t = 72\kappa \bar{L}^2h_t^2+ \frac{d^3\sigma^2\kappa}{h^2_t}
    ,\qquad\text{and}\qquad
    {\sf V}_1 = 36d\kappa\enspace,
\end{align*}
with $\ell = \floor{\beta}$. Using these bounds in~\eqref{eq:proof_only_smooth2} we get
\begin{equation}
\label{gre3}
\begin{aligned}
    \Exp\left[\norm{\nabla f(\bx_S)}^2\right]\leq {(d\Xi_T{+}1)}&\left(\cst_6d^{1-\beta}h_T^{2(\beta-1)}{+}\Xi_T\left(\cst_7h_T^2{+}\cst_8d^3h_T^{-2}\right)\right)\\
    &+\left({\color{black}\frac{d}{\mathfrak{y}}}{+}{\Xi_T^{-1} }\right)\frac{4\delta_1}{T}
    \enspace,
\end{aligned}
\end{equation}
where the constants are defined as
\begin{align*}
\cst_6 &= 2({\color{purple}c_{\beta}}\kappa_\beta \ell L)^2\big(\mathfrak{y}^{-1}+1\big),\quad
\cst_7 = 144\kappa\bar{L}^3\big(\mathfrak{y}^{-1}+1\big),\quad
\cst_8 = 2{\bL\sigma^2\kappa}\big(\mathfrak{y}^{-1}+1\big)\enspace.
\end{align*}
Substituting the expressions for $\Xi_T$ and $h_T$ in~\eqref{gre3}, we deduce that
\begin{align*}
    \Exp\left[\norm{\nabla f(\bx_S)}^2\right]\leq \frac{4 d}{T \mathfrak{y}}\delta_1+\left\{4\delta_1+\left(\left(\frac{d}{T^{\beta}}\right)^{\frac{1}{2\beta-1}} +1\right)\left(\cst_6+\cst_8+\cst_7\left(\frac{d^{5-2\beta}}{T^2}\right)^{\frac{1}{2\beta-1}}\right)\right\}\left(\frac{d^2}{T}\right)^{\frac{\beta-1}{2\beta-1}}\enspace.
\end{align*}
Finally, we assumed that $T\geq d^{\frac{1}{\beta}}$, which implies that both $\frac{d}{T^{\beta}}$ and $\frac{d^{5-2\beta}}{T^2}$ are less than or equal to one. Thus, we have
\begin{align*}
    \Exp\left[\norm{\nabla f(\bx_S)}^2\right]\leq \left(\cst_1\delta_1+\cst_2\right)\left(\frac{d^2}{T}\right)^{\frac{\beta-1}{2\beta-1}}\enspace,
\end{align*}
where $\cst_1 = 4(\mathfrak{y}^{-1}+1)$, and $\cst_2=2\left(\cst_6+\cst_7+\cst_8\right)$.
% \end{proof}

\end{proof}

{\color{purple}\begin{proof}[Proof of Theorem~\ref{thm:OSsigma=0}]
As in the proof of Theorem \ref{nonc} we use Lemma \ref{HighSF}, cf. \eqref{eq:proof_only_smooth}:
\begin{align*}
\Exp\left[\norm{\nabla f(\bx_S)}^2\right] \leq \frac{2\delta_1 + \sum_{t = 1}^T\eta_t(b_t^2 + \bL\eta_t v_t)}{\sum_{t = 1}^T\eta_t(1 - \bL\eta_t{\sf V}_1)}.
\end{align*}
From this inequality and the fact that, by assumption, $\eta_t = (2\bL{\sf V}_1)^{-1}$ we obtain:
\begin{align*}
\Exp\left[\norm{\nabla f(\bx_S)}^2\right] \leq 
\frac{8\bL{\sf V}_1\delta_1 + 2\sum_{t = 1}^T(b_t^2 + (2{\sf V}_1)^{-1} v_t)}{T}.
\end{align*}
Since the gradient estimators~\eqref{eq:grad_l2} and~\eqref{eq:grad_l1} satisfy Assumption~\ref{ass:grad_general} and we 
consider the case $\sigma = 0$, $\beta=2$, the values $b_t=bLh_t$ and $v_t={\sf V}_2\bL^2h_t^2$ can be made as small as possible by choosing $h_t$ small enough.  Thus, we can take $h_t$ sufficiently small to have $\sum_{t = 1}^T(b_t^2 + (2{\sf V}_1)^{-1} v_t) \leq \bL{\sf V}_1$. Under this choice of $h_t$,
\begin{align*}
\Exp\left[\norm{\nabla f(\bx_S)}^2\right] \leq (8\delta_1 + 2)\frac{\bL{\sf V}_1}{T}\,.
\end{align*}
Using the values of ${\sf V}_1$ for the gradient estimators~\eqref{eq:grad_l2} and~\eqref{eq:grad_l1} (see Table~\ref{table:par}) we obtain the result.
\end{proof}
}
\subsection{Upper bounds: Smoothness and $\alpha$-gradient dominance}

\begin{proof}[Proof of Theorem~\ref{thm:master_dominant}]For brevity, we write $\Exp_t[\cdot]$ in place of $\Exp[\cdot \mid \bx_t]$.
    Using Lipschitz continuity of $\nabla f$ \citep[see e.g.][Lemma 3.4]{bubeck2014}  and the definition of the algorithm in~\eqref{eq:algo_general} with $\Theta=\mathbb{R}^d$ we have
    \begin{align*}
        \Exp_t[f(\bx_{t+1})]
        &\leq
        f(\bx_t) - \eta_t \scalar{\nabla f(\bx_t)}{\Exp_t[\bg_t]} + \frac{\bL\eta_t^2}{2}\Exp_t\left[\norm{\bg_t}^2\right]\\
        &\leq
        f(\bx_t) - \eta_t \norm{\nabla f(\bx_t)}^2 + \eta_t\norm{\nabla f(\bx_t)}\norm{\Exp_t[\bg_t] - \nabla f(\bx_t)}
        +
        \frac{\bL\eta_t^2}{2}\Exp_t\left[\norm{\bg_t}^2 \right]\enspace.
    \end{align*}
    Next, invoking Assumption~\ref{ass:grad_general} on the bias and variance of $\bg_t$ and using the elementary inequality $2ab \leq a^2 + b^2$ we get that, for the iterative procedure~\eqref{eq:algo_general} with $\Theta=\mathbb{R}^d$,
\begin{align*}
    \delta_{t+1} \leq \delta_{t}-\frac{\eta_t}{2}(1-\bar{L}\eta_{t}{\sf V}_1)\Exp{[\norm{\nabla f(\bx_t)}^2]}+\frac{\eta_t}{2}\left(b^2L^2h_t^{2(\beta-1)}+\bar{L}\eta_t\left({\sf V}_2{\color{purple}\bL^2}h_t^{2}+{\sf V}_3\sigma^2h_t^{-2}\right)\right)\enspace,
\end{align*}
where $\delta_t = \Exp[f(\bx_t)-f^{\star}]$. Furthermore, our choice of the step size $\eta_t$ ensures that $1 - \bL\eta_t {\sf V}_1 \geq \tfrac{1}{2}$. Using this inequality and the fact that $f$ is $\alpha$-gradient dominant we deduce that
\begin{align}\label{eq:PL_general_1}
    \delta_{t+1} \leq \delta_{t}\left(1 - \frac{\eta_t\alpha}{2}\right)+\frac{\eta_t}{2}\left(b^2L^2h_t^{2(\beta-1)}+\bar{L}\eta_t\left({\sf V}_2{\color{purple}\bL^2}h_t^{2}+{\sf V}_3\sigma^2h_t^{-2}\right)\right)\enspace.
\end{align}
We now analyze this recursion according to the cases $T > T_0$ and $T \leq T_0$,
where $T_0: = \floor{\frac{8\bar{L}{\sf V}_1}{\alpha}}$ is the value of $t$, where $\eta_t$ switches its regime. 

\medskip

\paragraph{First case: $T > T_0$.}
In this case, the recursion \eqref{eq:PL_general_1} has two different regimes, depending on the value of $\eta_t$.
In the first regime, for any $t = T_0+1,\dots,T$, we have $\eta_t = \frac{4}{\alpha t}$ and \eqref{eq:PL_general_1} takes the form
\begin{align}\label{eq:PL_general_2}
    \delta_{t+1}&\leq \delta_t\left(1-\frac{2}{t}\right) + 2b^2L^2\cdot\frac{h_{t}^{2(\beta-1)}}{{\alpha t}} +\frac{8\bar{L}}{\alpha^2t^2}\Big({\sf V}_2{\color{purple}\bL^2}h_t^2+{\sf V}_3\sigma^2h_{t}^{-2}\Big)\enspace.
\end{align}
Additionally in this regime of $t$, we have $h_t = \big(\frac{4\bar{L}{\color{purple}\sigma^2}{\sf V}_3}{b^2{\color{purple}L^2} \alpha t}\big)^\frac{1}{2\beta}$. Using this expression for $h_t$ in~\eqref{eq:PL_general_2} we obtain that
\begin{equation}
\begin{aligned}
    \delta_{t+1}\leq \delta_t\left(1-\frac{2}{t}\right) &+ \cst_3\cdot \frac{1}{\alpha}\left(\frac{\bL{\color{purple}\sigma^2}}{\alpha}\right)^{1-\frac{1}{\beta}}{\sf V}_3\left(\frac{{\sf V}_3}{b^2{\color{purple}L^2}}\right)^{-\frac{1}{\beta}}t^{-\frac{2\beta-1}{\beta}} \\&\phantom{\leq}+\cst_4\cdot \frac{1}{\alpha}\left(\frac{\bL}{\alpha}\right)^{1+\frac{1}{\beta}}{\sf V}_2{\color{purple}\bL^2}\left(\frac{{\sf V}_3{\color{purple}\sigma^2}}{b^2{\color{purple}L^2}}\right)^{\frac{1}{\beta}} t^{-\frac{2\beta+1}{\beta}},
\end{aligned}
\end{equation}
where $\cst_3 = 2^{4-\frac{2}{\beta}}$, and $\cst_4 = 2^{3 +\frac{2}{\beta}}$. Applying Lemma \ref{lemma:bisev} to the above recursion we get
\begin{align}\label{eq:PL_general_3}
    \delta_{T}&\leq \frac{2T_0}{T}\delta_{T_{0}+1} + \frac{\beta\cst_3}{(\beta{+}1)\alpha}\cdot {\sf V}_3\left(\frac{{\sf V}_3}{b^2{\color{purple}L^2}}\right)^{-\frac{1}{\beta}}\hspace{-.2truecm}\left(\frac{\alpha T}{\bL{\color{purple}\sigma^2}}\right)^{-\frac{\beta-1}{\beta}} 
    \\ \nonumber &\quad +\frac{\beta\cst_4}{(3\beta{+}1)\alpha}\cdot {\sf V}_2{\color{purple}\bL^2}\left(\frac{{\sf V}_3{\color{purple}\sigma^2}}{b^2{\color{purple}L^2}}\right)^{\frac{1}{\beta}} \hspace{-.2truecm}\left(\frac{\alpha T}{\bL}\right)^{-\frac{\beta+1}{\beta}}.
    %\enspace.
\end{align}

% \paragraph{Case 2: $t \in [1, T_0]$.}
{\color{purple}If $T_0 = 0$, we conclude the proof for the case $T>T_0$}.  Otherwise, we consider the second regime that corresponds to $t \in [1, T_0]$. In this regime, we have $h_t = \big(\frac{4\bar{L}{\color{purple}\sigma^2}{\sf V}_3}{b^2{\color{purple}L^2}\alpha T}\big)^\frac{1}{2\beta}$, $\eta_t = \frac{1}{2\bar{L}{\sf V}_1}$, and $\frac{4}{(T_0+1)\alpha}\leq \eta_t \leq \frac{4}{T_0\alpha}$. Using these expressions for $h_t$ and $\eta_t$ in~\eqref{eq:PL_general_1} we get that, for $1\leq t\leq T_0$,
\begin{align*}
    \delta_{t+1} \leq \delta_t\left(1 - \frac{2}{T_0+1}\right)&+\frac{2^{4-\frac{2}{\beta}}}{T_0}\cdot \frac{1}{\alpha}\left(\frac{\bL{\color{purple}\sigma^2}}{\alpha}\right)^{1-\frac{1}{\beta}}{\sf V}_3\left(\frac{{\sf V}_3}{b^2{\color{purple}L^2}}\right)^{-\frac{1}{\beta}}\left(T^{-\frac{\beta-1}{\beta}} + \frac{T^{\frac{1}{\beta}}}{T_0}\right)
    \\&\quad\quad+\frac{2^{3+\frac{2}{\beta}}}{T_0^2}\cdot \frac{1}{\alpha}\left(\frac{\bL}{\alpha}\right)^{1+\frac{1}{\beta}}{\sf V}_2{\color{purple}\bL^2}\left(\frac{{\sf V}_3{\color{purple}\sigma^2}}{b^2{\color{purple}L^2}}\right)^{\frac{1}{\beta}}T^{-\frac{1}{\beta}}\enspace.
\end{align*}
Using the rough bound $1 - \tfrac{2}{T_0+1}\leq 1$ and unfolding the above recursion we obtain: 
\begin{align*}
%\label{eq:PL_general_4}
     \delta_{T_0+1} &\leq \delta_1+2^{4-\frac{2}{\beta}}\cdot\frac{1}{\alpha}\left(\frac{\bL{\color{purple}\sigma^2}}{\alpha}\right)^{1-\frac{1}{\beta}} {\sf V}_3\left(\frac{{\sf V}_3}{b^2{\color{purple}L^2}}\right)^{-\frac{1}{\beta}}\left(T^{-\frac{\beta-1}{\beta}} + \frac{T^{\frac{1}{\beta}}}{T_0}\right)
 \\&\phantom{\leq}+\frac{2^{3+\frac{2}{\beta}}}{T_0}\cdot \frac{1}{\alpha}\left(\frac{\bL}{\alpha}\right)^{1+\frac{1}{\beta}}{\sf V}_2{\color{purple}\bL^2}\left(\frac{{\sf V}_3{\color{purple}\sigma^2}}{b^2{\color{purple}L^2}}\right)^{\frac{1}{\beta}}T^{-\frac{1}{\beta}}\enspace.
\end{align*}
Taking into account the definition of $T_0$, and the fact that $T_0 \leq T$ we further obtain:
\begin{align}\label{eq:PL_general_5}
    \frac{2T_0}{T}\delta_{T_0+1}&\leq \frac{16\bar{L}{\sf V}_1}{\alpha T}\delta_1 + \frac{2\cst_3}{\alpha}\cdot {\sf V}_3\left(\frac{{\sf V}_3}{b^2{\color{purple}L^2}}\right)^{-\frac{1}{\beta}}\left(\frac{\alpha T}{\bL{\color{purple}\sigma^2}}\right)^{-\frac{\beta-1}{\beta}} 
    \\ \nonumber &\quad+ \frac{2\cst_4}{\alpha}\cdot \left(\frac{\bL}{\alpha}\right)^{1+\frac{1}{\beta}}{\sf V}_2{\color{purple}\bL^2}\left(\frac{{\sf V}_3{\color{purple}\sigma^2}}{b^2{\color{purple}L^2}}\right)^{\frac{1}{\beta}}\left(\frac{\alpha T}{\bL}\right)^{-\frac{\beta+1}{\beta}}\enspace.
\end{align}
Finally, combining~\eqref{eq:PL_general_3} and~\eqref{eq:PL_general_5} yields:
\begin{align*}
    \delta_T \leq \cst_1\cdot\frac{\bL{\sf V}_1}{\alpha T}\delta_1 + \frac{\cst_2}{\alpha}\cdot\left({\sf V}_3\left(\frac{{\sf V}_3}{b^2{\color{purple}L^2}}\right)^{-\frac{1}{\beta}}+{\sf V}_2{\color{purple}\bL^2}\left(\frac{{\sf V}_3}{b^2{\color{purple}L^2}}\right)^{\frac{1}{\beta}}\left(\frac{\alpha T}{\bL{\color{purple}\sigma^2}}\right)^{-\frac{2}{\beta}}{\color{purple}\sigma^{-2}}\right)\left(\frac{\alpha T}{\bL{\color{purple}\sigma^2}}\right)^{-\frac{\beta-1}{\beta}}\enspace,
\end{align*}
where $\cst_1=16$ and $\cst_2 = \left(2 + \frac{\beta}{\beta+1}\right)\cst_3+\left(2 + \frac{\beta}{3\beta+1}\right)\cst_4$.

\medskip

\paragraph{Second case: $T \leq T_0$.} In this case, we have  $h_t = \big(\frac{4\bar{L}{\color{purple}\sigma^2}{\sf V}_3}{b^2{\color{purple}L^2}\alpha T}\big)^\frac{1}{2\beta}$ and thus~\eqref{eq:PL_general_1} takes the form
\begin{align*}
    \delta_{T+1}
    &\leq
    \delta_1 \left(1- \frac{2}{T_0+1}\right)^{T}+\frac{2^{4-\frac{2}{\beta}}}{T_0}\cdot \frac{1}{\alpha}\left(\frac{\bL{\color{purple}\sigma^2}}{\alpha}\right)^{1-\frac{1}{\beta}}{\sf V}_3\left(\frac{{\sf V}_3}{b^2{\color{purple}L^2}}\right)^{-\frac{1}{\beta}}\sum_{t=1}^{T}\left(T^{-\frac{\beta-1}{\beta}} + \frac{T^{\frac{1}{\beta}}}{T_0}\right)
    \\&\quad\quad+\frac{2^{3+\frac{2}{\beta}}}{T_0^2}\cdot \frac{1}{\alpha}\left(\frac{\bL}{\alpha}\right)^{1+\frac{1}{\beta}}{\sf V}_2{\color{purple}\bL^2}\left(\frac{{\sf V}_3{\color{purple}\sigma^2}}{b^2{\color{purple}L^2}}\right)^{\frac{1}{\beta}}\sum_{t=1}^{T}T^{-\frac{1}{\beta}}
    \\ &\leq \delta_1 \left(1- \frac{2}{T_0+1}\right)^{T} +\frac{\cst_3}{\alpha}\cdot {\sf V}_3\left(\frac{{\sf V}_3}{b^2{\color{purple}L^2}}\right)^{-\frac{1}{\beta}}\left(\frac{\alpha T}{\bL{\color{purple}\sigma^2}}\right)^{-\frac{\beta-1}{\beta}} +\frac{\cst_4}{\alpha}\cdot {\sf V}_2{\color{purple}\bL^2}\left(\frac{{\sf V}_3{\color{purple}\sigma^2}}{b^2{\color{purple}L^2}}\right)^{\frac{1}{\beta}}\left(\frac{\alpha T}{\bL}\right)^{-\frac{\beta+1}{\beta}}.
\end{align*}
Note that, for any $\rho, T> 0$, we have $(1 - \rho)^{T} \leq \exp(-\rho T) \leq \tfrac{1}{\rho T}$. Using this inequality for $\rho = \tfrac{2}{T_0 + 1}$, the definition of $T_0$ and the fact that $T+1\leq 2T$ we obtain: 
\begin{align*}
    \delta_{T+1} &\leq \cst_1\frac{\bL{\sf V}_1}{\alpha (T+1)}\delta_1 
    \\ & \quad + \frac{\cst_2}{\alpha}\left({\sf V}_3\left(\frac{{\sf V}_3}{b^2{\color{purple}L^2}}\right)^{-\frac{1}{\beta}}+{\sf V}_2{\color{purple}\bL^2}\left(\frac{{\sf V}_3}{b^2{\color{purple}L^2}}\right)^{\frac{1}{\beta}}\left(\frac{\alpha(T+1)}{\bL{\color{purple}\sigma^2}}\right)^{-\frac{2}{\beta}}{\color{purple}\sigma^{-2}}\right)\left(\frac{\alpha(T+1)}{\bL{\color{purple}\sigma^2}}\right)^{-\frac{\beta-1}{\beta}}.
\end{align*}
\end{proof}
\begin{proof}[Proof of Theorem~\ref{thm:PL-sigma=0}]
As in the proof of Theorem \ref{thm:master_dominant}, we consider separately the cases $T>T_0$ and   $T\leq T_0$, where $T_0:=\floor{\frac{8\bL {\sf V}_1}{\alpha}}$.
\paragraph{The case $T > T_0$.} First, consider the algorithm at steps $t= T_0,\dots, T$, where we  have $\eta_t = \frac{4}{\alpha t}$. Since $\sigma = 0$ and $\beta=2$, from~\eqref{eq:PL_general_1} we have
\begin{equation}
\label{eq:PL_sigma=0_1}
\begin{aligned}
    \delta_{t+1}&\leq \delta_t\left(1-\frac{2}{t}\right) + \left(\frac{2b^2L^2}{\alpha t}+\frac{8{\color{purple}\bL^3}}{\alpha^2 t^2}{\sf V}_2\right)h_t^2\enspace.
\end{aligned}
\end{equation}
{\color{purple}Since $\beta=2$, we have $L = \bL$}. Thus, using the assumption that $h_t \le \left(\frac{{\color{purple}\bL\vee 1}}{\alpha\wedge 1}{\color{purple}}T\left(2b^2\bL + \frac{8 {\color{purple}\bL^2}{\sf V}_2}{\alpha}\right)\right)^{-\frac{1}{2}}$  we deduce from~\eqref{eq:PL_sigma=0_1} that
\begin{align}
    \delta_{t+1} \leq \delta_t\left(1-\frac{2}{t}\right) + \frac{\bL}{\alpha t^2}\enspace.
\end{align}
Applying Lemma \ref{lemma:bisev} to the above recursion gives
\begin{align}\label{eq:PL_sigma=0_2}
    \delta_{T} \leq \frac{2T_0}{T}\delta_{T_{0}+1}  + \frac{\bL}{\alpha T}\enspace.
\end{align}
{\color{purple}If $T_0 = 0$, we conclude the proof for the case $T>T_0$}. Otherwise, we consider the algorithm at steps $t = 1,\dots,T_0$, where $\eta_t = \frac{1}{2\bL{\sf V}_1}$ and $\frac{4}{(T_0+1)\alpha}\leq \eta_t\leq \frac{4}{T_0 \alpha}$. From~\eqref{eq:PL_general_1} with $\sigma =0$ and $\beta=2$ we obtain
\begin{align}\label{eq:PL_sigma=0_4}
\delta_{t+1} \leq \delta_t\left(1-\frac{2}{T_0+1}\right) + \frac{\bL}{\alpha T_0}\left(2b^2{\color{purple}\bL} + \frac{8}{\alpha}{\color{purple}\bL^2}{\sf V}_2\right)h_t^{2}\enspace.
\end{align}
Using here the assumption that $h_t\le\left(\frac{{\color{purple}\bL\vee 1}}{\alpha\wedge 1}T\left(2b^2\bL + \frac{8 {\color{purple}\bL^2}{\sf V}_2}{\alpha}\right)\right)^{-\frac{1}{2}}$ and a rough bound $1 - \tfrac{2}{T_0+1}\leq 1$ and summing up both sides of the resulting inequality from $t=1$ to $t=T_0$ we get
\begin{align*}
\delta_{T_0+1} \leq \delta_1 + \frac{{\color{purple}\bL}(\alpha \wedge 1)}{\alpha(\bL\vee 1)T}\enspace.
\end{align*}
Combining this inequality with~\eqref{eq:PL_sigma=0_2} and using the definition of $T_0$ and the fact that $T_0 \leq T$ we obtain the bound 
%Eq.~\eqref{eq:PL_general_4} 
\begin{align*}
    \delta_T \leq \frac{16\bL{\sf V}_1}{\alpha T}\delta_1 + \frac{16\bL{\sf V}_1}{\alpha T^2} + \frac{\bL}{\alpha T}\enspace.
\end{align*}
It remains to note that ${\sf V}_1\leq 36d\kappa$, cf. Table \ref{table:par}. This implies the theorem for the case $T<T_0$ with $\cst_1 = 576\kappa$ and $\cst_2 =1/T + 1/\cst_1 d$.
\paragraph{Second case: $T \leq T_0$.} Using the fact that $h_t\le\left(\frac{T}{\alpha\wedge 1}\left(2b^2\bL + \frac{8 {\color{purple}\bL^2}{\sf V}_2}{\alpha}\right)\right)^{-\frac{1}{2}}$  and unfolding the recursion in~\eqref{eq:PL_sigma=0_4} gives
\begin{align*}
    \delta_T \leq \delta_1\left(1- \frac{2}{T_0+1}\right)^T + \frac{{\color{purple}\bL}}{\alpha T}\enspace.
\end{align*}
From the elementary inequality  $(1 - \rho)^{T}  \leq \tfrac{1}{\rho T}$, which is valid for all $\rho, T> 0$, we obtain with $\rho = \tfrac{2}{T_0 + 1}$ that
\begin{align*}
\delta_T \leq \frac{T_0+1}{2T}\delta_1 + \frac{\bL}{\alpha T} \leq 
\frac{T_0}{T}\delta_1 + \frac{\bL}{\alpha T}
\le
\cst_1\frac{\bL d}{{\alpha}T}\left(\delta_1+\cst_2\right)\enspace,
\end{align*}
where $\cst_1 = 288\kappa$, $\cst_2 = {\color{purple}1}/\cst_1d$, and the last inequality follows from the facts that $T_0 \le \frac{8\bL {\sf V}_1}{\alpha}$ and ${\sf V}_1\leq 36d\kappa$, cf. Table \ref{table:par}. 
\end{proof}

\subsection{Smoothness and $\alpha$-strong convexity: unconstrained minimization}

We will use the following basic lemma.

\begin{lemma}\label{scuntunc_gammatrial}
    Consider the iterative algorithm defined in~\eqref{eq:algo_general}.
    Let $f$ be $\alpha$-strongly convex on $\mathbb{R}^d$, let Assumption~\ref{ass:grad_general} be satisfied. {\color{brown}Let the minimizer $\bx^\star$ of $f$ on $\com$ be such that $\nabla f(\bx^\star)=0$.}
    % Assume that there exist $b, {\sf V}_1, {\sf V}_2, {\sf V}_3>0$ such that for all $t \geq 1$ it holds almost surely that
    % \begin{align*}
    %     \norm{\Exp[\bg_t \mid \bx_t] - \nabla f(\bx_t)} \leq bLh_t^{\beta-1}\qquad\text{and}\qquad \Exp\norm{\bg_t}^2  \leq {\sf V}_1  \Exp\norm{\nabla f(\bx_t)}^2 + {\sf V}_2\bL^2h_t^2 + {\sf V}_3\sigma^2h_t^{-2}\enspace.
    % \end{align*}
    Then we have
    \begin{small}
    \begin{align}\label{simpl_st}
       \Exp[f(\bx_t)-f^\star]\leq\frac{r_t-r_{t+1}}{2\eta_t}-r_t\left(\frac{\alpha}{4}-\frac{\eta_t}{2}\bL^2{\sf V}_1\right)+\frac{(bLh_t^{\beta-1})^2}{\alpha}+\frac{\eta_t}{2}\left({\sf V}_2\bar{L}^2h_t^{2} + {\sf V}_3\sigma^2h_t^{-2}\right)\enspace,
    \end{align}
    \end{small}
    where $r_t = \Exp[\|\bx_t - \bx^\star\|^2]$.
\end{lemma}
\begin{proof}
Recall the notation $\Exp_t[\cdot] = \Exp[\cdot \mid \bx_t]$.
For any $\bx \in \com$, by the definition of projection, 
%we have 
\begin{align}\label{nl1}
    \norm{\bx_{t+1}-\bx}^2 = \norm{\proj_{\com}\big(\bx_{t} - \eta_t\bg_t\big) - \bx}^2 \leq \norm{\bx_t - \eta_t\bg_{t} - \bx}^2.
\end{align}
%where the above inequality is obtained from the definition of Euclidean projection on the set $\com$ and the fact that $\bx \in \com$. 
Expanding the squares and rearranging the above inequality, we deduce that~(\ref{nl1}) is equivalent to 
\begin{align}\label{nl2}
    \langle \bg_t, \bx_t - \bx \rangle \leq \frac{\norm{\bx_t - \bx}^2-\norm{\bx_{t+1}-\bx}^2}{2\eta_t}+\frac{\eta_t}{2}\norm{\bg_t}^2\enspace.
\end{align}
On the other hand, since $f$ is a $\alpha$-strongly convex function on $\com$, we have
\begin{align}\label{nl3}
    f(\bx_t) - f(\bx) \leq \langle \nabla f(\bx_t), \bx_t -\bx \rangle -\frac{\alpha}{2}\norm{\bx_t -\bx}^2\enspace.
\end{align}
Combining~(\ref{nl2}) with~(\ref{nl3}) and introducing the notation $a_t = \norm{\bx_t - \bx^\star}^2$ we deduce that
\begin{align}
\label{nl40}
    \Exp_t[f(\bx_t) {-} f(\bx^\star)] &\leq \norm{\Exp_t[\bg_t] {-} \nabla f(\bx_t)}\norm{\bx_t{-}\bx^\star} {+} \frac{1}{2\eta_t}\Exp_t[a_t{-}a_{t+1}]+\frac{\eta_t}{2}\Exp_t\norm{\bg_t}^{2}{-}\frac{\alpha}{2}\Exp_t[a_t]\\\nonumber&\leq
    b L h_t^{\beta-1}\norm{\bx_t-\bx^\star}+\frac{1}{2\eta_t}\Exp_t[a_{t} - a_{t+1}]\\\label{nl4}
    &\phantom{\leq}
    + \frac{\eta_t}{2}\Big({\sf V}_1  \Exp[\norm{\nabla f(\bx_t)}^2] + {\sf V}_2\bL^2h_t^2 + {\sf V}_3\sigma^2h_t^{-2}\Big) - \frac{\alpha}{2}\Exp_t[a_t]\enspace.
\end{align}
Using the elementary inequality $2ab\leq a^2+b^2$ we have
\begin{align}\label{nl5}
    bLh_t^{\beta-1}\norm{\bx_t - \bx^\star} \leq \frac{(b L h_t^{\beta-1})^2}{\alpha}+\frac{\alpha}{4}\norm{\bx_t-\bx^\star}^2\enspace.
\end{align}
Substituting~(\ref{nl5}) in~(\ref{nl4}), setting $r_t = \Exp[a_t]$, {\color{brown}using the fact that $$\norm{\nabla f(\bx_t)}^2\le \bL^2\norm{\bx_t - \bx^\star}^2= \bL^2 a_t^2$$} and taking the total expectation from both sides of the resulting inequality 
yields the lemma.
%\begin{align*}
 %   \Exp[f(\bx_t)-f(\bx)]\leq\frac{r_t-r_{t+1}}{2\eta_t}-r_t\left(\frac{\alpha}{4}-\frac{\eta_t}{2}\bL^2{\sf V}_1\right)+\frac{(bLh_t^{\beta-1})^2}{\alpha}+\frac{\eta_t}{2}\left({\sf V}_2\bar{L}^2h_t^{2} + {\sf V}_3\sigma^2h_t^{-2}\right)\enspace,
%\end{align*}
%where the last display is obtained from $\alpha$ strong convexity of $f$. 
\end{proof}

\begin{proof}[Proof of Theorem \ref{thm:master_str_convex_uncostrained}]
By definition, $\eta_t \leq  \frac{\alpha}{4\bL^2 {\sf V}_1}$, so that $\frac{\alpha}{4} - \frac{\eta_t}{2}\bL^2 {\sf V}_1\geq\frac{\alpha}{8}$ and \eqref{simpl_st} implies that
\begin{align}\label{generalscuns}
    \Exp[f(\bx_t) - f^{\star}] \leq \frac{r_t - r_{t+1}}{2\eta_t} -\frac{\alpha}{8}r_t+\frac{(bLh_t^{\beta-1})^2}{\alpha}+\frac{\eta_t}{2}({\sf V}_2\bar{L}^2h_t^{2} + {\sf V}_3\sigma^2h_t^{-2})\enspace.
\end{align}
 Set {\color{purple}$T_0: = \max\left(\floor{\frac{32\bL^2{\sf V}_1}{\alpha^2}-1},0\right)$}. This is the value of $t$, where $\eta_t$ switches its regime. 
 We analyze the recursion \eqref{generalscuns} separately for the cases $T > T_0$ and $T \leq T_0$. We will use the fact that, by the convexity of $f$ and Jensen's inequality, 
\begin{equation}\label{eq:strongly_convex_l2_1}
\begin{aligned}
    f(\hat{\bx}_T) - f^{\star} \leq \frac{2}{T(T+1)}\sum_{t = 1}^T{\color{purple}t}(f(\bx_t) - f^\star). %&= \frac{{\color{purple}2}}{T(T+1)}\underbrace{\sum_{t = T_0 + 1}^{T}{\color{purple}t}(f(\bx_t) - f^\star)}_{\text{case 1}} \\&\phantom{\leq}+ \frac{2}{T(T+1)}\underbrace{\sum_{t = 1}^{T_0}{\color{purple}t}(f(\bx_t) - f^\star)}_{\text{case 2}} \enspace.
\end{aligned}
\end{equation}

\paragraph{First case: $T > T_0$.}
In this case, we decompose the sum in \eqref{eq:strongly_convex_l2_1} into the sum over  $t\in [T_0 +1,  T]$ and the sum over $t\in [1,  T_0]$. 

We first evaluate the sum $\sum_{t=T_0+1}^{T}{\color{purple}t}\Exp[f(\bx_t) - f^{\star}]$.
For any $t\in [T_0 +1,  T]$, we have {\color{purple}$\eta_t = \frac{8}{\alpha (t+1)}$ and $h_t = \left(\frac{4{\color{purple}\sigma^2}{\sf V}_3}{b^2{\color{purple}L^2}t}\right)^{\frac{1}{2\beta}}$}. Using in~\eqref{generalscuns} these values of $\eta_t$ and $h_t$, multiplying by $t$, and summing both sides of the resulting inequality from $T_0+1$ to $T$  we deduce that
\begin{align*}
    \sum_{t=T_0+1}^{T}{\color{purple}t\Exp[f(\bx_t) - f^{\star}] \leq
    \underbrace{\frac{\alpha}{16}\sum_{t=T_0+1}^{T}{\color{purple}t}\left( \left(r_t -r_{t+1}\right)(t+1) -2r_t\right)}_{=:{\sf I}}
    }&+
    \underbrace{\frac{\cst_4}{\alpha}(b{\color{purple}L})^{\frac{2}{\beta}}({\sf V}_3{\color{purple}\sigma^2})^{\frac{\beta-1}{\beta}}\sum_{t=T_0+1}^{T}t^{\frac{1}{\beta}}}_{=:{\sf II}}\\
    &+
    \underbrace{ \frac{\cst_5}{\alpha}{\color{purple}\bL^2}{\sf V}_2\left(\frac{{\sf V}_3{\color{purple}\sigma^2}}{b^2{\color{purple}L^2}}\right)^{\frac{1}{\beta}}\sum_{t=T_0+1}^{T}t^{-\frac{1}{\beta}}}_{=:{\sf III}}\enspace,
\end{align*}
where $\cst_4 = 2^{\frac{3\beta-2}{\beta}}$, $\cst_5 = 2^{\frac{2\beta+2}{\beta}}$ and we defined the terms ${\sf I}, {\sf II}$, and ${\sf III}$ that will be evaluated separately.

It is not hard to check that ${\sf I} \leq \tfrac{\alpha}{16}{\color{purple}T} T_0 r_{T_0+1}$ since the summation in term ${\sf I}$ is telescoping. Next, we have 
$${\sf II}
    \leq \frac{\cst_4}{\alpha}(b{\color{purple}L})^{\frac{2}{\beta}}({\sf V}_3{\color{purple}\sigma^2})^{\frac{\beta-1}{\beta}}\sum_{t=1}^{T}t^{\frac{1}{\beta}}
    \leq \frac{\cst_6}{\alpha}(b{\color{purple}L})^{\frac{2}{\beta}}({\sf V}_3{\color{purple}\sigma^2})^{\frac{\beta-1}{\beta}}T^{\frac{1}{\beta}+1}.
    $$
% \begin{align*}
%     {\sf II}
%     \leq \frac{\cst_4}{\alpha}b^{\frac{2}{\beta}}{\sf V}_3^{\frac{\beta-1}{\beta}}\sum_{t=1}^{T}t^{-\frac{\beta-1}{\beta}}
%     \leq \frac{\beta\cst_4}{\alpha}b^{\frac{2}{\beta}}{\sf V}_3^{\frac{\beta-1}{\beta}}T^{\frac{1}{\beta}}\enspace.
% \end{align*}
Finally, 
\begin{align*}
    {\sf III}
    %\leq \frac{\cst_5}{\alpha }{\color{purple}\bL^2}{\sf V}_2\left(\frac{{\sf V}_3{\color{purple}\sigma^2}}{b^2{\color{purple}L^2}}\right)^{\frac{1}{\beta}}\sum_{t=T_0+1}^{T}t^{-\frac{1}{\beta}} 
    % &\leq
    % \frac{4\cst_4}{\alpha }d^{1+\frac{2}{\beta}}\left(\parent{T_0 + 1}^{-\frac{1}{\beta}-1} + \int_{T_0}^{T} t^{-\frac{1}{\beta}-1} \d t\right)\\
    % &\leq \frac{4\cst_4}{\alpha }d^{1+\frac{2}{\beta}}\left(\parent{T_0 + 1}^{-\frac{1}{\beta}-1} + \beta {T_0}^{-\frac{1}{\beta}}\right)\\
    % &\leq \frac{4(\beta+1)\cst_4}{\alpha }d^{1+\frac{2}{\beta}} {T_0}^{-\frac{1}{\beta}}\\
    % &\leq \frac{4(\beta+1)\cst_4}{\alpha }d^{1+\frac{1}{\beta}} \frac{\alpha^{\frac{2}{\beta}}}{(64\kappa\bL^2)^{\frac{1}{\beta}}}
    % \leq \frac{\cst_{6}}{\alpha}d^{1+\frac{1}{\beta}}
    \leq \frac{\cst_5}{\alpha }{\color{purple}\bL^2}{\sf V}_2\left(\frac{{\sf V}_3{\color{purple}\sigma^2}}{b^2{\color{purple}L^2}}\right)^{\frac{1}{\beta}}\sum_{t=1}^{T}t^{-\frac{1}{\beta}} 
    \leq \frac{\cst_{7}}{\alpha}{\color{purple}\bL^2}{\sf V}_2\left(\frac{{\sf V}_3{\color{purple}\sigma^2}}{b^2{\color{purple}L^2}}\right)^{\frac{1}{\beta}}T^{1-\frac{1}{\beta}}\enspace,
\end{align*}
where $\cst_6 = \frac{\beta+1}{\beta}\cst_4$ and $\cst_7 = \frac{\beta-1}{\beta}\,\cst_5$. Combining these bounds on ${\sf I}, {\sf II}$, and ${\sf III}$ we obtain
\begin{align}\label{sfar1}
    \sum_{t=T_0+1}^{T}{\color{purple}t}\Exp[f(\bx_t)-f^{\star}] & \leq \frac{\alpha}{16}{\color{purple}T}T_0r_{T_0+1} 
    \\ \nonumber
    & \quad +  \Big(\cst_6(b{\color{purple}L})^{\frac{2}{\beta}}({\sf V}_3{\color{purple}\sigma^2})^{\frac{\beta-1}{\beta}}+\cst_{7}{\sf V}_2{\color{purple}\bL^2}\left(\frac{{\sf V}_3{\color{purple}\sigma^2}}{b^2{\color{purple}L^2}}\right)^{\frac{1}{\beta}}{\color{purple}T^{-\frac{2}{\beta}}}\Big)\frac{{\color{purple}T^{\frac{1}{\beta}+1}}}{\alpha}\,.
\end{align}
If {\color{purple}$T_0=0$} then combining \eqref{eq:strongly_convex_l2_1} and  \eqref{sfar1} proves the theorem. If $T_0\ge 1$ we need additionally to control the value $r_{T_0+1}$ on the right hand side of \eqref{sfar1}. It follows from~\eqref{generalscuns} that, for $1 \leq t \leq T_0$, 
\begin{align*}
    r_{t+1} \leq r_t + 2\eta_t\frac{(bLh_t^{\beta-1})^2}{\alpha}+\eta_t^2\left({\sf V}_2\bar{L}^2h_t^{2} + {\sf V}_3\sigma^2h_t^{-2}\right).
\end{align*}
Moreover, for $1\leq t\leq T_0$ we have $\eta_t = \frac{\alpha}{4\bL^2 {\sf V}_1}$ and $\eta_t\leq \frac{8}{\alpha (T_0+1)}$. Therefore, unfolding the above recursion we get
\begin{align*}
    r_{T_0+1}\leq r_1 + \sum_{t=1}^{T_0}\left(\frac{16}{ \alpha^2T_0}\left(bLh_t^{\beta-1}\right)^2 + \frac{64}{\alpha^2 T_0^2}\left({\sf V}_2\bar{L}^2h_t^{2} + {\sf V}_3\sigma^2h_t^{-2}\right)\right)\,.
\end{align*}
For $1\leq t\leq T_0$ we have $h_t = \left(\frac{{4\color{purple}\sigma^2}{\sf V}_3}{b^2{\color{purple}L^2}T}\right)^{\frac{1}{2\beta}}$, which yields
\begin{align*}
    r_{T_0+1} \leq r_1 + 16\left(\cst_4(b{\color{purple}L})^{\frac{2}{\beta}}({\sf V}_3{\color{purple}\sigma^2})^{\frac{\beta-1}{\beta}}+\cst_5{\sf V}_2{\color{purple}\bL^2}\left(\frac{{\sf V}_3{\color{purple}\sigma^2}}{b^2{\color{purple}L^2}}\right)^{\frac{1}{\beta}}T^{-\frac{2}{\beta}}\right)\frac{T^{\frac{1}{\beta}}}{\alpha^2 T_0}\,,
\end{align*}
so that
\begin{align}\label{sfar2}
    \frac{\alpha}{16}{\color{purple}T} T_0 r_{T_0+1}\leq \frac{2\bL^2T{\sf V}_1}{\alpha}r_1 + \left(\cst_4(b{\color{purple}L})^{\frac{2}{\beta}}({\sf V}_3{\color{purple}\sigma^2})^{\frac{\beta-1}{\beta}}+\cst_5{\sf V}_2{\color{purple}\bL^2}\left(\frac{{\sf V}_3{\color{purple}\sigma^2}}{b^2{\color{purple}L^2}}\right)^{\frac{1}{\beta}}T^{-\frac{2}{\beta}}\right)\frac{{\color{purple}T^{\frac{1}{\beta}+1}}}{\alpha}\,.
\end{align}
It follows from~\eqref{sfar1} and~\eqref{sfar2} that
\begin{equation}
\label{sfar3}
\begin{aligned}
    \sum_{t=T_0+1}^{T}{\color{purple}t}\Exp[f(\bx_t)-f^{\star}]
    \leq
    \frac{2\bL^2{\color{purple}T}{\sf V}_1}{\alpha}r_1
    % \\
      +
    \Bigg\{&{\color{purple}(\cst_4+\cst_6)}(b{\color{purple}L})^{\frac{2}{\beta}}({\sf V}_3{\color{purple}\sigma^2})^{\frac{\beta-1}{\beta}}\\
    &\quad+
    {\color{purple}(\cst_5+\cst_7)}{\sf V}_2{\color{purple}\bL^2}\left(\frac{{\sf V}_3{\color{purple}\sigma^2}}{b^2{\color{purple}L^2}}\right)^{\frac{1}{\beta}}T^{-\frac{2}{\beta}}\Bigg\}\frac{{\color{purple}T^{\frac{1}{\beta}+1}}}{\alpha}\, ,    
\end{aligned}
\end{equation}
We now evaluate the sum $\sum_{t=1}^{T_0} t\Exp[f(\bx_t) -f^{\star}]$. Recall that for $t \in [1, T_0]$ the parameters $h_t,\eta_t$ take constant values: {\color{purple}$h_t = \left(\frac{4{\color{purple}\sigma^2}{\sf V}_3}{b^2{\color{purple}L^2}T}\right)^{\frac{1}{2\beta}}$} and $\eta_t = \frac{\alpha}{4\bL^2 {\sf V}_1}$. Omitting in \eqref{generalscuns} the term $-\alpha r_t/8$, summing the resulting recursion from $1$ to $T_0$ and using the inequality {\color{purple}$\eta_t\leq \frac{8}{\alpha (T_0+1)}$} we obtain
{\color{brown}
\begin{align} \label{sfar4}
    \sum_{t=1}^{T_0} {\color{purple}t}\Exp[f(\bx_t) -f^{\star}] 
    &
  \leq   T \sum_{t=1}^{T_0}\Exp[f(\bx_t) - f^{\star}] 
    \\
    &
     \le
     \nonumber
    \frac{T r_1}{2\eta_1}
    +T  \sum_{t=1}^{T_0}\Big(\frac{(bLh_t^{\beta-1})^2}{\alpha}+\frac{\eta_t}{2}({\sf V}_2\bar{L}^2h_t^{2} + {\sf V}_3\sigma^2h_t^{-2})\Big)
    \\
    \nonumber
    &
    \le
    \frac{T r_1}{2\eta_1}
    +T^2\frac{(bLh_1^{\beta-1})^2}{\alpha}+\frac{2T}{\alpha}({\sf V}_2\bar{L}^2h_1^{2} + {\sf V}_3\sigma^2h_1^{-2})
    \\
    \nonumber
    &
    =
    \frac{2\bL^2{\color{purple}T}{\sf V}_1}{\alpha}r_1 + \Big(\cst_4(b{\color{purple}L})^{\frac{2}{\beta}}({\sf V}_3{\color{purple}\sigma^2})^{\frac{\beta-1}{\beta}}+ \cst_5{\sf V}_2{\color{purple}\bL^2}\Big(\frac{{\sf V}_3{\color{purple}\sigma^2}}{ b^2{\color{purple}L^2}}\Big)^{\frac{1}{\beta}}T^{-\frac{2}{\beta}}\Big)\frac{{\color{purple}T^{\frac{1}{\beta}+1}}}{\alpha}\,.
\end{align}
}
Summing up~\eqref{sfar3} and~\eqref{sfar4} and using~\eqref{eq:strongly_convex_l2_1} we obtain the bound of the theorem:  
\begin{align*}
    \Exp[f(\hat{\bx}_T)-f^{\star}] \leq \cst_1\frac {{\color{purple}\bL^2}{\sf V}_1}{\alpha T}r_1+\left(\cst_2(b{\color{purple}L})^{\frac{2}{\beta}}({\sf V}_3{\color{purple}\sigma^2})^{\frac{\beta-1}{\beta}}+\cst_3{\sf V}_2{\color{purple}\bL^2}\left(\frac{{\sf V}_3{\color{purple}\sigma^2}}{b^2{\color{purple}L^2}}\right)^{\frac{1}{\beta}}T^{-\frac{2}{\beta}}\right)\frac{T^{-\frac{\beta-1}{\beta}}}{\alpha}\,,
\end{align*}
where ${\color{brown}\cst_1 = 8}$, {\color{purple}$\cst_2 = 4\cst_4+2\cst_6$, and $\cst_3 =4\cst_5+2\cst_7$}.

\paragraph{Second case: $T \leq T_0$.} {\color{purple}
In this scenario, the summation $\sum_{t=1}^{T}{\color{purple}t}\Exp[f(\bx_t) - f^{\star}]$ is treated as in \eqref{sfar4}, with the only difference that $T_0$ is replaced by $T$. As a result, we obtain the same bound as in \eqref{sfar4}.
}
\end{proof}

\subsection{Smoothness and $\alpha$-strong convexity: constrained minimization}
\begin{proof}[Proof of Theorem~\ref{scunt}]
% Let $\bx^\star_T$ be any minimizer of $\sum_{t = 1}^T f(\bx)$ over $\com$.
Since  $\sup_{\bx \in \com}\norm{\nabla f(\bx)} \leq G$ we get from~\eqref{nl4} that, for any $t = 1, \ldots, T$,
\begin{align}
    \label{eq:lemma43_0new}
 \Exp[f(\bx_t)- f^{\star}]\leq \frac{r_t-r_{t+1}}{2\eta_t}-\frac{\alpha}{4}r_t+\frac{\left(bLh_t^{\beta-1}\right)^2}{\alpha}+\frac{\eta_t}{2}\big({\sf V}_1G^2+{\sf V}_2\bL^2h_t^2+{\sf V}_3\sigma^2h_t^{-2}\big).
\end{align}
Multiplying both sides of \eqref{eq:lemma43_0new} by $t$, summing up from $t=1$ to $T$ and using the fact that 
$$
\sum_{t=1}^{T}\left(\frac{t(r_{t}-r_{t+1})}{2\eta_{t}} - \frac{\alpha}{4}t r_{t}\right) \le 0\qquad  \text{if}\ \ \ \eta_{t} = \frac{4}{\alpha (t+1)}
$$
we find that
\begin{align*}
 %   \label{eq:lemma43_1new}
 \sum_{t=1}^{T}t\Exp[f(\bx_t)- f^{\star}]
 \leq 
 \frac{1}{\alpha}\sum_{t=1}^{T}\left[t(bLh_t^{\beta-1})^2+\frac{2t}{t+1}\big({\sf V}_1G^2+{\sf V}_2\bL^2h_t^2+{\sf V}_3\sigma^2h_t^{-2}\big)\right].
\end{align*}
Since $h_t = \left(\frac{{\color{purple}\sigma^2}{\sf V}_3}{b^2{\color{purple}L^2} t}\right)^{\frac{1}{2\beta}}$ we obtain
\begin{align}
    \label{eq:lemma43_1new}
 \sum_{t=1}^{T}t\Exp[f(\bx_t)- f^{\star}]
 \leq 
 \frac{2{\sf V}_1 G^2 T}{\alpha } +
\frac{\cst_3}{\alpha}\left({\sf V}_3{\color{purple}\sigma^2}\left(\frac{{\sf V}_3{\color{purple}\sigma^2}}{b^2{\color{purple}L^2}}\right)^{-\frac{1}{\beta}} + {\sf V}_2{\color{purple}\bL^2}\left(\frac{{\sf V}_3{\color{purple}\sigma^2}}{b^2{\color{purple}L^2}}\right)^{\frac{1}{\beta}}T^{-\frac{2}{\beta}}\right)T^{1+\frac{1}{\beta}},
\end{align}
where $\cst_3 = {\color{purple}2}$. To complete the proof, we multiply both sides of \eqref{eq:lemma43_1new} by $\frac{2}{T(T+1)}$ and use \eqref{eq:strongly_convex_l2_1}. 
\end{proof}

\section{Proof of the lower bounds}

{\color{brown}Set for brevity $A_t = \{(\bz_i,y_i)_{i=1}^{t},(\bz'_i,y'_i)_{i=1}^{t}\}$ for $t\geq 1$. %and $A_0 = \{(\bz_1,y_1)_{i=1}^{t},(\bz'_i,y'_i)_{i=1}^{t}\}$. 
Without loss of generality, we assume that $\bz_1$ and $\bz'_1$ are fixed and we prove the result for any fixed
$(\bz_1,\bz'_1)$.
For any $f : \bbR^d \to \bbR$ and any sequential strategy in the class $\Pi_T$, such that
$\bz_t = \Phi_t(A_{t-1},\btau_t)$ with $y_t = f(\bz_t) + \xi_t$ and $\bz'_t = \Phi'_t(A_{t-1},\btau_t)$ for $t\ge 2$, with $y'_t = f(\bz'_t) + \xi'_t$ for $t = 1, \ldots, T$, 
we will denote by $\mathbf{P}_{f}$ the joint distribution of $(A_T, (\btau_i)_{i=2}^T)$. }

 We start with the following lemma that will be used in the proof of Theorem \ref{lb}.
\begin{lemma}
\label{lem:hellinger_lower}
Let Assumption \ref{ass:lower-bound} be satisfied. Then, for any functions $f, f' : \bbR^d \to \bbR$ such that $\norm{f-f'}_{\infty} =\max_{\bu\in\bbR^d}|f(\bu)-f'(\bu)| \leq v_0$ it holds that
\begin{align*}
     \frac{1}{2}H^2(\mathbf{P}_{f} , \mathbf{P}_{f'}) \leq 1-\left(1-\frac{I_0}{2}\norm{f-f'}_{\infty}^2\right)^{T}.
\end{align*}
\end{lemma}
\begin{proof}%[Proof of Lemma~\ref{lem:hellinger_lower}]
 {\color{brown}{Since, for each $t\ge 2$,the noise $(\xi_t,\xi'_t)$ is independent of $(A_{t-1},(\btau_i)_{i=1}^{t})$ and $\btau_t$ is independent of 
$A_{t-1}$ we have 
\begin{align*}
    \d\mathbf{P}_f = \d F_1(y_1 - f(\bz_1),y'_1 - f(\bz'_1))\prod_{t = 2}^T \d F_t\bigg(y_t - f\Big(\Phi_{t}(A_{t-1},\btau_t)\Big),y'_t - f\Big(\Phi'_{t}(A_{t-1},\btau_t)\Big) \bigg)\d\mathbb{P}_t(\btau_t)\,
\end{align*}
}
where $\mathbb{P}_t$ is the probability measure corresponding to the distribution of $\btau_t$. Set for brevity {\color{purple}$\d F_{f, 1} \triangleq \d F_1\big(y_1 - f(\bz_1),y'_1 - f(\bz'_1)\big)$ and $\d F_{f, t} \triangleq \d F_t\big(y_t - f(\bz_t),y'_t - f(\bz'_t)\big) \d\mathbb{P}_t(\btau_t)$, $t\ge 2$.}
% \begin{align*}
%     \d F_{f, t} \triangleq \d F\bigg(y_t - f\Big(\Phi_{t}(\bz_1, y_1, \ldots, y_{t - 1})\Big) \bigg)\enspace,
% \end{align*}
With this notation, we have $\d\mathbf{P}_f = \prod_{t = 1}^T \d F_{f, t}$.}
Using the definition of Hellinger distance we obtain
\begin{align*}
   1 {-} \frac{1}{2}H^2(\mathbf{P}_{f} , \mathbf{P}_{f'})
    &=
    \int\sqrt{\d\mathbf{P}_{f}\d\mathbf{P}_{f'}}
    =\prod_{t=1}^{T}\int \sqrt{\d F_{f, t}}\sqrt{\d F_{f', t}}
    =
    \prod_{t = 1}^T\parent{1 - \frac{H^2\parent{\d F_{f, t}, \d F_{f', t}}}{2}}.
    % \\
    % &=\Big(1- \prod_{i=1}^{T}\big(1 - \frac{H^2(\d F_{\omega, i}, \d F_{\omega', i})}{2}\big)\Big)\\&\leq \Big(1-\min_{1 \leq i \leq T}\big(1-\frac{H^{2}(\d F_{\omega,i},dF_{\omega^{'},i})}{2}\big)^{T}\Big)\\&\leq \Big(1-\min_{u \in \com^{'}}\big(1-\frac{I_{0}|f_{\omega}(u)-f_{\omega'}(u)|^2}{2}\big)^{T}\Big) \\&= \Big(1-\big(1-\max_{u\in\com^{'}}\frac{I_{0}|f_{\omega}(u)-f_{\omega'}(u)|^2}{2}\big)^{T}\Big)\\&\AK{\leq \Big(1 - \big(1-\frac{I_{0}r^{2}h^{2\beta}\eta(1/\sqrt{d})}{2}\big)\Big)}
\end{align*}
Finally, invoking Assumption \ref{ass:lower-bound}, we get
\begin{align*}
    \prod_{t = 1}^T\parent{1 - \frac{H^2\parent{\d F_{f, t}, \d F_{f', t}}}{2}}
    &\geq
    \min_{1 \leq t \leq T}\parent{1 - \frac{H^2\parent{\d F_{f, t},  \d F_{f', t}}}{2}}^T\\
    &\geq
    \parent{1 - \frac{I_0\norm{f-f'}_{\infty}^2}{2}}^T,
\end{align*}
which implies the lemma.
\end{proof}

\begin{proof}[Proof of Theorem~\ref{lb}] The proof follows the general lines given in~\cite{akhavan2020}, so that we omit some details that can be found in that paper.
We first assume that $\alpha\ge T^{-1/2+1/\beta}$.

Let $\eta_{0} : \mathbb{R} \to \mathbb{R}$ be an infinitely many times differentiable function such that
\begin{equation*}
   \eta_{0}(x) = \begin{cases}
      =1 & \text{if $|x|\leq 1/4$},\\
      \in (0,1) & \text{if $1/4 < |x| < 1$},\\
      =0 & \text{if $|x| \geq 1$}.
    \end{cases} 
\end{equation*}
Set $\eta(x) = \int_{-\infty}^{x} \eta_{0}(\tau)d\tau$. Let $\Omega = \big\{-1,1\big\}^{d}$ be the set of binary sequences of length $d$. 
Consider the finite set of functions $f_{\omega}: \mathbb{R}^{d}\to \mathbb{R}, \bomega = (\omega_1, \ldots, \omega_d)\in\Omega$, defined as follows:
\[
f_{\bomega}(\bu) = \alpha(1+\delta) \norm{\bu}^{2}/2 + \sum_{i=1}^{d}\omega_{i}rh^{\beta}\eta(u_{i}h^{-1}),
\qquad \bu=(u_1,\dots,u_d),
\]
where $\omega_i\in \{-1,1\}$, 
 $h =\min\big((\alpha^2/d)^{\frac{1}{2(\beta-1)}}, T^{-\frac{1}{2\beta}}\big)$ and $r>0, \delta >0$ are fixed numbers that will be chosen small enough.

 It is shown in \cite{akhavan2020} that if $\alpha\ge T^{-1/2+1/\beta}$ then $f_{\bomega}\in\mathcal{F}'_{\alpha,\beta}$
 for $r>0$ and $\delta >0$ small enough, and the minimizers of functions $f_{\bomega}$ belong to $\com$ and are of the form \[\bx_{\bomega}^{*} = (x^\star(\omega_1), \dots, x^\star(\omega_d))\enspace, \]
where $x^\star(\omega_i)=-\omega_{i}\alpha^{-1}(1+\delta)^{-1}r h^{\beta-1}$.

{\color{brown} For any fixed $\bomega\in\Omega$, we denote by $\mathbf{P}_{\bomega,T}$ the probability measure corresponding to the joint distribution of {$(A_T,(\btau_i)_{i=2}^T)$} where $y_{t}=f_{\bomega} (\bz_{t})+\xi_{t}$ and {$y'_{t}=f_{\bomega} (\bz'_{t})+\xi'_{t}$} with 
%independent identically distributed 
{$(\xi_{t},\xi'_{t})$'s satisfying Assumption \ref{ass:lower-bound}, 
%namely, \eqref{distribution} and the fact that $(\xi_t,\xi'_t)$ is independent of $(A_{t-1},\btau_t)$ for each $t$, 
and $(\bz_t,\bz'_t)$'s chosen by a sequential strategy in $\Pi_T$.}}
Consider the statistic
$$\hat{\bomega} \in \argmin_{\bomega \in \Omega} \norm{\tilde\bx_{T}-\bx^{*}_{\bomega}}\enspace.$$
Classical triangle inequality based arguments yield
\[
\max_{\bomega \in \Omega}\mathbf{E}_{\bomega,T}\big[\norm{\tilde\bx_{T}-\bx^{*}_{\bomega}}^{2}\big] \geq \alpha^{-2}r^{2} h^{2\beta-2}\inf_{\hat{\bomega}}\max_{\bomega \in \Omega}\mathbf{E}_{\bomega,T} \big[\rho(\hat{\bomega},\bomega) \big]\enspace.\]
Note that for all $\bomega, \bomega'\in \Omega$ such that $\rho(\bomega, \bomega')=1$ we have
\[
\max_{\bu\in\mathbb{R}^d}{|f_{\bomega}(\bu)-f_{\bomega'}(\bu)|}\leq 2rh^{\beta} \eta(1) \leq 2rT^{-1/2} \eta(1).
\]
{\color{brown}
Thus, choosing $r$ small enough to satisfy $2r\eta(1)< \min(v_{0}, I_0^{-1/2})$ we ensure $2rT^{-1/2} \eta(1) \leq v_0$ to apply Lemma~\ref{lem:hellinger_lower} and deduce for the considered $\bomega, \bomega' \in \Omega$ that
\begin{align*}
    H^2(\mathbf{P}_{\bomega,T} , \mathbf{P}_{\bomega',T})&\leq 2\Big(1 - \big(1 - (2T)^{-1}\big)^T\Big)
    \leq
    1,
\end{align*}
where we have used the fact that $1-x\ge 4^{-x}$ for $0<x\le 1/2$.
Applying \cite[Theorem 2.12]{Tsybakov09} we deduce that
\[
\inf_{\hat{\bomega}}\max_{\bomega \in \Omega}\mathbf{E}_{\bomega, T} [\rho(\hat{\bomega},\bomega)]\geq 0.3 \, d.
\]
Therefore, we have proved that if $\alpha\ge T^{-1/2+1/\beta}$ there exist $r>0$ and $\delta >0$ 
%small enough 
such that  
\begin{equation}\label{eq3:lb}
 \max_{\bomega \in \Omega}\mathbf{E}_{\bomega,T}\big[\norm{\tilde\bx_{T}-\bx^{*}_{\bomega}}^{2}\big]\geq 0.3 \, d\alpha^{-2}r^{2} h^{2\beta-2} = 0.3 \, r^{2}\min\Big(1, 
\,\frac{d}{\alpha^{2}}T^{-\frac{\beta-1}{\beta}}
\Big).   
\end{equation}
}This implies \eqref{eq2:lb} for $\alpha\ge T^{-1/2+1/\beta}$. 
In particular, if $\alpha=\alpha_0:= T^{-1/2+1/\beta}$ the bound \eqref{eq3:lb} is of the order $\min\Big(1, {d}{T^{-\frac{1}{\beta}}}
\Big)$. Then for $0<\alpha<\alpha_0$ we also have the bound of this order since the classes $\mathcal{F}_{\alpha,\beta}$ are nested: $\mathcal{F}_{\alpha_0,\beta}\subset \mathcal{F}_{\alpha,\beta}$. This completes the proof of \eqref{eq2:lb}.

We now prove \eqref{eq1:lb}. From \eqref{eq3:lb} and $\alpha$-strong convexity of $f$ we get that, for $\alpha\ge T^{-\frac{\beta+2}{2\beta}}$,
\begin{equation}\label{eq4:lb}
\max_{\bomega \in \Omega}\mathbf{E}_{\bomega,T}
\big[f(\tilde\bx_T)-f(x_{\bomega}^*)\big]\geq
{\color{brown}0.15} \, r^{2}\min\Big(\alpha, 
\,\frac{d}{\alpha}T^{-\frac{\beta-1}{\beta}}
\Big)\enspace.  
\end{equation}
This implies \eqref{eq1:lb} in the zone $\alpha\ge T^{-\frac{\beta+2}{2\beta}}=\alpha_0 $ since for such $\alpha$ we have
\[
\min\Big(\alpha, 
\,\frac{d}{\alpha}T^{-\frac{\beta-1}{\beta}}
\Big)
=
\min\Big(\max(\alpha, T^{-\frac{\beta+2}{2\beta}}), \frac{d}{\sqrt{T}}, \,\frac{d}{\alpha}T^{-\frac{\beta-1}{\beta}}\Big)\enspace.
\]
On the other hand, 
\(\min\big(\alpha_0, 
\,\frac{d}{\alpha_0}T^{-\frac{\beta-1}{\beta}}
\big)
=
\min\big(T^{-\frac{\beta+2}{2\beta}}, d/\sqrt{T} \big)
\).
The same lower bound holds for $0<\alpha<\alpha_0$  by the nestedness argument that we used to prove \eqref{eq2:lb} in the zone $0<\alpha<\alpha_0$. Thus, \eqref{eq1:lb} follows.
\end{proof}

%\begin{proof}[Proof of Theorem~\ref{thm:lower_grad}]
%Clearly, $\mathcal{F}'_{1, \beta} \subset \mathcal{F}_{\beta}$, where $\mathcal{F}'_{1, \beta}$ is the class of functions considered in Theorem \ref{lb} with $\alpha=1$. As the functions in $\mathcal{F}'_{1, \beta}$ are 1-strongly convex, for $t\geq 1$, $\norm{\nabla f(\bz_t)}^2\geq 2(f(\bz_t) -f^{\star}).$ 
%\begin{align*}   \norm{\nabla f(\bz_t)}^2\geq 2(f(\bz_t) -f^{\star}). \end{align*}
%Using Theorem \ref{lb}, for $T\geq 1$, we have
%\begin{align}\label{lbpook1}
%    \sup_{f \in \mathcal{F}'_{1,\beta}}\Exp[f(\bz_t)-f^{\star}]\geq {C}{'}dT^{-\frac{\beta-1}{\beta}},
%\end{align}
%where ${C}{'}>0$ does not depend on $d, T$, and $\beta$.
%For any random variable $S$ with values in $\{1,\dots,T\}$ independent of other sources of randomness, we have
%$\Exp[\norm{\nabla f(\bz_S)}^2]=\sum_{t=1}^{T}p_{t}\Exp[\norm{\nabla f(\bz_t)}^2],$
%where $(p_1,\dots,p_T)$ is a probability vector: $p_t\geq 0$, and $\sum_{t=1}^{T}p_t=1$. Thus,
%\begin{align}\label{lbpook2}
 %  \Exp[\norm{\nabla f(\bz_S)}^2]&\geq 2\sum_{t=1}^{T}p_t\Exp[f(\bz_t)-f^{\star}]
%   \geq
%   2\Exp\left[f\left(\sum_{t=1}^{T}p_t\bz_t\right)-f^{\star}\right].
%\end{align}
%Note that $\tilde{\bz}_T=\sum_{t=1}^{T}p_t\bz_t$ is an estimator depending only on the past, that is in the same class of estimators as $\bz_T$. Therefore, the bound (\ref{lbpook1}) holds for $\tilde{\bz}_T$ as well. Combining (\ref{lbpook1}) and (\ref{lbpook2}) gives (\ref{lb}).
%\end{proof}

\end{document}